\documentclass[a4paper,10pt]{article}

\usepackage{a4wide}
\usepackage{amsmath, amsthm, amsfonts, amssymb,mathtools}
\usepackage{mathrsfs}
\usepackage{dsfont}
\usepackage[T1]{fontenc}

\def\br{\mathbb R}

\def\bz{\mathbb Z}
\def\bc{\mathbb C}
\newcommand{\bra}{\langle}
\newcommand{\ket}{\rangle}

\def\e{\mathop{\mathrm{e}}\nolimits}

\newtheorem{Theorem}{Theorem}[section]
\newtheorem{Remark}[Theorem]{Remark}
\newtheorem{Lemma}[Theorem]{Lemma}
\newtheorem{Corollary}[Theorem]{Corollary}
\newtheorem{Proposition}[Theorem]{Proposition}
\newtheorem{Definition}[Theorem]{Definition}
\newtheorem{Assumption}[Theorem]{Assumption}

\makeatletter
 \@addtoreset{equation}{section}
 \makeatother
 
 \usepackage{hyperref}

\begin{document}


\title{On the singularities of the spectral shift function for some tight-binding models }

\author{M. Assal$^1$, O. Bourget$^{2}$, D. Sambou$^{3}$, A. Taarabt$^{2}$}

\date{\today}
\maketitle


\begin{quote}
\begin{itemize}
\item[$^1$] Departamento de Matem\'atica y Ciencia de la Computati\'on, \\
Universidad de Santiago de Chile, Las sophoras 173, Santiago, Chile.

\smallskip

\emph{E-mail: marouane.assal@usach.cl}

\medskip

\item[$^2$] Facultad de Matem\'aticas, Pontificia Universidad Cat\'olica de Chile,\\
Av. Vicu\~na Mackenna 4860, Santiago, Chile.

\smallskip

\emph{E-mails: bourget@uc.cl, amtaarabt@uc.cl}

\medskip

\item[$^3$] Institut Denis Poisson, Universit\'e d'Orl\'eans, UMR CNRS 7013, \\
45067 Orl\'eans cedex 2, France.

\smallskip

\emph{E-mail: diomba.sambou@univ-orleans.fr}
\end{itemize}
\end{quote}


\begin{abstract}
We consider perturbed discrete tight-binding models in $\ell^2(\bz_h,\mathcal{G})$ describing union of quantum particles with localized interactions, where $\bz_h$ is the 1D lattice $h\bz$, $h > 0$, and $\mathcal G$ is a separable Hilbert space. The perturbations play the role of self-adjoint relatively compact (matrix-valued) electric potentials with $\mathcal B(\mathcal G)$-valued coefficients decaying polynomially at infinity. We analyze the Spectral Shift Function (SSF) associated to the pair of the perturbed and the unperturbed operators. On the one hand, we show that the SSF is bounded near the spectral thresholds of the essential spectrum if $\dim(\mathcal G) < +\infty$. On the other hand, if $\dim(\mathcal G) = +\infty$, we show that it may have singularities at some thresholds points $\mu$ of the essential spectrum. In particular, new mechanisms allowing the SSF to have singularities at the thresholds are exhibited, based on the degeneracy of the spectrum of the unperturbed operator. Moreover, we give the main terms of the asymptotic behaviors of the SSF near $\mu$ described in terms of some explicit effective Berezin-Toeplitz type operators. These results are completed by Levinson type formulas and examples of eigenvalues asymptotics for power-like and exponential decay potentials.
\end{abstract}

\medskip

\textbf{AMS 2010 Mathematics Subject Classification:} 35J10, 81Q10, 35P20, 35P25, 47A10, 

47A11, 47A55, 47F05.

\smallskip

\textbf{Keywords:} Tight-binding models, discrete Schr\"odinger operators, spectral shift function, 

thresholds asymptotics, Levinson formula, discrete spectrum.

\pagebreak

\tableofcontents


\section{Introduction}\label{section_intro}
\setcounter{equation}{0}

\subsection{General setting and motivations} 

\begin{figure}
\begin{center}
\includegraphics[scale=0.3]{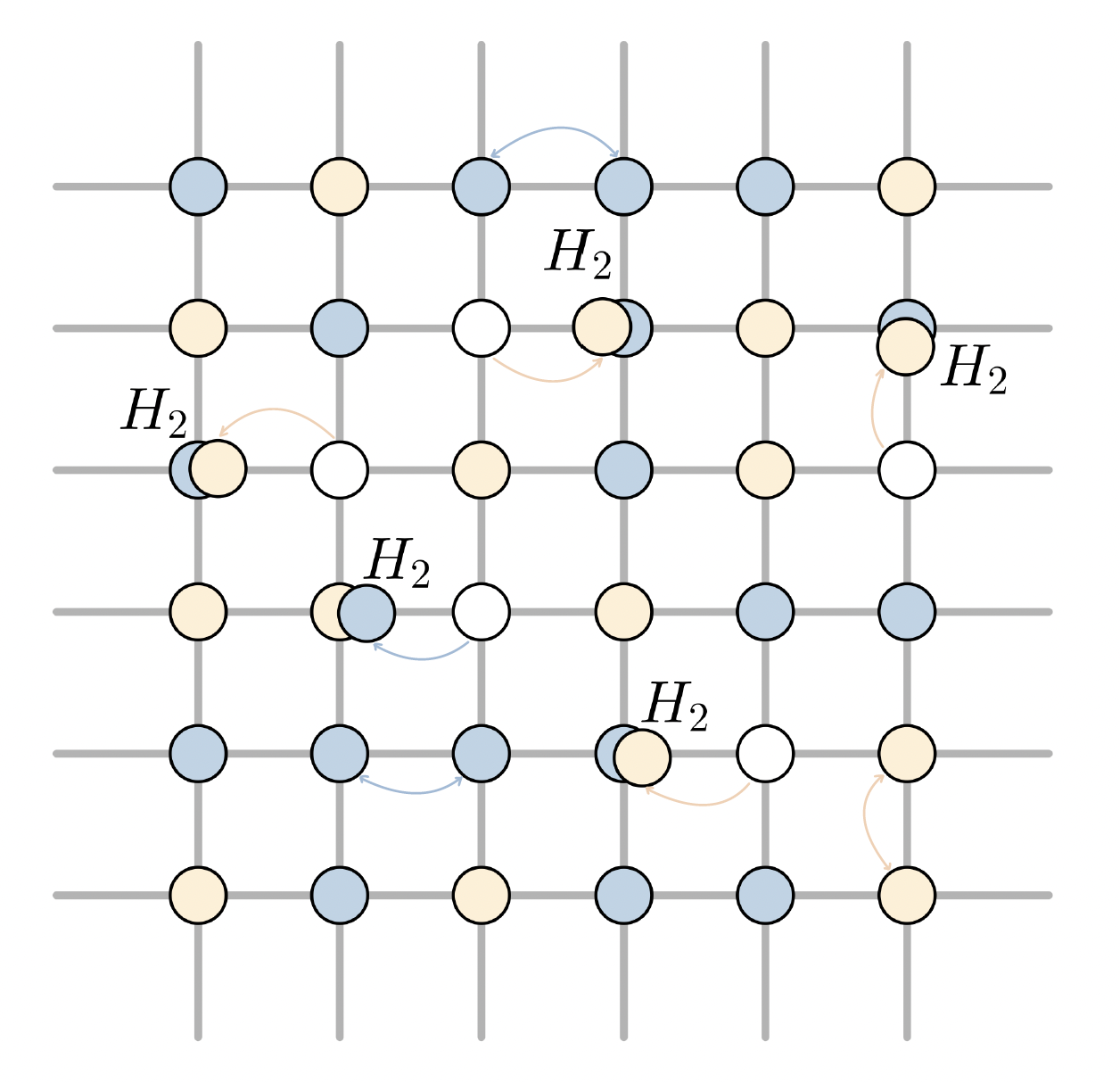}
\caption{Illustration of the Hubbard model.} \label{fig2}
\end{center}
\end{figure}

We consider operators of the form  $H_1\otimes I+I\otimes H_2$ that combine translational motion with internal configuration space. Such models are used in the physical literature to describe union of quantum particles with localized interactions where $H_1$ acts only on the first part of the system. 
The operator $H_2$ describes the dynamics within the second part of the system \cite{lewin}. 
Other tight-binding models \cite{DFF} have been studied in metal-insulator transitions
to capture the competing tendencies toward electron localization and mobility in materials (see also \cite{masc,mat}). For instance, the Hubbard model (see Figure \ref{fig2}) describes the behavior of interacting quantum particles in an atomic lattice where particle can hop between lattice sites and if two particles occupy the same site, they interact with each other through the operator $H_2$ \cite{hubbard1,hubbard2}. See also \cite{fef} and the references therein for recent works on tight-binding approximations for continuum magnetic two-dimensional crystalline structures.

This paper is devoted, in particular, to improve our understanding of the mechanisms involved in the distribution, the creation and the accumulation of bound states under relatively compact self-adjoint perturbations, in the vicinity of the thresholds points of the spectrum of some discrete tight-binding hamiltonians. 
For continuous models, these mechanisms have mainly been studied when a threshold coincides or is induced by an eigenvalue of infinite multiplicity (see \cite{rai,rawa1,rawa2,tda,ferrai,brra,brmi,sata} and references therein). 
Such phenomena are also related to long range perturbations at the threshold of the absolutely continuous component of the spectrum as is the case of the hydrogen atom model \cite{teschl}.
However, there are few results showing spectral accumulation phenomena near the thresholds points of the spectrum under relatively compact self-adjoint perturbations of discrete models (see the recent work \cite{mpr}). 


We consider the 1D lattice $$\bz_h := \{ hn : n \in \bz \},$$ with mesh size $h>0$. Let $\mathcal{G}$ be a separable Hilbert space and $\ell^2(\bz_h,\mathcal{G})$ be the Hilbert space endowed with the scalar product 
$
\langle \varphi,\phi \rangle := \sum_{n \in \bz} \langle \varphi(hn),\phi(hn) \rangle_\mathcal{G},
$
so that
$$
\ell^2(\bz_h,\mathcal{G}) = \Big\lbrace \varphi : {\bz_h} \to \mathcal{G} : \| \varphi \|^2 = \sum_{n \in \bz} \| \varphi(hn) \|_\mathcal{G}^2 < +\infty \Big\rbrace.
$$
For $\varphi \in \ell^2 ({\mathbb Z}_h,\mathcal{G})$, we define the finite-difference bounded operator
\begin{equation*}
(\partial \varphi) (hn) := \frac{1}{h^2} \big( \varphi(h(n+1)) - \varphi(hn) \big),
\end{equation*}
whose adjoint $\partial^\ast$ is given by
\begin{equation*}
(\partial^* \varphi) (hn) = \frac{1}{h^2} \big( \varphi(h(n-1)) - \varphi(hn) \big).
\end{equation*} 
We define the bounded self-adjoint Schr\"odinger operator 
\begin{equation}
\label{eq:g+1}
H_0 = -\partial - \partial^\ast\quad\text{on}\quad \ell^2(\bz_h,\mathcal{G}).
\end{equation}
Identifying $\ell^2(\bz_h,\mathcal{G})$ with $\ell^2(\bz_h)\otimes \mathcal{G}$, $H_0$ (see Section \ref{sec2} for more details) can be rewritten as
\begin{equation}
H_0 = -\Delta_h \otimes I_\mathcal{G},
\end{equation}
where $-\Delta_h$ is the 1D Schr\"odinger operator acting in $\ell^2(\bz_h) := \ell^2(\bz_h,\mathbb C)$ as
\begin{equation}
\label{eq:g+2}
(-\Delta_h \phi)(n) = \frac1{h^2} \big( 2 \phi(hn) - \phi(h(n+1)) - \phi(h(n-1)) \big).
\end{equation}
Then, the spectrum of $H_0$ is purely absolutely continuous and satisfies 
 \begin{equation}
\label{eq:g+4}
\sigma (H_0) = \sigma_{\textup{ac}} (H_0) = \sigma_{\textup{ess}} (H_0) = [0,\tfrac{4}{h^2}],
\end{equation}
where the points $\{0,\frac4{h^2}\}$ are the thresholds of this spectrum. 
The operator $H_0$ generalizes the Schr\"odinger operator $-\Delta_h$ on the discrete line $\bz_h$. When $\mathcal{G} \cong \mathbb C^m$, $m \ge 1$, it may be considered as the Hamiltonian of the system describing the behavior of a free particle moving in the strip $\mathbb Z_h \times \{1,\cdots,m\}$ (for more general settings see \cite{rabinovich}).
We define in $\ell^2(\bz_h)\otimes \mathcal{G}$ the operators 
\begin{equation}
\label{eq:mgten}
H := H_Q + V, \qquad H_Q := -\Delta_h \otimes I_\mathcal{G} + I_{\ell^2(\bz_h)} \otimes Q,
\end{equation}
where $V$ is a relatively compact self-adjoint electric potential and  $Q$ is a self-adjoint operator acting on $\mathcal{G}$. 
For $Q \neq 0$, we assume that its spectrum is given by a set of real eigenvalues (counted with multiplicity)
$$
\sigma(Q) = \{ \mu_{s} : s \in S \subseteq \mathbb Z_+ \},
$$ 
such that 
\begin{equation}
\label{decG}
\mathcal G = \bigoplus_{s=1}^d {\rm Ker}(Q - \mu_s), \quad d < +\infty.
\end{equation} 
Notice that if $\dim(\mathcal{G})<+\infty$, then \eqref{decG} holds trivially by the spectral theorem with $d$ equal to the number of distinct eigenvalues of $Q$. Otherwise, if $\dim(\mathcal{G})=+\infty$, then \eqref{decG} implies that $Q$ is 
unitarily diagonalizable on a orthonormal eigenbasis of $\mathcal{G}$, and there exists $1 \le s \le d$ such that $\dim {\rm Ker}(Q - \mu_s) = +\infty$. By Weyl's criterion, we have
\begin{equation}
\sigma_{\rm ess}(H) = \sigma(H_Q) = [0,\tfrac{4}{h^2}] + \sigma(Q) = 
\begin{cases} 
[0,\tfrac{4}{h^2}] \quad {\rm if} \quad Q = 0, \\
\bigcup_{s = 1}^d [\mu_{s},\tfrac{4}{h^2} + \mu_{s}] \quad {\rm if} \quad Q \neq 0,
\end{cases},
\end{equation}
so that 
\begin{equation}
\mathscr E_Q := 
\begin{cases} 
\{0,\frac{4}{h^2} \} \quad {\rm if} \quad Q = 0, \\
\{\mu_{s},\mu_{s}+\frac{4}{h^2} \}_{s=1}^d \quad {\rm if} \quad Q \neq 0,
\end{cases}
\end{equation}
plays the role of the spectral thresholds set of the spectrum $\sigma(H_Q)$. 
For $1 \le s \le d$, one denotes by $\pi_s$ the projection onto ${\rm Ker}(Q - \mu_{s})$. 

\begin{Definition}
If $Q=0$, we set $\mu_0 := 0$.
If $Q \neq 0$, $\mu_s \in \mathscr E_Q$ (resp. $\frac4{h^2}+\mu_s)$ is non-degenerate if $\mu_s \neq \frac4{h^2}+\mu_{s'}$ (resp. $\frac4{h^2}+\mu_s \neq \mu_{s'})$ for all $s \neq s'$. If not, it is said to be degenerate.
\end{Definition}

 We consider self-adjoint matrix-valued electric perturbations $V$ with coefficients decaying polynomially at infinity,
 and we are interested in the spectral properties of the operators 
$$
H^\pm := H_Q \pm V, 
$$ 
where  $V \in \mathcal B(\ell^2(\bz_h,\mathcal{G}))$  is a positive matrix-valued electric potential such that
\begin{equation}
\label{eq:pert}
V = \{ V_h(n,m) \}_{(n,m) \in \bz^2}, \quad V_h(n,m) \in \mathcal{B}(\mathcal{G}).
\end{equation}
 Here, $\mathcal{B}(\mathcal{G})$ denotes the set of bounded linear operators in $\mathcal G$. This electric potential $V$ can be interpreted as a summation kernel operator whose kernel is given by the operator-valued function $$(n,m) \in \bz^2 \mapsto V_h(n,m) \in \mathcal{B}(\mathcal{G}).$$ Namely, for any $\varphi \in \ell^2(\bz_h,\mathcal{G})$, one has
 \begin{equation}\label{eq:g+9}
(V \varphi)(hn) = \sum_{m \in \bz} V_h(n,m) \varphi(hm), \quad n \in \bz.
\end{equation}
We also assume that 
$\Vert V_h(n,m)\Vert_{\mathcal B(\mathcal{G})}$ 
decays more rapidly than $\Vert (n,m) \Vert^{-3}$ as $\Vert (n,m) \Vert \to +\infty$ (see Assumption \ref{eq:hyppert0}, (3) of Remark \ref{r:spd+} and \eqref{hyp:decV} for more details). 
So, we show that
\begin{equation}\label{TOP}
H^\pm - H_Q \in \mathfrak{S}_1,
\end{equation}
where $\mathfrak S_1$ denotes the trace class operators. Then, there  exists (see \cite{krei} or e.g. \cite[Theorem 8.3.3]{yaf}) a unique function $\xi(\cdot;H^{\pm},H_Q) \in {\rm L}^1(\br)$ such that the \textit{Lifshits-Krein trace formula}
\begin{equation}
\label{eq:fctdec}
 {\rm Tr}(f(H^\pm)-f(H_Q)) = \int_\br \xi (\lambda;H^\pm,H_Q) f'(\lambda) d\lambda,
\end{equation}
holds for every $f \in C_0^\infty(\br)$. The function $\xi(\cdot, H^\pm, H_Q)$ is called the \textit{Spectral Shift Function} (SSF) for the pair $(H^\pm,H_Q)$. It can be related to the number of eigenvalues of the operators $H^\pm$ in $\mathbb R \smallsetminus \sigma_{\rm ess}(H^\pm)$ (see formula \eqref{eq:relssfvp}). It is also related to the scattering matrix $S(\lambda;H^\pm,H_Q)$ for the pair $(H^\pm,H_Q)$ by the Birman-Krein formula 
\begin{equation}
\label{eq:scatdet}
\det S(\lambda;H^\pm,\mathcal H_Q) = e^{2i\pi\xi(\lambda;H^\pm,H_Q)}, \quad a.e. \: \lambda \in \sigma_{\rm ess}(H^\pm).
\end{equation}

We purpose to analyze in particular the distribution of the bound states near $\sigma_{\rm ess}(H_Q \pm V)$, adapting techniques present in the literature (see \cite{rai,rawa1,rawa2,tda,ferrai}). In particular, we improve and we generalize previous results (see \cite{bst,bst+}) on a class of discrete Laplace type operators on the 1D lattice and on strips, established under self-adjoint exponential decay matrix-valued perturbations at infinity. 
The techniques developed and used in \cite{bst,bst+} are based on resonances theory and complex scaling arguments.
However, this strategy is not adapted to analyze matrix-valued perturbations that decay polynomially at infinity as in our case here. For such perturbations, the spectral analysis can be performed using the SSF \cite{lifshits,krei}, a useful notion for spectral analysis and scattering theory of quantum systems. 
Technically, the formation of cluster of eigenvalues is somehow encoded in the behavior of the SSF. In contrast to the scattering matrix, the spectral shift function is meaningful both on the continuous and discrete spectra.

\subsection{Description of the main results}

The main results of this article concern the asymptotic behavior of the SSF $\xi(\lambda;H^\pm,H_Q)$ as $\lambda \to \lambda_0 \in \mathscr E_Q$, for matrix-valued electric potentials $\pm V \ge 0$. 

We will first identify $\xi(\cdot;H^{\pm},H_Q)$ with a representative of its equivalence class described explicitly in Section \ref{sec4}, assuming that the electric matrix-valued potential $V$ has a definite sign. Then, we show the boundedness of $\xi(\cdot;H^\pm,H_Q)$  on compact subsets of $\br \smallsetminus \mathscr E_Q$ independently on the dimension of $\mathcal{G}$ (see Theorem \ref{theo0}). In Theorem \ref{theo1}, we establish the asymptotic behavior of $\xi(\lambda;H^\pm,H_Q)$ as $\lambda \nearrow \mu_s$ and as $\lambda \searrow \frac{4}{h^2} + \mu_s$. In Theorems \ref{theo2}, \ref{theo2+}, we determine the asymptotic behavior of $\xi(\lambda;H^\pm,H_Q)$ as $\lambda \searrow \mu_s$ and as $\lambda \nearrow \frac{4}{h^2} + \mu_s$. Several consequences can be deduced from 
these results.

In the finite-dimensional case $\dim(\mathcal G) < +\infty$, Theorem \ref{theo1:cor1} implies that the SSF $\xi(\cdot;H^\pm,H_Q)$ is bounded in $\br \smallsetminus \mathscr E_Q$,
which improves Theorem \ref{theo0}. In particular, if $\mathcal G = \mathbb C$, then the operator $H_0 = -\Delta_h$ is the 1D discrete Laplacian on the lattice $\bz_h$ and  Corollary \ref{cor:specdisfini} shows the finiteness of the discrete spectrum of $-\Delta_h \pm V$ for polynomial decay perturbations at infinity. This extends results of \cite{bst,bst+} where the finiteness of the discrete spectrum of $-\Delta_1+V$ has been proved for self-adjoint exponential decay perturbations at infinity. 

Otherwise, if $\dim(\mathcal G) = +\infty$, we prove that $\xi(\cdot;H^\pm,H_Q)$ 
may have singularities at the spectral thresholds $\mu_s$ and $\frac4{h^2}+\mu_s$, $0 \le s \le d$, with $\dim {\rm Ker}(Q - \mu_{s}) = +\infty$ for $s \ge 1$, under generic assumptions on $V$ (see Theorems \ref{cor3:mr}, \ref{theo:fin} and Corollaries \ref{cor4:mr}, \ref{cor:lev-}, \ref{cor:fin1}). More precisely, for $V > 0$, we have for such thresholds
\smallskip
\begin{center}
$\begin{cases}
\xi(\lambda;H^+,H_Q) = \mathcal O(1) \quad as \quad \lambda \nearrow \mu_s, \\ 
\xi(\lambda;H^+,H_Q) \to +\infty \quad as \quad \lambda \searrow \mu_s,
\end{cases} {\rm while} \: \:
\xi(\lambda;H^+,H_Q) \to +\infty \quad as \quad \lambda \to \frac4{h^2} + \mu_s,
$
\end{center}
\smallskip
and
\smallskip
\begin{center}
$\begin{cases}
\xi(\lambda;H^-,H_Q) = \mathcal O(1) \quad as \quad \lambda \searrow \frac4{h^2} + \mu_s, \\
\xi(\lambda;H^-,H_Q) \to -\infty \quad as \quad \lambda \nearrow \frac4{h^2}+\mu_s,
\end{cases} {\rm while} \: \:  
\xi(\lambda;H^-,H_Q) \to -\infty \quad as \quad \lambda \to \mu_s.
$
\end{center}
Actually, the singularities of the SSF at the spectral thresholds are described in terms of some explicit effective "Berezin-Toeplitz" type operators (see \eqref{eq:raj+} for a precise definition). Hence and  under suitable condition, we give the main terms of the asymptotic expansions of $\xi(\lambda;H^\pm,H_Q)$ as $\lambda \to z_0 \in \{\mu_s,\frac4{h^2}+\mu_s\}$ (see Corollary \ref{cor4:mr} and Theorem \ref{theo:fin} for the general case, and  Corollaries \ref{cor:lev-}, \ref{cor:fin1} for power-like and exponential decay perturbations). In particular, if $V > 0$, then 
\begin{equation*}
\lim_{\lambda \searrow 0} \frac{\xi(\mu_s+\lambda;H^-,H_Q)}{\xi(\mu_s-\lambda;H^-,H_Q)}
\quad\text{and}\quad
\lim_{\lambda \searrow 0} \frac{\xi(\frac4{h^2}+\mu_s-\lambda;H^+,H_Q)}{\xi(\frac4{h^2}+\mu_s+\lambda;H^+,H_Q)},
\end{equation*}
exist and are equal to positive constants depending on the decay rate of $V$ at infinity (see Theorem \ref{theo:levin1} and Corollary \ref{cor:levin1}). This can be interpreted as generalized Levinson formulae (see the original work \cite{lev} or the survey article \cite{rob}).

\subsection{Comments on the literature} 

Our results extend to a class of discrete tight-binding models those established in \cite{rai,tda,ferrai} for continuous models. More precisely, in \cite{rai,tda} the asymptotic behavior of the SSF has been considered near the low ground energy and near $\pm m$ for  2D Pauli and 3D Dirac operators with non-constant magnetic fields, respectively. In \cite{ferrai} the asymptotic behavior of the SSF has been considered near the Landau levels for 3D Schr\"odinger operators with constant magnetic fields. Similar results can be also found in \cite{brra,brmi}. However, in the discrete case, there are few results concerning the asymptotics expansions of the SSF at spectral thresholds. 
To our best knowledge, the most recent work in this direction showing spectral accumulation phenomena near the essential spectrum for self-adjoint perturbations, seems to be \cite{mpr}. 
It is important to highlight that in the papers mentioned above, the singularities of the SSF near the spectral thresholds are induced by infinitely degenerated eigenvalues. This is similar to our situation in the case $Q \neq 0$, where we show that the SSF may have singularities at the spectral thresholds $\mu_s$ and $\frac4{h^2}+\mu_s$ when ${\rm Rank} \, \pi_s = +\infty$.
This is rather different compared to our particular case $Q = 0$ where the singularities of the SSF near the thresholds $\{0,\frac4{h^2}\}$ are produced by a highly degenerated absolutely continuous component. Note that the singularities of the SSF are probably due to an accumulation of resonances near the spectral thresholds. However, this aspect of the problem will not be addressed here and will be consider in a further work. 
In the finite-dimensional case $\dim(\mathcal G) < +\infty$, this issue is considered in \cite{abms} for (non)-selfadjoint exponential decay perturbations of $H_Q$.

\smallskip

The article is organized as follows. In Section \ref{sec2}, we perform the spectral analysis of the operator $H_0$ and introduce some standard tools needed so far. In Section \ref{sec3}, we state and discuss our main assumptions concerning the perturbed operators $H^\pm$. In Section \ref{sec4}, we recall some abstract results due to A. Pushnitski on the representation of the spectral shift function for a pair of self-adjoint operators. The Section \ref{sec5} is devoted to the formulation of our main results, some corollaries of them, as well as examples of explicit eigenvalues asymptotics. In Section \ref{sec6}, we compute a suitable decomposition of the potential $V$ satisfying a main assumption given in Section \ref{sec3}. The Section \ref{sec7} contains auxiliary material such as extensions of the convolution kernel of $-\Delta_h$, and estimates of appropriate weighted resolvents. In Section \ref{sec8}, we prove Theorems \ref{theo1} and \ref{theo2} while in Section \ref{sec9} we prove the asymptotics identities \eqref{as:lifs-kre0} and \eqref{as:lifs-kre4}.

\section{Spectral properties of the hamiltonian $H_0$}\label{sec2}

Let $\Lambda \subseteq \mathbb Z_+$ and consider an orthonormal basis $(\e_j)_{j \in \Lambda}$ of ${\mathcal{G}}$ (formed by orthonormal basis of ${\rm Ker}(Q - \mu_{s})$, $1 \le s \le d$, if $Q \neq 0$, in accordance with \ref{decG}). Of course $\# \Lambda = \dim(\mathcal{G})$ if $\dim(\mathcal{G}) < +\infty$ and we take $\Lambda = \mathbb Z_+$ if $\dim(\mathcal{G}) = +\infty$. 
Let $(\delta_{n} )_{n\in {\mathbb Z}}$ be the canonical orthonormal basis of  $\ell^2 ({\mathbb Z}_h),$ where $\delta_{n} (kh) = \delta_{nk}$ for $k\in {\mathbb Z}$. Then, it is useful to identify the spaces $\ell^2(\bz_h,\mathcal{G})$ and $\ell^2(\bz_h)\otimes \mathcal{G}$ so that
$
\ell^2(\bz_h,\mathcal{G}) \cong \ell^2(\bz_h)\otimes \mathcal{G}
$
and it follows that $(\delta_n \otimes \e_j)_{(n,j)\in {\mathbb Z} \times \Lambda}$ is an orthonormal basis of $\ell^2(\bz_h)\otimes \mathcal{G}$. For $j \in \Lambda$, one defines 
$
{\mathcal{G}}_j = {\rm span} \{ x \otimes \e_j : x\in \ell^2(\bz_h)\},
$
together with its corresponding orthogonal projection $Q_j := I_{ \ell^2(\bz_h)}  \otimes |\e_j \ket \bra \e_j |$, so that
$$
\ell^2(\bz_h,\mathcal{G}) \cong \ell^2(\bz_h) \otimes  \mathcal{G} = \bigoplus_{j \in \Lambda} \mathcal{G}_j.
$$
Hence, we notice that for every $j \in \Lambda$, $\mathcal{G}_j$ is $H_0$-invariant and thus
$$
H_0 = \bigoplus_{j \in \Lambda} Q_j H_0 Q_j = \bigoplus_{j \in \Lambda} -\Delta_h \otimes |\e_j \ket \bra \e_j | = -\Delta_h \otimes I_\mathcal{G},
$$
where $-\Delta_h$ is the 1D Schr\"odinger operator given by \eqref{eq:g+2}.

Let $\tau>0$ be such that $h\tau =2\pi$ and
$$
{\mathbb T} = {\mathbb R}/\tau{\mathbb Z} \sim [-\tfrac{\tau}{2},\tfrac{\tau}{2}].
$$
In view of the bijection between $\ell^2 ({\mathbb Z}_h)$ and ${\rm L}^2 ({\mathbb T}) := {\rm L}^2 ({\mathbb T},\bc)$, one defines the discrete Fourier transform ${\mathscr F} : \ell^2 ({\mathbb Z}_h) \rightarrow {\rm L}^2 ({\mathbb T})$ by
\begin{equation}
\label{eq:Tfour}
({\mathscr F} \phi) (\theta) := \sum_{n\in {\mathbb Z}} \phi(hn)e^{-ihn\theta}, \quad
\phi(hn) =\frac1\tau \int_{{\mathbb T}} ({\mathscr F} \phi) (\theta) e^{ihn\theta} d\theta.
\end{equation}
Since the operator ${\mathscr F}$ is unitary, then by using the partial transform $\mathscr F \otimes I_\mathcal{G}$ acting in $\ell^2(\bz_h) \otimes {\mathcal{G}}$, one can show that $H_0$ is unitarily equivalent to the operator $\widehat {-\Delta_h} \otimes I_\mathcal{G}$ acting in ${\rm L}^2({\mathbb T}, {\mathcal{G}}) \cong {\rm L}^2({\mathbb T}) \otimes {\mathcal{G}}$, where $\widehat {-\Delta_h}$ is the multiplication operator in ${\rm L}^2({\mathbb T})$ by the function $f$ defined by
\begin{equation}
\label{eq:fonctf}
f(\theta) := \frac{2 - 2\cos(h\theta)}{h^2} = \frac4{h^2} \sin^2 \Big( \frac{h\theta}{2} \Big), \quad \theta \in {\mathbb T}.
\end{equation} 
Therefore, $[0,\tfrac{4}{h^2}] = \sigma (-\Delta_h) = \sigma (H_0)$ and the spectrum of the operators $-\Delta_h$ and $H_0$ are purely absolutely continuous so that \eqref{eq:g+4} holds.

\section{The electric potentials} \label{sec3}

Recall that the potential $V$ acting in $\ell^2(\bz_h,\mathcal{G})$ is a bounded matrix-valued $V = \{ V_h(n,m) \}_{(n,m) \in \bz^2}$ with coefficients $V_h(n,m) \in \mathcal{B}(\mathcal{G})$. In the basis $(\e_j)_{j \in \Lambda}$ of $\mathcal{G}$, for each $(n,m) \in \bz^2$, the operator $V_h(n,m)$ has the matrix representation
\begin{equation}
\label{eq:g+9+}
V_h(n,m) = \{ v_{jk}^h(n,m) \}_{j,k \in \Lambda}, \qquad v_{jk}^h(n,m) := \langle \e_j,V_h(n,m) \e_k \rangle_\mathcal{G}.
\end{equation}
Hence, one has 
\begin{equation}
\label{eq:g+9,0}
V_h(n,m) = \sum_{(j,k) \in \Lambda^2} v_{jk}^h(n,m) | \e_j \rangle \langle \e_k |.
\end{equation} 
The operator $V_h(n,m)$ viewed as a matrix $\{ v_{jk}^h(n,m) \}_{(j,k) \in \Lambda^2}$ belongs to $\mathcal{M}_r(\bc)$ if $r = \dim(\mathcal{G})<+\infty$. 
So, in $\ell^2(\bz_h) \otimes  \mathcal{G}$, $V$ has a canonical representation given by
$$
V = \sum_{(n,m)\in\bz^2} | \delta_n \rangle \langle \delta_m | \otimes V_h(n,m) = \sum_{(n,m)\in\bz^2} \sum_{(j,k) \in \Lambda^2} | \delta_n \rangle \langle \delta_m | \otimes v_{jk}^h(n,m) | \e_j \rangle \langle \e_k |.
$$

In the sequel, for $y = (y_1,\ldots,y_d) \in \br^r$, one sets $\bra y \ket := (1+\vert y \vert^2)^{1/2}$. Bearing in mind \eqref{eq:pert} and \eqref{eq:g+9+}, we introduce the following polynomial decay assumption on $V$.
\begin{Assumption}
\label{eq:hyppert0}
$V = \{ V_h(n,m) \}_{(n,m) \in \bz^2}$ is of definite sign ($V \ge 0)$ such that 
\begin{equation}
\label{hyp:decV0}
\vert v_{jk}^h(n,m) \vert \leq {\rm Const} \cdot G_1(j,k) \bra n \ket^{-\nu_1} \bra m \ket^{-\nu_2}, \quad (n,m) \in \bz^2, 
\end{equation}
for some $\nu_1, \nu_2 > 1$, where $0 \le G_1$ defined in $\Lambda^2$ satisfies
\begin{center}
$\begin{cases}
G_1 \in {\rm L}^\infty(\Lambda^2) \quad if \quad \dim(\mathcal{G}) < +\infty, \\
G_1(j,k) \le {\rm Const.} \bra j \ket^{-\beta_1} \bra k \ket^{-\beta_2} \quad if \quad \dim(\mathcal{G}) = +\infty,
\end{cases}$
\end{center}
$(j,k) \in \Lambda^2$, for some constants $\beta_1, \beta_2 > 1$.
\end{Assumption}

Let us make some comments on Assumption \ref{eq:hyppert0}.
\begin{itemize}
\item If $\dim(\mathcal{G}) < + \infty$, then typical examples of potentials satisfying \eqref{hyp:decV0} are $V$ such that
\begin{equation}
\label{dec1}
\vert v_{jk}^h(n,m) \vert \le {\rm Const.} \bra (hn,hm) \ket^{-\nu}, 
\end{equation}
$(n,m) \in \bz^2$, $(j,k) \in \Lambda^2$, $\nu > 2$. Indeed \eqref{dec1} implies that for every $(n,m) \in \bz^2$, 
$$
\vert v_{jk}^h(n,m) \vert \leq {\rm Const.} \bra hn \ket^{-\nu/2} \bra hm \ket^{-\nu/2} \leq {\rm Const.} \frac{\bra n \ket^{-\nu/2} \bra m \ket^{-\nu/2}}{(\min(1,h^2))^{\nu}}.
$$
\item If $\dim(\mathcal{G}) = +\infty$, then \eqref{hyp:decV0} holds for instance if for $(j,k) \in \Lambda^2$, 
\begin{equation}
\label{eqaj+}
\vert v_{jk}^h(n,m) \vert \le {\rm Const} \cdot \langle (j,k) \rangle^{-\beta} \bra (n,m) \ket^{-\nu}, 
\end{equation}
$(n,m) \in \mathbb{Z}^2$, $\beta>2$, $\nu > 2$.
For example, \eqref{eqaj+} is satisfied if $\beta>2$, $\nu > 2$ and 
$$
V_h(n,m) = \bra (hn,hm) \ket^{-\nu} \sum_{(j,k)\in\Lambda^2} \langle (j,k) \rangle^{-\beta} |\e_j \rangle \langle \e_k|.
$$
\end{itemize}

Now, one sets
\begin{equation}
\label{exp:nu0}
\nu_0 := \min(\nu_1,\nu_2) > 1 \quad {\rm and} \quad \beta_0 := \min(\beta_1,\beta_2) > 1.
\end{equation}
Consider the function 
\begin{equation}
\label{def:psi}
\psi := \bra (\cdot)h^{-1} \ket^{-\nu_0/2} : hn \in \bz_h \mapsto \bra n \ket^{-\nu_0/2} \in \br_+^\ast,
\end{equation}
and define in $\ell^2(\bz_h)$ the multiplication operator $M_\psi$ by the function $\psi$.
Namely, for $\phi \in \ell^2(\bz_h)$, $(M_\psi \phi)(hn) =  
\bra n \ket^{-\nu_0/2} \phi(hn)$ or
\begin{equation}
\label{def:opn}
M_\psi := \sum_{n \in \bz} \bra n \ket ^{-\nu_0/2} | \delta_n \rangle \langle \delta_n |.
\end{equation}
Similarly, one defines $p$ the operator acting in $\mathcal{G}$ by 
\begin{equation}
\label{def:opp}
p := \sum_{j \in \Lambda} \bra j \ket ^{-\beta_0/2} | \e_j \rangle \langle \e_j |.
\end{equation}

\begin{Remark}
\label{rem:perV}
The matrix representations of $M_\psi$ and $p$ are diagonal. Moreover, $M_\psi$ and $p$ belong to $\mathfrak S_2$ the Hilbert-Schmidt class since $\sum_{n \in \bz} \bra n \ket ^{-\nu_0} < \infty$ and $\sum_{j \in \Lambda} \bra j \ket ^{-\beta_0} < \infty$. Hence, they belong to $\mathfrak S_\infty$ the class of compact linear operators. Of course, if $\dim(\mathcal{G}) < \infty$, then $p \in \mathfrak S_2$. 
\end{Remark}

In Lemma \ref{l:sdecV}, one proves the following decomposition of $V$. More precisely, if $V$ satisfies the polynomial decay Assumption \ref{eq:hyppert0}, we show that there exists $\mathscr V \in \mathcal B (\ell^2(\bz_h,\mathcal{G}))$, $\mathscr V \ge 0$, such that 
\begin{equation}
\label{V}
V = (M_\psi \otimes p) {\mathscr V} (M_\psi \otimes p).
\end{equation}
Moreover, $V = {\mathscr M}^\ast {\mathscr M}$ is trace class with ${\mathscr M} = \mathscr V^{1/2} (M_\psi \otimes p)$ and 
$
\Vert V \Vert_{\mathfrak S_1} \le \Vert {\mathscr M} \Vert_{\mathfrak S_2}^2.
$

\begin{Remark}
\label{rem0}
\begin{enumerate}
\item Under Assumption \ref{eq:hyppert0}, the factorization \eqref{V} of $V$ is not unique and other choices can be more suitable. For instance, if $\dim(\mathcal{G}) < +\infty$, one can deal in our analysis with the decomposition \eqref{rep:V+} introduced in the  proof of Lemma \ref{l:sdecV}.
\item If $\dim(\mathcal{G}) = +\infty$, suppose moreover that there exists $n_0 > 0$ (fixed) such that for all $(n,m) \in \bz^2$, $v_{jk}^h(n,m) = 0$ for each $j >n_0$ and $k > n_0$. That is, $V_h(n,m)$ is of the form
\begin{equation}
\label{eq:decdiminfranf}
  V_h(n,m) = \left(\begin{array}{ c | c }
    M_{n_0}^h(n,m) & {\bf 0} \\
    \hline
    {\bf 0} & {\bf 0}
  \end{array}\right).
\end{equation}
Then, the same argument used in proof of Lemma \ref{l:sdecV} part b), allows to replace the operator $p$ in \eqref{V}  by the finite-rank operator 
\begin{equation}
\widetilde p := \sum_{j=0}^{n_0} \bra j \ket ^{-\beta_0/2} | \e_j \rangle \langle \e_j |.
\end{equation}
\end{enumerate}
\end{Remark}

So, bearing in mind \eqref{V} and Remark \ref{rem0}, we will consider the perturbed operators
$H^\pm$
with self-adjoint perturbations $V$ satisfying the next assumption which extends Assumption \ref{eq:hyppert0}. 

\begin{Assumption}
\label{eq:hyppert}
$V= (M_\psi \otimes K^\ast) {\mathscr V} (M_\psi \otimes K)$, $\nu_0 > 1$, where
$0 \le \mathscr V \in \mathcal B (\ell^2(\bz_h,\mathcal{G}))$ and
$K$ acting in $\mathcal{G}$ satisfies $K \in {\mathfrak S}_2 (\mathcal{G})$.
\end{Assumption}

Potentials $V$ satisfying Assumption \ref{eq:hyppert} belong to the trace class $\mathfrak S_1$ with
\begin{equation}
\label{est1V}
\Vert V \Vert_{\mathfrak S_1} \le \Vert {\mathscr M} \Vert_{\mathfrak S_2}^2 \le \Vert \mathscr V \Vert \Vert M_\psi \otimes K \Vert_{\mathfrak S_2}^2.
\end{equation}

In order to fix ideas, let us point out some important remarks on Assumption \ref{eq:hyppert}.

\begin{Remark}
\label{r:spd+}
\begin{enumerate}
\item If $\dim(\mathcal{G}) = +\infty$, then Assumption \ref{eq:hyppert0} implies Assumption \ref{eq:hyppert} with $K = K^\ast = p$, according to \eqref{V}. Since $p$ is of infinite rank, then  Assumption \ref{eq:hyppert} includes the class of finite-rank operators $K$ (as $\widetilde p$) in $\mathcal{G}$. 
\item Under Assumption \ref{eq:hyppert}, one has
\begin{equation}
\label{eq:Vmm} 
V = {\mathscr M}^\ast {\mathscr M}, \quad \quad {\mathscr M} := {\mathscr V}^{1/2} (M_\psi \otimes K) \in \mathfrak S_2.
\end{equation}
\item Our main results will be formulated under a more restrictive assumption, namely with $\nu_0 > 3$.
\end{enumerate}
\end{Remark}



\section{Representation of the spectral shift function} \label{sec4}

In this section, one recalls some abstract results due to A. Pushnitski on the representation of the spectral shift function for a pair of self-adjoint operators.
Let us define the sandwiched resolvent 
 \begin{equation}
 \label{def:Tpond}
 T(z) := {\mathscr M} (H_Q-z)^{-1} {\mathscr M}^\ast, \quad z \in \bc \smallsetminus \sigma(H_Q),
 \end{equation}
 where ${\mathscr M}$ is given by \eqref{eq:Vmm} and 
 \begin{equation}
 (H_Q-z)^{-1} = \sum_{s=0}^d (-\Delta_h + \mu_{s} -z)^{-1} \otimes \pi_s,
 \end{equation}
 where $\mu_0 = 0$ together with the following {\bf conventions}.
 \begin{itemize}
 \item For $Q=0$, we set $\pi_0 = I_{\mathcal G}$ and $\pi_s = 0$ for $s \ge 1$, so that $(H_0-z)^{-1} = (-\Delta_h - z)^{-1} \otimes I_{\mathcal G}$.
 \item For $Q \neq 0$, one sets $\pi_0 = 0$ so that $(H_Q-z)^{-1} = \sum_{s=1}^d (-\Delta_h + \mu_{s} -z)^{-1} \otimes \pi_s$.
 \end{itemize}
Denote by 
 \begin{equation}
 A(z) := {\rm Re} \, T(z) \quad {\rm and} \quad B(z) := {\rm Im} \, T(z),
 \end{equation}
 the real and the imaginary parts of the operator $T(z)$ respectively. Then, under \eqref{eq:Vmm}, it is well known that for a.e. $\lambda \in \br$, the limit
 \begin{equation}
 T(\lambda+i0) := \lim_{\varepsilon \searrow 0} T(\lambda + i\varepsilon),
 \end{equation}
exists in the $\mathfrak S_2$-norm (and even in the $\mathfrak S_p$-norm for any $p>1$). Moreover $0\le B(\lambda+i0) \in \mathfrak S_1$. See \cite{yaf,bir} and \cite{nab} for the case $p> 1$. 
Let $\mathcal T= \mathcal T^\ast \in \mathfrak S_\infty(\mathcal G)$.
Define
\begin{equation}
\label{eq:fcompt}
\mathscr N_\pm(r,\mathcal T) := {\rm Rank} \: \mathds 1_{(r,\infty)}(\pm \mathcal T), \quad r > 0,
\end{equation}
the counting functions of the positive eigenvalues of $\pm \mathcal T$. Then, by \cite[Theorem 1.1]{pus} we have:

\begin{Theorem}
\label{theorepPush}
Let Assumption \ref{eq:hyppert} holds. Then, for a.e. $\lambda \in \br$, the SSF $\xi(\cdot;H^\pm,H_Q)$ admits the representation via the converging integral
\begin{equation}
\label{rep:Push}
\xi(\lambda;H^\pm,H_Q) = \pm \int_{\br} \mathscr N_\mp(1,A(\lambda+i0)+tB(\lambda+i0))  \frac{dt}{\pi(1+t^2)}.
\end{equation}
\end{Theorem}

For further use, let us recall the following estimates, useful in the study of the convergence of the r.h.s. of \eqref{rep:Push}.

\begin{Lemma}[Lemma 2.1 of \cite{pus}]
\label{lem:pushCV}
Let $T_1 = T_1^\ast \in \mathfrak S_\infty$ and $T_2 = T_2^\ast \in \mathfrak S_1$ acting in the same Hilbert space. Then, for any $x_1$, $x_2 > 0$, one has
$$
\frac1\pi \int_{\br} \mathscr N_\pm(x_1 + x_2,T_1+tT_2)  \frac{dt}{1+t^2} \le \mathscr N_\pm(x_1,T_1) + \frac{1}{\pi x_2} \Vert T_2 \Vert_{\mathfrak S_1}.
$$
\end{Lemma}
 In Corollary \ref{cor:exTs1}, one establishes that $T(\lambda +i0)$ belongs to $\mathfrak S_1$ for every $\lambda \in \br \smallsetminus \mathscr E_Q$. It follows from Lemma \ref{lem:pushCV} that the r.h.s. of \eqref{rep:Push} is well-defined for each $\lambda \in \br \smallsetminus \mathscr E_Q$. So, one can consider the function $\widetilde \xi(\cdot;H^\pm,H_Q)$ defined in $\br \smallsetminus \mathscr E_Q$ by
 \begin{equation}
 \label{eq:foncdecgen}
 \lambda \in \br \smallsetminus \mathscr E_Q \mapsto \widetilde \xi(\lambda;H^\pm,H_Q) = \pm \frac1\pi \int_{\br} \mathscr N_\mp(1,A(\lambda+i0)+tB(\lambda+i0))  \frac{dt}{1+t^2}.
 \end{equation}
 By Theorem \ref{theorepPush}, 
 $
 \widetilde \xi(\lambda;H^\pm,H_Q) = \xi(\lambda;H^\pm,H_Q)$, a.e. $\lambda \in \br.
 $
 Then, in the sequel, we identify these two functions. If Assumption \ref{eq:hyppert} is fulfilled, then the potential $V$ is relatively compact w.r.t. $H_Q$ and by Weyl's criterion on the invariance of the essential spectrum, it follows that
 \begin{equation}
 \sigma_{\rm ess}(H^\pm) = \sigma_{\rm ess}(H_Q) = [0,\tfrac{4}{h^2}]  + \sigma(Q).
 \end{equation}
 However in $\br \smallsetminus \sigma_{\rm ess}(H^\pm)$, the spectrum of $H^\pm$ is purely discrete. Let $\lambda_1 < \lambda_2$ with $[\lambda_1,\lambda_2] \subset \br \smallsetminus \sigma_{\rm ess}(H^\pm)$ and $\lambda_1$, $\lambda_2 \notin \sigma(H^\pm)$. Then, thanks to \cite[Theorem 9.1]{push}, the SSF $\xi(\cdot;H^\pm,H_Q)$ is related to the number of eigenvalues of $H^\pm$ through the formula
 \begin{equation}
 \label{eq:relssfvp}
 \xi(\lambda_1;H^\pm,H_Q) - \xi(\lambda_2;H^\pm,H_Q) = {\rm Rank} \: \mathds 1_{[\lambda_1,\lambda_2)}(H^\pm).
 \end{equation}

\section{Main results}\label{sec5}

\subsection{Statement of the main results} 

Our first theorem is the next simple result which is an immediate by-product of \eqref{eq:foncdecgen}, Lemma \ref{lem:pushCV}, ii) of Proposition \ref{prop:Ts2}, \eqref{eq:pim1}, \eqref{eq:pim2}, Lemma \ref{lem:hold2} and Weyl's inequality \eqref{eq:inegsp}.

\begin{Theorem}
\label{theo0}
Let $V$ satisfy Assumption \ref{eq:hyppert}. Then, the SSF 
is bounded on compact subsets 
$\Gamma \subset \br \smallsetminus \mathscr E_Q$. That is, 
$
\sup_{\lambda \in \Gamma} \xi(\lambda;H^\pm,H_Q) < +\infty.
$
\end{Theorem}

The above result will be useful in Section \ref{sec:corr}. 
In what follows below, to simplify the presentation, our results will be stated for non-degenerate thresholds. However, note that they can be extended to degenerate thresholds (see Remark \ref{seuil:deg}).
Let us introduce some notations.

Recall that the function $\psi \in \ell^2(\mathbb{Z}_h)$ is given by \eqref{def:psi} and let us define the operator $\bra \psi \vert : \ell^2(\bz_h) \to \bc$ so that
\begin{equation}
\label{eq:n0}
\bra \psi \vert^\ast : \zeta \in \bc \mapsto \zeta \psi \in \ell^2(\bz_h).
\end{equation}
We associate to a spectral threshold $\mu_s$, $0 \le s \le d$, the compact operator $L_{s} : \ell^2(\bz_h) \otimes \mathcal{G} \to \bc \otimes \mathcal{G}$ defined by
\begin{equation}
\label{eq:L0}
L_{s} := ( \bra \psi \vert \otimes \pi_s K^\ast ) \mathscr V^{1/2}, \quad \implies L_{s}^\ast = \mathscr V^{1/2} ( \bra \psi \vert^\ast \otimes K \pi_s ) :  \bc \otimes \mathcal{G} \to \ell^2(\bz_h) \otimes \mathcal{G}.
\end{equation}
Let $J$ be the self-adjoint unitary operator defined in $\ell^2(\bz_h)$ by 
\begin{equation}
\label{op:rel0}
(J \varphi)(hn) := (-1)^{n} \varphi(hn).
\end{equation}
Note that $J$ commutes with any multiplication operator. Moreover, it relates both thresholds $0$ and $\frac4{h^2}$ through the relation $J(-\Delta_h)J^\ast = \Delta_h + \frac4{h^2}$. As above, for $0 \le s \le d$, we associate to a spectral threshold $\mu_s + \frac{4}{h^2}$ the compact operator 
\begin{equation}
\label{eq:L4}
L_{4,s} := ( \bra \psi \vert J^\ast \otimes \pi_s K^\ast ) \mathscr V^{1/2} : \ell^2(\bz_h) \otimes \mathcal{G} \to \bc \otimes \mathcal{G},
\end{equation}
so that
\begin{equation*}
L_{4,s}^\ast = \mathscr V^{1/2} ( J \bra \psi \vert^\ast \otimes K \pi_s ) :  \bc \otimes \mathcal{G} \to \ell^2(\bz_h) \otimes \mathcal{G}.
\end{equation*}

\begin{Definition}
For two real-valued functionals $F_1(V,\lambda)$ and $F_2(V,\lambda)$ of $V$ depending on $\lambda \in \br \smallsetminus \mathscr E_Q$, we write 
$$
F_1(V,\lambda) \sim F_2(V,\lambda), \quad \lambda \to \lambda_0 \in \mathscr E_Q,
$$
if for every $\varepsilon \in (0,1)$, we have the estimates
$$
F_2((1-\varepsilon)^{-1}V,\lambda) + \mathcal O_\varepsilon(1) \le F_1(V,\lambda) \le F_2((1+\varepsilon)^{-1}V,\lambda) + \mathcal O_\varepsilon(1), \qquad \lambda \to \lambda_0.
$$
\end{Definition}

\subsubsection{The case $\lambda \nearrow \mu_s$ and $\lambda \searrow \frac4{h^2}+\mu_s$, $0 \le s \le d$}

Our second theorem concerns the asymptotic behavior of the SSF $\xi(\lambda;H^\pm,H_0)$ as $\lambda \to \mu_s$ from below and as $\lambda \to \frac4{h^2}+\mu_s$ from above. Define the operators
\begin{equation*}
\label{proj}
P_s := \bra \psi \vert \otimes \pi_s : \ell^2(\bz_h) \otimes \mathcal{G} \to \bc \otimes \mathcal{G}, \quad 0 \le s \le d,
\end{equation*}
where $\bra \psi \vert$ is defined by \eqref{eq:n0}, and 
\begin{equation}
\label{op:nu+}
\textbf{\textup{V}} := (I_{\ell^2(\bz_h)} \otimes K^\ast) {\mathscr V} (I_{\ell^2(\bz_h)} \otimes K).
\end{equation}
Our results are closely related to the trace class operators 
\begin{equation}
\label{eq:raj+}
P_s \textbf{\textup{V}} P_s^\ast = L_{s} L_{s}^\ast \quad {\rm and} \quad P_s \textbf{\textup{V}}_J P_s^\ast = L_{4,s} L_{4,s}^\ast,
\end{equation}
acting from $\bc \otimes \mathcal{G}$ onto $\bc \otimes \mathcal{G}$, where
\begin{equation}
\label{eq:vj}
\textbf{\textup{V}}_J := (J \otimes I_{\mathcal G})^{\ast} \textbf{\textup{V}} (J \otimes I_{\mathcal G}).
\end{equation}
Therefore, $\textbf{\textup{V}}_J$ is unitarily equivalent to $\textbf{\textup{V}}$. Next, one sets
\begin{equation}
\label{eq:opdom04}
\omega_s(\lambda) :=  \frac{h}{2} \frac{P_s \textbf{\textup{V}} P_s^\ast }{\sqrt{\vert \lambda - \mu_s\vert}} \quad {\rm and} \quad 
\omega_{4,s}(\lambda) := - \frac{h}2 \frac{P_s \textbf{\textup{V}}_J P_s^\ast }{\sqrt{\vert 4/{h^2}+\mu_s-\lambda \vert}}, \quad \lambda \in  \br \smallsetminus \{\mu_s,\tfrac{4}{h^2}+\mu_s\}.
\end{equation}

The following result holds.

\begin{Theorem}
\label{theo1}
Let $V$ satisfy Assumption \ref{eq:hyppert} with $\nu_0 > 3$. Then, for all thresholds $\mu_s$ and $\frac{4}{h^2} + \mu_s \in \mathscr E_Q$, $0 \le s \le d$, non-degenerate for $s \ge 1$, we have:
\begin{itemize}
\item As $\lambda \nearrow \mu_s$,
\begin{equation}
\label{est:h+0-}
 \xi(\lambda;H^+,H_Q) =  \mathcal O(1),
 \end{equation}
 \begin{equation}
 \label{est:h-0-}
\xi(\lambda;H^-,H_Q) \sim - {\rm Tr} \: \mathds 1_{(1,+\infty)} ( \omega_s(\lambda)).
\end{equation}
\item As $\lambda \searrow \frac4{h^2}+\mu_s$,
\begin{equation}
\label{est:h+4+}
-\xi(\lambda;H^+,H_Q) \sim - {\rm Tr} \: \mathds 1_{(1,+\infty)} (-\omega_{4,s}(\lambda)),
 \end{equation}
 \begin{equation}
 \label{est:h-4+}
 \xi(\lambda;H^-,H_Q) =  \mathcal O(1).
\end{equation}
\end{itemize}
\end{Theorem}


\begin{Remark}
\label{rem:th1}
If there exists $0 < \epsilon \ll 1$ with $\sigma_{\rm ess}(H_Q) \cap [\mu_s - \epsilon,\mu_s) = \varnothing$, then \eqref{est:h+0-} and \eqref{eq:relssfvp} imply that the bound states of $H+V$ do not accumulate at $\mu_s$ from the left. Otherwise, \eqref{est:h-0-} 
implies that the problem of counting the number of bound states of the operator $H-V$ near $\mu_s$ from the left, is reduced to the problem of counting the number of eigenvalues of the positive trace class operator $P_s \textbf{\textup{V}} P_s^\ast$ near $0$. If there exists $0 < \alpha \ll 1$ such that $(\frac4{h^2}+\mu_s,\frac4{h^2}+ \mu_s + \alpha] \cap \sigma_{\rm ess}(H_Q) = \varnothing$, then \eqref{est:h+4+} and \eqref{est:h-4+} lead to similar conclusions on the number of bound states of the operators $H \pm V$ near $\frac4{h^2}+\mu_s$ from above. In particular, 
the problem of counting the number of bound states of the operator $H+V$ near $\frac4{h^2}+\mu_s$ from the right, is reduced to the problem of counting the number of eigenvalues of the positive trace class operator $P_s \textbf{\textup{V}}_J P_s^\ast$ near $0$.
\end{Remark}

\subsubsection{The case $\lambda \searrow \mu_s$ and $\lambda \nearrow \frac4{h^2}+\mu_s$, $0 \le s \le d$}

Our third theorem concerns the asymptotic behavior of the SSF $\xi(\lambda;H^\pm,H_Q)$ as $\lambda \to \mu_s$ from above and as $\lambda \to \frac4{h^2}+\mu_s$ from below. One needs first to introduce some notations. Set
\begin{equation}
\label{eq:fctgh}
g_s(\lambda) := \arcsin \Big( \frac{h}{2}\sqrt{\lambda-\mu_s} \Big), \quad \lambda \in (\mu_s,\tfrac{4}{h^2}+\mu_s).
\end{equation}
Let us introduce the operators $\textup{cos}_{\psi,s}$, $\textup{sin}_{\psi,s} : \ell^2(\bz_h) \to \bc$ defined by
\begin{equation}
\label{op:cos}
\textup{cos}_{\psi,s}(\lambda) := \big \bra \psi \cos[2(\cdot)h^{-1} g_s(\lambda)] \big \vert, 
\end{equation}
and
\begin{equation}
\label{op:sin}
\textup{sin}_{\psi,s}(\lambda) := \big \bra \psi \sin[2(\cdot)h^{-1} g_s(\lambda)] \big \vert.
\end{equation}
The adjoints $\textup{cos}_{\psi,s}^\ast$, $\textup{sin}_{\psi,s}^\ast : \bc \to \ell^2(\bz_h)$ are the rank one operators given by
\begin{equation}
\textup{cos}_{\psi,s}(\lambda)^\ast \zeta = \zeta \psi \cos[2(\cdot)h^{-1} g_s(\lambda)],
\end{equation}
and 
\begin{equation}
\textup{sin}_{\psi,s}(\lambda)^\ast \zeta = \zeta \psi \sin[2(\cdot)h^{-1} g_s(\lambda)].
\end{equation}
Define the operator $Y_s(\lambda) : \ell^2(\bz_h) \to \bc^2$ given by
\begin{equation}
\label{eq:dim}
Y_s(\lambda) \phi = \begin{pmatrix}
\textup{cos}_{\psi,s} (\lambda) \phi \\
\textup{sin}_{\psi,s}(\lambda) \phi
\end{pmatrix},
\end{equation}
so that its adjoint $Y_s(\lambda)^\ast : \bc^2 \to \ell^2(\bz_h)$ is given by
\begin{equation}
\label{eq:dimm}
Y_s(\lambda)^\ast \begin{pmatrix}
\zeta_1 \\
\zeta_2
\end{pmatrix} = 
\textup{cos}_{\psi,s}(\lambda)^\ast \zeta_1 + \textup{sin}_{\psi,s}(\lambda)^\ast \zeta_2.
\end{equation}

The next result is closely related to the trace class positive operator 
\begin{equation}
\label{eq:opeomega}
\Omega_s(\lambda) = \frac{1}{\sqrt{\lambda-\mu_s} \sqrt{4/{h^2}+\mu_s-\lambda}} (Y_s(\lambda) \otimes \pi_s) \textbf{\textup{V}} (Y_s(\lambda)^\ast \otimes \pi_s) : \bc^2 \otimes \mathcal{G} \to \bc^2 \otimes \mathcal{G},
\end{equation}
where $\textbf{\textup{V}}$ is given by \eqref{op:nu+}.

\begin{Theorem}
\label{theo2}
Let $V$ satisfy Assumption \ref{eq:hyppert} with $\nu_0 > 3$. Then, for all thresholds $\mu_s$ and $\frac{4}{h^2} + \mu_s \in \mathscr E_Q$, $0 \le s \le d$, non-degenerate for $s \ge 1$, we have
\begin{equation}
\mp \xi(\lambda;H^\pm,H_Q) \sim -\frac1{\pi} {\rm Tr} \arctan(\Omega_s(\lambda)),
 \end{equation}
as $\lambda \searrow \mu_s$ and $\lambda \nearrow \frac4{h^2}+\mu_s$.
\end{Theorem}


\begin{Remark}
\label{rem:corr}
Under the conditions of Theorem \ref{theo2} and for $x > 0$, the operator $\Omega_s(\lambda)$ satisfies
\begin{equation}
\label{as:lifs-kre0}
{\rm Tr} \arctan(x^{-1} \Omega_s(\lambda)) = {\rm Tr} \arctan(x^{-1} \Omega_{0,s}(\lambda)) + \mathcal O(1), \qquad \lambda \searrow \mu_s,
\end{equation}
where the operator $\Omega_{0,s}(\lambda)$ is given by
\begin{equation}
\label{eq:Y0nu}
\Omega_{0,s}(\lambda) := \frac{h}{2\sqrt{\lambda-\mu_s}} (Y_{0} \otimes \pi_s) \textbf{\textup{V}} (Y_{0}^\ast \otimes \pi_s), \quad Y_{0} = \begin{pmatrix}
\bra \psi \vert \\
0
\end{pmatrix} : \ell^2(\bz_h) \to \bc^2,
\end{equation}
with $\bra \psi \vert$ defined by \eqref{eq:n0}. The estimate \eqref{as:lifs-kre0} follows from the Lifshits-Krein trace formula \eqref{eq:fctdec}.
 The details of the proof are given in Section \ref{sec9}
  and the argument is analogous to the one of \cite[Corollary 2.2]{ferrai}. Now, using \eqref{eq:eqarctan}, one gets
\begin{equation}
\label{as:lifs-kre0asy}
{\rm Tr} \arctan(x^{-1} \Omega_{0,s}(\lambda)) = {\rm Tr} \arctan (x^{-1} \omega_s(\lambda)),
\end{equation}
the operator $\omega_s(\lambda)$ being defined by \eqref{eq:opdom04}. 
In particular, it follows from \eqref{as:lifs-kre0} and \eqref{as:lifs-kre0asy} that Theorem \ref{theo2} can be formulated when $\lambda \searrow \mu_s$ in terms of the operator $\omega_s(\lambda)$, as in Theorem \ref{theo1}.
In a similar way, one has for $x > 0$
\begin{equation}
\label{as:lifs-kre4}
{\rm Tr} \arctan(x^{-1} \Omega_s(\lambda)) = {\rm Tr} \arctan(x^{-1} \Omega_{4,s}(\lambda)) + \mathcal O(1), \qquad \lambda \nearrow \tfrac4{h^2}+\mu_s,
\end{equation}
where the operator $\Omega_{4,s}(\lambda)$ is given 
\begin{equation*}
\Omega_{4,s}(\lambda) := \frac{h}{2\sqrt{4/{h^2}+\mu_s-\lambda}} (Y_{4} \otimes \pi_s) \textbf{\textup{V}} (Y_{4}^\ast \otimes \pi_s), \quad Y_{4} = \begin{pmatrix}
\bra \psi \vert J^{\ast} \\
0
\end{pmatrix} : \ell^2(\bz_h) \to \bc^2,
\end{equation*}
with $J$ defined by \eqref{op:rel0}. By using \eqref{eq:eqarctan}, one obtains
\begin{equation}
\label{as:lifs-kre4asy}
{\rm Tr} \arctan(x^{-1} \Omega_{4,s}(\lambda)) = {\rm Tr} \arctan (-x^{-1} \omega_{4,s}(\lambda)),
\end{equation}
where $\omega_{4,s}(\lambda)$ is defined by \eqref{eq:opdom04}. In particular, it follows from \eqref{as:lifs-kre4} and \eqref{as:lifs-kre4asy} that Theorem \ref{theo2} can be formulated when $\lambda \nearrow \frac4{h^2}+\mu_s$ in terms of the operator $\omega_{4,s}(\lambda)$, as in Theorem \ref{theo1}.
\end{Remark}

It follows from Theorem \ref{theo2} and Remark \ref{rem:corr} the following result.

\begin{Theorem}
\label{theo2+}
Let $V$ satisfy Assumption \ref{eq:hyppert} with $\nu_0 > 3$. Then, for all thresholds $\mu_s$ and $\frac{4}{h^2} + \mu_s \in \mathscr E_Q$, $0 \le s \le d$, non-degenerate for $s \ge 1$, we have
$$
\mp \xi(\lambda;H^\pm,H_Q) \sim -\frac1{\pi} {\rm Tr} \arctan(\omega_s(\lambda)), \quad \lambda \searrow \mu_s,
$$ 
and
$$
\mp \xi(\lambda;H^\pm,H_Q) \sim -\frac1{\pi} {\rm Tr} \arctan(-\omega_{4,s}(\lambda)), \quad \lambda \nearrow \tfrac4{h^2}+\mu_s.
$$ 
\end{Theorem}

\begin{Remark}
\label{seuil:deg}
For $Q \neq 0$, i.e. $s \ge 1$, similar results to Theorems \ref{theo1}, \ref{theo2} and \ref{theo2+} can be established for degenerate thresholds. For such thresholds $\mu_s = \frac{4}{h^2} + \mu_{s'}$, $s \neq s'$, the asymptotics are given by expressions involving both the operators $\omega_s(\lambda)$ and $\omega_{4,s'}(\lambda)$ as $\lambda \to \mu_s = \frac{4}{h^2} + \mu_{s'}$.
\end{Remark}

\subsection{Corollaries}\label{sec:corr}

In this section, we present some consequences of the above results gathered in two parts.
We will see that in the first part (the case $\dim(\mathcal G) < + \infty$), the SSF is bounded at the spectral thresholds of the essential spectrum while it may have singularities in the second one (the case $\dim(\mathcal G) = +\infty$).

\subsubsection{Boundedness of the SSF at the spectral thresholds}

We assume that $\dim(\mathcal G) < +\infty$. Then, the operators $P_s \textbf{\textup{V}} P_s^\ast$ et $P_s \textbf{\textup{V}}_J P_s^\ast$ with $0 \le s \le d$ acting from $\bc \otimes \mathcal G$ onto $\bc \otimes \mathcal G$ are of finite rank. Otherwise, for $x > 0$ we have
\begin{equation}
\label{eq:supp}
{\rm Tr} \arctan (x^{-1} \omega_{s}(\lambda)) = \int_{\br_+} \mathscr N_+\big( x\sqrt{\vert \lambda - \mu_s \vert}t,(h/2) P_s \textbf{\textup{V}} P_s^\ast \big)  \frac{dt}{1+t^2}, 
\end{equation}
\begin{equation}
\label{eq:supp+}
{\rm Tr} \arctan (-x^{-1} \omega_{4,s}(\lambda)) = \int_{\br_+} \mathscr N_+\big( x\sqrt{\vert 4/h^2+\mu_s-\lambda \vert}t,(h/2) P_s \textbf{\textup{V}}_J P_s^\ast \big)  \frac{dt}{1+t^2}.
\end{equation}
Together with Theorems \ref{theo1}, \ref{theo2+}, \ref{theo0} and Remark \ref{seuil:deg}, this implies the following result.

\begin{Theorem}
\label{theo1:cor1}
Let $V$ satisfy Assumption \ref{eq:hyppert} with $\nu_0 > 3$. Suppose that $\dim(\mathcal G) < + \infty$. Then, we have
$\sup_{\lambda \in \br \smallsetminus \mathscr E_Q} \xi(\lambda;H^\pm,H_Q) < +\infty.$
\end{Theorem}


\begin{Corollary}
\label{cor:specdisfini}
If $V$ satisfies Assumption \ref{eq:hyppert} with $\nu_0 > 3$ and $\dim(\mathcal G) < +\infty$, then:
\begin{itemize}
\item $\sigma_{\rm ess}(H_Q \pm V) =  \sigma_{\rm ess}(H_Q) = [0,\frac4{h^2}] + \sigma(Q)$. 
\item The bound states of the operators $H_Q \pm V$ do not accumulate to any point of $\mathscr E_Q$. In particular, it follows that $\# \sigma_{\rm disc}(H_Q \pm V) < \infty.$ 
\end{itemize}
\end{Corollary}

\begin{Remark}
\label{r:finit:spec}
Thanks to Lemma \ref{l:sdecV}, it follows from Corollary \ref{cor:specdisfini} that the bound states of the perturbed operators $H_Q \pm V$ do not accumulate to $\mathscr E_Q$, under Assumption \ref{eq:hyppert0} with $\nu_i > 3$, $i=1$, $2$. In particular, 
the Schr\"odinger operator $-\Delta_h \pm V$ corresponding to $Q=0$  satisfies
$\# \sigma_{\rm disc}(-\Delta_h \pm V) < \infty$.
This and Corollary \ref{cor:specdisfini} 
can be compare to \cite[Corollary 2.1]{bst} where we prove the finiteness of $\sigma_{\rm disc}(-\Delta_1 \pm V)$ 
for $V$ self-adjoint exponentially decaying at infinity.
\end{Remark}

\subsubsection{Thresholds singularities and asymptotic behaviors of the SSF}

Here, we assume that $\dim(\mathcal G) = + \infty$. To exhibit singularities of the SSF, we focus our analysis near the thresholds $\mu_s$ and $\frac{4}{h^2} + \mu_s$ such that $P_s \textbf{\textup{V}} P_s^\ast$ and $P_s \textbf{\textup{V}}_J P_s^\ast$ are of infinite rank. Indeed, as in the previous section, notice that for $Q \neq 0$ the SSF is bounded near the spectral thresholds $\{\mu_s,\frac{4}{h^2} + \mu_s\}$ such that ${\rm Rank} \: \pi_s < +\infty$, i.e. if $\mu_s$ is a bounded state. More generally, we have:

\begin{Remark}
\label{eq:rem5-}
Theorem \ref{theo1:cor1} and Corollary \ref{cor:specdisfini} remain valid if we assume that the positive trace class operators $P_s \textbf{\textup{V}} P_s^\ast$ and $P_s \textbf{\textup{V}}_J P_s^\ast$, $0 \le s \le d$, acting from $\bc \otimes \mathcal G$ onto $\bc \otimes \mathcal G$ are finite-rank. For instance, this holds when $K$ acting in $\mathcal G$ is finite-rank (see also Remark \ref{rem0}). 
\end{Remark}

\medskip

\noindent
$-$ \textbf{\textit{Case}} $\lambda \nearrow \mu_s$ and $\lambda \searrow \frac4{h^2}+\mu_s$:
A direct consequence of Theorem \ref{theo1} is the following result.


\begin{Theorem}
\label{cor3:mr}
Under the assumption of Theorem \ref{theo1}, fix $0 \le s \le d$ such that ${\rm Rank} \:P_s \textbf{\textup{V}} P_s^\ast = {\rm Rank} \:P_s \textbf{\textup{V}}_J P_s^\ast = +\infty$. Then, the SSF admits singularities at the thresholds $\mu_s$ and $\frac4{h^2}+\mu_s$ with
$$
\xi(\lambda;H^-,H_Q) 
\sim \begin{cases}
-{\rm Tr} \: \mathds 1_{(\frac2h\sqrt{\mu_s-\lambda},+\infty)} (P_s \textbf{\textup{V}} P_s^\ast) \quad & as \quad \lambda \nearrow \mu_s, \\
\mathcal O(1) \quad & as \quad \lambda \searrow \frac4{h^2}+\mu_s,
\end{cases}
$$
and 
$$
-\xi (\lambda;H^+,H_Q) 
\sim \begin{cases}
\mathcal O(1) \quad & as \quad \lambda \nearrow \mu_s, \\
- {\rm Tr} \: \mathds 1_{(\frac2h\sqrt{\lambda-4/{h^2}-\mu_s},+\infty)} (P_s \textbf{\textup{V}}_J P_s^\ast) \quad & as \quad \lambda \searrow \frac4{h^2}+\mu_s.
\end{cases}
$$
\end{Theorem}

\begin{Remark}[Accumulation of bound states]
\label{rem3:mr}
Remark \ref{rem:th1} and Theorem \ref{cor3:mr} show that: 
\begin{enumerate}
\item If there exists $0 < \epsilon \ll 1$ such that $\sigma_{\rm ess}(H_Q) \cap [\mu_s - \epsilon,\mu_s) = \varnothing$,
then the operator $H^-$ has infinitely many discrete eigenvalues below $\mu_s$, with a rate of accumulation close to
 ${\rm Tr} \: \mathds 1_{(\frac2h\sqrt{\mu_s-\lambda},+\infty)} (P_s \textbf{\textup{V}} P_s^\ast)$ as $\lambda \nearrow \mu_s$.
\item If there exists $0 < \alpha \ll 1$ such that $(\frac4{h^2}+\mu_s,\frac4{h^2}+ \mu_s + \alpha] \cap \sigma_{\rm ess}(H_Q) = \varnothing$,
then the operator $H^+$ has infinitely many discrete eigenvalues above $\frac4{h^2}+\mu_s$, with a rate of accumulation close to 
${\rm Tr} \: \mathds 1_{(\frac2h\sqrt{\lambda-4/{h^2}-\mu_s},+\infty)} (P_s \textbf{\textup{V}}_J P_s^\ast)$ as $\lambda \searrow \tfrac4{h^2}+\mu_s$.
\item For $Q=0$ i.e. $s=0$ so that $\mu_0 = 0$ and $\sigma_{\rm ess}(H_Q) = [0,\frac4{h^2}]$, one has $\xi(\lambda;H^-,H_0) = - {\rm Tr} \: \mathds 1_{(-\infty,\lambda)}(H^-)$, $\lambda < 0$, and $\xi(\lambda;H^+,H_0) = {\rm Tr} \: \mathds 1_{(\lambda,+\infty)}(H^+)$, $\lambda > \frac4{h^2}$.
\end{enumerate}
\end{Remark}

It is not difficult to construct potentials $V$ such that ${\rm Rank} \:P_s \textbf{\textup{V}} P_s^\ast = {\rm Rank} \: P_s \textbf{\textup{V}}_J P_s^\ast = +\infty$. But, it is more interesting to investigate cases where we can obtain a more precise description of the asymptotic behavior of $\xi (\lambda;H^\pm,H_Q)$ near the spectral thresholds $\mu_s$ and $\frac{4}{h^2}+\mu_s$. So, in what follows below, one sets for $r > 0$
\begin{equation}
\label{eq:fctphi}
\Phi_1(r) := {\rm Tr} \: \mathds 1_{(r,+\infty)} (P_s \textbf{\textup{V}} P_s^\ast) \quad {\rm and} \quad \Phi_2(r) := {\rm Tr} \: \mathds 1_{(r,+\infty)} (P_s \textbf{\textup{V}}_J P_s^\ast ).
\end{equation}

Actually, under additional conditions on the functions $\Phi_i$, $i=1, 2$, Theorems \ref{cor3:mr}
 produces the next more precise result that gives the main terms of the asymptotic behaviors of $\xi(\lambda;H^\pm,H_Q)$, $\lambda \nearrow \mu_s$ and $\lambda \searrow \frac4{h^2}+\mu_s$.

\begin{Corollary}
\label{cor4:mr}
Under the assumptions of Theorem \ref{cor3:mr}, suppose in addition that for any $\varepsilon \in (0,1)$ small, $\Phi_i(r(1 \pm \varepsilon)) = \Psi_i(r)(1 + o(1) + \mathcal O(\varepsilon))$ as $r \searrow 0$, $i=1$, $2$, with $\Psi_i(r) \to + \infty$ as $r \to 0$. Then, one has the asymptotics
$$
\xi(\lambda;H^-,H_Q) 
=\begin{cases}
-\Psi_1((2/h)\sqrt{\mu_s-\lambda})(1 + o(1)) \quad & as \quad \lambda \nearrow \mu_s, \\
\mathcal O(1) \quad & as \quad \lambda \searrow \frac4{h^2}+\mu_s,
\end{cases}
$$
and 
$$
\xi (\lambda;H^+,H_Q) 
= \begin{cases}
\mathcal O(1) \quad & as \quad \lambda \nearrow \mu_s, \\
\Psi_2((2/h)\sqrt{\lambda-4/{h^2}-\mu_s})(1 + o(1)) \quad & as \quad \lambda \searrow \frac4{h^2}+\mu_s.
\end{cases}
$$
\end{Corollary}

\begin{Remark}
Examples of such $\Psi_i$, $i = 1, 2$ of Corollary \ref{cor4:mr} are given by
\begin{equation}
\label{eq:fctpsi}
\Psi_i(r) = \Psi(r) = 
\begin{cases}
C_0 r^{-\alpha}, \quad \alpha > 0 \\
C_0 \vert \ln r \vert^\alpha, \quad \alpha > 0 \\
C_0 (\ln\vert \ln r \vert)^\alpha, \quad \alpha > 0 \\
C_0 \vert \ln r \vert (\ln\vert \ln r \vert)^{-1},
\end{cases}, \quad r > 0,
\end{equation}
where $C_0 > 0$ is a constant (see \cite[Proof of Corollary 3.11]{bo}). For more details, we give examples of explicit computations of $\Phi_1(r)$, $\Phi_2(r)$ and $\Psi(r)$ in Section \ref{sec:exple}, including polynomial and exponential decay potentials along the component $\mathcal G$ of $\ell^2(\bz_h,\mathcal G)$ (see Propositions \ref{prop:exem1} and \ref{exemp2}).
\end{Remark}

Taking into account the previous remark, the next result holds. 

\begin{Corollary}
\label{cor:lev-}
Set $z_- := \mu_s$, $z_+ := \frac4{h^2}+\mu_s$. The following holds w.r.t. $\pm$. 
\begin{enumerate}
\item If $V$ satisfies the assumptions of Proposition \ref{prop:exem1} with $\nu_0 > 3$, then 
$$
\xi(z_\pm \pm \lambda;H^\pm,H_Q) = \pm \Big( \sum_{n \in \bz} \bra n \ket^{-\nu_0} \Big)^{1/\beta_0} ( 2/h )^{-1/\beta_0} (\sqrt\lambda)^{-1/\beta_0} (1 + o(1)), \quad \lambda \searrow 0.
$$
\item If $V$ satisfies the assumptions of Proposition \ref{exemp2} with $\nu_0 > 3$, then: 

{\rm i)} If $\xi(j) = \eta j^\beta$, $\eta > 0$ and $\beta > 0$, we have
$$
\xi(z_\pm \pm \lambda;H^\pm,H_Q) = \pm (2/\eta)^{1/\beta} \vert \ln \sqrt{\lambda} \vert^{1/\beta} (1 + o(1)), \quad \lambda \searrow 0.
$$

{\rm ii)} If $\xi(j) = e^{\eta j^\beta}$, $\eta > 0$ and $\beta > 0$, we have
$$
\xi(z_\pm \pm \lambda;H^\pm,H_Q) = \pm \eta^{-1/\beta} ( \ln \vert \ln \sqrt{\lambda} \vert )^{1/\beta} (1 + o(1)), \quad \lambda \searrow 0.
$$

{\rm iii)} If $\xi(j) = \chi_\eta^{-1}(j)$, $\eta > 0$, we have
$$
\xi(z_\pm \pm \lambda;H^\pm,H_Q) = \pm 2\eta^{-1} \vert \ln \sqrt{\lambda} \vert (\ln \vert \ln \sqrt{\lambda} \vert )^{-1} (1 + o(1)), \quad \lambda \searrow 0.
$$
\end{enumerate}
\end{Corollary}

\medskip

\noindent
$-$ \textbf{\textit{Case}} $\lambda \searrow \mu_s$ and $\lambda \nearrow \frac4{h^2}+\mu_s$:
Formulas \eqref{eq:supp} and \eqref{eq:supp+} can be rewritten as
\begin{equation}
\label{eq:supp1}
{\rm Tr} \arctan (x^{-1} \omega_s(\lambda)) = \int_{\br_+} \Phi_1\big((2/h) x\sqrt{\lambda - \mu_s}t \big)  \frac{dt}{1+t^2},
\end{equation}
\begin{equation}
\label{eq:supp2}
{\rm Tr} \arctan (-x^{-1} \omega_{4,s}(\lambda)) = \int_{\br_+} \Phi_2 \big( (2/h)x\sqrt{4/h^2+\mu_s-\lambda }t \big)  \frac{dt}{1+t^2}.
\end{equation}

As above, if the functions $\Phi_i$, $i=1, 2$ verify some asymptotics behaviors near $0$, then Theorem \ref{theo2+} together with \eqref{eq:supp1} ans \eqref{eq:supp2} produces the next more precise result that gives the main terms of the asymptotics of $\xi(\lambda;H^\pm,H_Q)$ as $\lambda \searrow \mu_s$ and $\lambda \nearrow \frac4{h^2}+\mu_s$.

\begin{Theorem}
\label{theo:fin}
Let $V$ satisfy Assumption \ref{eq:hyppert} with $\nu_0 > 3$. Suppose in addition that $\Phi_i(r) = \Psi_i(r)(1 + o(1))$ as $r \searrow 0$, $i=1$, $2$, with $\Psi_i(r)$ given by \eqref{eq:fctpsi}. Set $z_1:= \mu_s$ and $z_2 := \frac4{h^2}+\mu_s$. Then, one has the following asymptotics near $\mu_s$ from above and $\frac4{h^2}+\mu_s$ from below.
\begin{itemize}
\item[\rm i)] If $\Psi_i(r) = C_0r^{-\alpha}$, $0 < \alpha < 1$, then
$$
\xi(\lambda;H^\pm,H_Q) = \pm \frac{1}{2\cos(\frac\pi\alpha)} \Psi_i \big( (2/h) \sqrt{\vert z_i - \lambda \vert} \big) (1 + o(1)), \quad \lambda \to z_i.
$$
\item[\rm ii)] If $\Psi_i(r) = C_0 \vert \ln r \vert^{\alpha}$, or $C_0 (\ln\vert \ln r \vert)^{\alpha}$, or $C_0 \vert \ln r \vert (\ln\vert \ln r \vert)^{-1}$, then
$$
\xi(\lambda;H^\pm,H_Q) = \pm \frac{1}{2} \Psi_i \big( (2/h) \sqrt{\vert z_i - \lambda \vert} \big) (1 + o(1)), \quad \lambda \to z_i.
$$
\end{itemize}
\end{Theorem}

\begin{Corollary}
\label{cor:fin1}
\begin{enumerate}
\item Let $V$ satisfy the assumptions of Proposition \ref{prop:exem1} with $\nu_0 > 3$. Then, as $\lambda \to z = \mu_s$ from above and $\lambda \to z = \frac4{h^2}+\mu_s$ from below, we have
$$
\xi(\lambda;H^\pm,H_Q) = \pm \frac{1}{2\cos(\pi\beta_0)} \Big( \sum_{n \in \bz} \bra n \ket^{-\nu_0} \Big)^{1/\beta_0} (2/h)^{-1/\beta_0} ( \sqrt{\vert z - \lambda \vert})^{-1/\beta_0} (1 + o(1)).
$$
\item Suppose that $V$ satisfies the assumptions of Proposition \ref{exemp2} with $\nu_0 > 3$. Then, as $\lambda \to z = \mu_s$ from above and $\lambda \to z = \frac4{h^2}+\mu_s$ from below, one has:

\noindent
{\rm i)} If $\xi(j) = \eta j^\beta$, $\eta > 0$ and $\beta > 0$, 
$$
\xi(\lambda;H^\pm,H_Q) = \pm\frac12 ( 2/\eta )^{1/\beta} \vert \ln \sqrt{\vert z - \lambda \vert} \vert^{1/\beta} (1 + o(1)).
$$
{\rm ii)} If $\xi(j) = e^{\eta j^\beta}$, $\eta > 0$ and $\beta > 0$, 
$$
\xi(\lambda;H^\pm,H_Q) = \pm\frac12 \eta^{-1/\beta} \big( \ln \vert \ln \sqrt{\vert z - \lambda \vert} \vert \big)^{1/\beta} (1 + o(1)).
$$
{\rm iii)} If $\xi(j) = \chi_\eta^{-1}(j)$, $\eta > 0$, 
$$
\xi(\lambda;H^\pm,H_Q) = \pm \eta^{-1} \vert \ln \sqrt{\vert z - \lambda \vert} \vert \big( \ln \vert \ln \sqrt{\vert z - \lambda \vert} \vert \big)^{-1} (1 + o(1)).
$$
\end{enumerate}
\end{Corollary}

\begin{Remark}
By \eqref{eq:scatdet}, Theorems \ref{theo1}, \ref{theo2}, \ref{theo2+}, \ref{cor3:mr} and \ref{theo:fin} concern  the asymptotics of the scattering phase $\arg \det S(\lambda;H^\pm,H_Q)$ near the spectral thresholds $\mu_s$ and $\frac4{h^2} + \mu_s$.
\end{Remark}

\subsubsection{Levinson type formulas}

Combining Corollary \ref{cor4:mr} and Theorem \ref{theo:fin}, one obtains the next result which can be interpreted as generalized Levinson formulae.

\begin{Theorem}
\label{theo:levin1}
Under the assumptions of Theorem \ref{theo:fin}, i), one has
\begin{equation}
\lim_{\lambda \searrow 0} \frac{\xi(\mu_s+\lambda;H^-,H_Q)}{\xi(\mu_s-\lambda;H^-,H_Q)} = \frac{1}{2\cos(\frac\pi\alpha)} = \lim_{\lambda \searrow 0} \frac{\xi(\frac4{h^2}+\mu_s-\lambda;H^+,H_Q)}{\xi(\frac4{h^2}+\mu_s+\lambda;H^+,H_Q)},
\end{equation}
while under the assumptions of Theorem \ref{theo:fin}, ii), one has
\begin{equation}
\lim_{\lambda \searrow 0} \frac{\xi(\mu_s + \lambda;H^-,H_Q)}{\xi(\mu_s-\lambda;H^-,H_Q)} = \frac{1}{2} = \lim_{\lambda \searrow 0} \frac{\xi(\frac4{h^2}+\mu_s-\lambda;H^+,H_Q)}{\xi(\frac4{h^2}+\mu_s+\lambda;H^+,H_Q)}.
\end{equation}
\end{Theorem}

In particular, Corollaries \ref{cor:lev-} and \ref{cor:fin1} give the following result.

\begin{Corollary}
\label{cor:levin1}
\begin{enumerate}
\item Let $V$ satisfy the assumptions of Proposition \ref{prop:exem1} with $\nu_0 > 3$. Then,
\begin{equation}
\lim_{\lambda \searrow 0} \frac{\xi(\mu_s+\lambda;H^-,H_Q)}{\xi(\mu_s-\lambda;H^-,H_Q)} = \frac{1}{2\cos(\pi\beta_0)} = \lim_{\lambda \searrow 0} \frac{\xi(\frac4{h^2}+\mu_s-\lambda;H^+,H_Q)}{\xi(\frac4{h^2}+\mu_s+\lambda;H^+,H_Q)}.
\end{equation}
\item Let $V$ satisfy the assumptions of Proposition \ref{exemp2} with $\nu_0 > 3$. Then,
\begin{equation}
\lim_{\lambda \searrow 0} \frac{\xi(\mu_s+\lambda;H^-,H_Q)}{\xi(\mu_s-\lambda;H^-,H_Q)} = \frac{1}{2} = \lim_{\lambda \searrow 0} \frac{\xi(\frac4{h^2}+\mu_s-\lambda;H^+,H_Q)}{\xi(\frac4{h^2}+\mu_s+\lambda;H^+,H_Q)}.
\end{equation}
\end{enumerate}
\end{Corollary}

\subsection{Examples of eigenvalues asymptotics for power-like and exponential decay potentials} \label{sec:exple}

In this section, one gives examples of asymptotics, inspired by \cite{rawa1,rawa2}, of the quantities $\Phi_1(r)$, $\Phi_2(r)$ and $\Psi(r)$ defined by \eqref{eq:fctphi} and \eqref{eq:fctpsi}. Here, $a_r \underset{r \to 0}{\sim} b_r$ means that $\frac{a_r}{b_r} \to 1$ as $r \to 0$.


\subsubsection{Polynomial decay potentials}

\begin{Proposition}
\label{prop:exem1}
Let $\nu_0 > 1$, $\beta_0 > 1$, 
$
V = \sum_{(n,m)\in\bz^2} | \delta_n \rangle \langle \delta_m | \otimes V_h(n,m)
$ 
with $V_h(n,m) = 0$ if $n \neq m$ and
$
V_h(n,n) = \bra n \ket^{-\nu_0} \sum_{j \in \bz_+} \langle j \rangle^{-\beta_0} |\e_j \rangle \langle \e_j|.
$ 
For $0 \le s \le d$ with ${\rm Rank} \: \pi_s = +\infty$, we have
$$\Phi_1(r) = \Phi_2(r) = {\rm Tr} \: \mathds 1_{(r,+\infty)} (P_s \textbf{\textup{V}} P_s^\ast),$$
and 
\begin{equation}
{\rm Tr} \: \mathds 1_{(r,+\infty)} (P_s \textbf{\textup{V}} P_s^\ast) = \Big( \sum_{n \in \bz} \bra n \ket^{-\nu_0} \Big)^{1/\beta_0} r^{-1/\beta_0} (1 + o(1)), \qquad r \searrow 0.
\end{equation}
\end{Proposition}

\begin{proof}
Clearly, $V$ satisfies Assumption \ref{eq:hyppert0} if $\nu_0 > 2$, $\beta_0 > 2$. Moreover, recalling that the operator $M_\psi$ and $p$ are respectively given by \eqref{def:opn} and \eqref{def:opp}, one obtains
\begin{equation*}
V = \sum_{(n,j)\in\bz \times \bz_+} \bra n \ket^{-\nu_0} \langle j \rangle^{-\beta_0} | \delta_n \otimes \e_j \rangle \langle \delta_n \otimes \e_j | = ( M_\psi \otimes p )^2.
\end{equation*}
Then, $V$ fulfills Assumption \ref{eq:hyppert} with $K^\ast = K = p$ and $\mathscr V = I$. It follows that 
$$
\textbf{\textup{V}}_J = \textbf{\textup{V}} = (I_{\ell^2(\bz_h)} \otimes p)^2 = I_{\ell^2(\bz_h)} \otimes \sum_{j \in \bz_+} \langle j \rangle^{-\beta_0} |\e_j \rangle \langle \e_j|, 
$$ 
so that $P_s \textbf{\textup{V}} P_s^\ast : \bc \otimes \mathcal{G} \to \bc \otimes \mathcal{G}$ is given by
\begin{equation*}
P_s \textbf{\textup{V}} P_s^\ast = \bra \psi \vert \bra \psi \vert^\ast \otimes \pi_s \Big( \sum_{j \in \bz_+} \langle j \rangle^{-\beta_0} |\e_j \rangle \langle \e_j| \Big) \pi_s
= \Vert \psi \Vert_{\ell^2(\bz_h)}^2 \Big( I_{\bc} \otimes \sum_{j \ge j_s} \langle j \rangle^{-\beta_0} |\e_j \rangle \langle \e_j| \Big),
\end{equation*}
where $(\e_{j})_{j \ge j_s}$ an orthonormal basis of ${\rm Ker}(Q - \mu_s)$ (if $Q \neq 0$).
We can see that the non-zero eigenvalues of the operator $P_0 \textbf{\textup{V}} P_0^\ast$ are simple and
\begin{equation}
\sigma(P_s \textbf{\textup{V}} P_s^\ast) = \bigr \{ \Vert \psi \Vert_{\ell^2(\bz_h)}^2 \langle j \rangle^{-\beta_0} : j \ge j_s \bigl \} \cup \{0\}.
\end{equation}
Hence, for $r > 0$ small enough, one has
\begin{align*}
\Phi_2(r) & = \Phi_1(r) = {\rm Tr} \: \mathds 1_{(r,+\infty)} (P_s \textbf{\textup{V}} P_s^\ast) \\
& = \# \bigr \{ \Vert \psi \Vert_{\ell^2(\bz_h)}^2 \langle j \rangle^{-\beta_0} : j \ge j_s : r < \Vert \psi \Vert_{\ell^2(\bz_h)}^2 \langle j \rangle^{-\beta_0} \bigl \} \\
& = \# \bigr \{ j \ge j_s : j < \big( \Vert \psi \Vert_{\ell^2(\bz_h)}^{4/\beta_0} r^{-2/\beta_0} - 1 \big)^{1/2}  \bigl \}. 
\end{align*}
By denoting $\lfloor x \rfloor$ the integer part of $x \in \br$, one obtains finally
\begin{align*}
\Phi_2(r) & = \Phi_1(r) = {\rm Tr} \: \mathds 1_{(r,+\infty)} (P_s \textbf{\textup{V}} P_s^\ast) \underset{r \to 0}{\sim} \big \lfloor \big( \Vert \psi \Vert_{\ell^2(\bz_h)}^{4/\beta_0} r^{-2/\beta_0} - 1 \big)^{1/2} \big\rfloor \\
& \underset{r \to 0}{\sim} \big( \Vert \psi \Vert_{\ell^2(\bz_h)}^{4/\beta_0} r^{-2/\beta_0} - 1 \big)^{1/2} \underset{r \to 0}{\sim} \Vert \psi \Vert_{\ell^2(\bz_h)}^{2/\beta_0} r^{-1/\beta_0}.
\end{align*}
The claim follows by noting that $\Vert \psi \Vert_{\ell^2(\bz_h)}^{2/\beta_0} = \big( \sum_{n \in \bz} \bra n \ket^{-\nu_0} \big)^{1/\beta_0}$.
\end{proof}

\subsubsection{Exponential decay potentials along the component $\mathcal G$ of $\ell^2(\bz_h,\mathcal G)$}

Consider $\xi$ an increasing unbounded real-valued function of the form
\begin{equation}
\label{eq:fctg}
\xi(x) = 
\begin{cases}
\eta x^\beta, \quad \eta > 0, \quad \beta > 0 \\
e^{\eta x^\beta}, \quad \eta > 0, \quad \beta > 0 \\
 \chi_\eta^{-1}(x), \quad \eta > 0
\end{cases}, \quad x > 0,
\end{equation}
where $\chi_\eta^{-1}$ is the inverse of the function $y \mapsto \chi_\eta(y) = \frac{\eta y}{\ln (y+2)}, y > 0$. 

\begin{Proposition}
\label{exemp2}
Let $\nu_0 > 1$ and consider the potential 
$
V = \sum_{(n,m)\in\bz^2} | \delta_n \rangle \langle \delta_m | \otimes V_h(n,m)
$ 
such that $V_h(n,m) = 0$ if $n \neq m$ and
$
V_h(n,n) = \bra n \ket^{-\nu_0} \sum_{j\in\bz_+} e^{-\xi(j)} |\e_j \rangle \langle \e_j|.
$
Then, for $0 \le s \le d$ such that ${\rm Rank} \: \pi_s = +\infty$, one has
$$
\Phi_1(r) = \Phi_2(r) = {\rm Tr} \: \mathds 1_{(r,+\infty)} (P_s \textbf{\textup{V}} P_s^\ast).$$
Furthermore, 
\begin{itemize}
\item If $\xi(j) = \eta j^\beta$, $\eta > 0$ and $\beta > 0$, 
\begin{equation}
{\rm Tr} \: \mathds 1_{(r,+\infty)} (P_s \textbf{\textup{V}} P_s^\ast) = ( 2/\eta )^{1/\beta} \vert \ln r \vert^{1/\beta} (1 + o(1)), \qquad r \searrow 0.
\end{equation}
\item If $\xi(j) = e^{\eta j^\beta}$, $\eta > 0$ and $\beta > 0$, 
\begin{equation}
{\rm Tr} \: \mathds 1_{(r,+\infty)} (P_s \textbf{\textup{V}} P_s^\ast) = \eta^{-1/\beta} (\ln \vert \ln r \vert )^{1/\beta} (1 + o(1)), \qquad r \searrow 0.
\end{equation}
\item If $\xi(j) = \chi_\eta^{-1}(j)$, $\eta > 0$, 
\begin{equation}
{\rm Tr} \: \mathds 1_{(r,+\infty)} (P_s \textbf{\textup{V}} P_s^\ast) = 2\eta^{-1} \vert \ln r \vert (\ln\vert \ln r \vert)^{-1} (1 + o(1)), \qquad r \searrow 0.
\end{equation}
\end{itemize}
\end{Proposition}

\begin{proof}
It can be checked, making the change of variable $x = \chi_\eta(y)$, that
$$
\lim_{x \to +\infty} \frac{\chi_\eta^{-1}(x)}{x \ln(x)} = \frac1{\eta}.
$$
Hence, $V$ satisfies Assumption \ref{eq:hyppert0} if $\nu_0 > 2$ and one has
\begin{equation*}
V = \sum_{(n,j)\in\bz\times\bz_+} \bra n \ket^{-\nu_0} e^{-\xi(j)} | \delta_n \otimes \e_j \rangle \langle \delta_n \otimes \e_j | = ( M_\psi \otimes K )^2,
\end{equation*}
where $K = \sum_{j\in\bz_+} e^{-\frac12 \xi(j)} |\e_j \rangle \langle \e_j|$. Thus, $V$ satisfies Assumption \ref{eq:hyppert} with $\mathscr V = I$ and 
\begin{equation*}
P_s \textbf{\textup{V}}_J P_s^\ast = P_s \textbf{\textup{V}} P_s^\ast = \Vert \psi \Vert_{\ell^2(\bz_h)}^2 \Big( I_{\bc} \otimes \sum_{j \ge j_s} e^{-\frac12 \xi(j)} |\e_j \rangle \langle \e_j| \Big),
\end{equation*}
where $(\e_{j})_{j \ge j_s}$ an orthonormal basis of ${\rm Ker}(Q - \mu_s)$ (if $Q \neq 0$). The non-zero eigenvalues of the operator $P_s \textbf{\textup{V}} P_s^\ast$ are simple and
\begin{equation}
\sigma(P_s \textbf{\textup{V}} P_s^\ast) = \bigr \{ \Vert \psi \Vert_{\ell^2(\bz_h)}^2 e^{-\frac12 \xi(j)} : j \ge j_s \bigl \} \cup \{0\}.
\end{equation}
Therefore, for $r > 0$ small enough, one has
\begin{equation}
\label{eq:decexp}
\begin{split}
\Phi_2(r) & = \Phi_1(r) = {\rm Tr} \: \mathds 1_{(r,+\infty)} (P_s \textbf{\textup{V}} P_s^\ast) \\
& = \# \bigr \{ \Vert \psi \Vert_{\ell^2(\bz_h)}^2 e^{-\frac12 \xi(j)} : j \ge j_s : r < \Vert \psi \Vert_{\ell^2(\bz_h)}^2 e^{-\frac12 \xi(j)} \bigl \} \\
& = \# \bigr \{ j \ge j_s : \xi(j) < 2 \ln \big( \Vert \psi \Vert_{\ell^2(\bz_h)}^2 r^{-1} \big)  \bigl \}. 
\end{split}
\end{equation}
The claim follows from \eqref{eq:decexp} and \eqref{eq:fctg}.
\end{proof}

\section{Decomposition of the potential}\label{sec6}

The aim of this section is to prove the next lemma that gives a suitable factorization of the perturbation $V$ satisfying Assumption \ref{eq:hyppert0}. In particular, this justifies our choice of the generalized Assumption \ref{eq:hyppert}.

\begin{Lemma}
\label{l:sdecV}
Let Assumption \ref{eq:hyppert0} holds, $M_\psi$ and $p$ be defined by \eqref{def:opn} and \eqref{def:opp} respectively. Then:
\begin{enumerate}
\item[{\rm i)}] There exists $\mathscr V \in \mathcal B (\ell^2(\bz_h,\mathcal{G}))$ such that 
\begin{equation}
\label{dec-V}
V = (M_\psi \otimes p) {\mathscr V} (M_\psi \otimes p).
\end{equation}
\item[{\rm ii)}] $\mathscr V \ge 0$ so that
\begin{equation}
V = {\mathscr M}^\ast {\mathscr M} \quad with \quad {\mathscr M} := \mathscr V^{1/2} (M_\psi \otimes p).
\end{equation}
In particular $V$ is trace class and 
\begin{equation}
\label{ntrV}
\Vert V \Vert_{\mathfrak S_1} \le \Vert {\mathscr M} \Vert_{\mathfrak S_2}^2.
\end{equation}
\item[{\rm iii)}] As a matrix, $\mathscr V = \{ a(n,m) \}_{(n,m) \in \bz^2}$ with 
\begin{equation}
\label{matrcalV}
\Vert a(n,m)\Vert_{\mathcal B(\mathcal{G})} \leq C \bra n \ket^{-\nu_1+\nu_0/2} \bra m \ket^{-\nu_2+\nu_0/2}. 
\end{equation}
\end{enumerate}
\end{Lemma}

\begin{proof}
i) Constants are generic, i.e. change from an estimate to another. One can write
\begin{equation}
\label{rep:V+}
V = (M_\psi \otimes I_{\mathcal G}) \widetilde V (M_\psi \otimes I_{\mathcal G}), 
\end{equation}
with $\widetilde V := \big\{ \bra n \ket^{\nu_0/2} V_h(n,m) \bra m \ket^{\nu_0/2} \big\}_{(n,m) \in \bz^2}$. Namely, one has
\begin{equation}
\label{rep:V}
\widetilde V = \sum_{(n,m)\in\bz^2} | \delta_n \rangle \langle \delta_m | \otimes \bra n \ket^{\nu_0/2} V_h(n,m) \bra m \ket^{\nu_0/2}.
\end{equation}

a) Firstly, let us show that the operator $\widetilde V$ is bounded. To see this, notice that Assumption \ref{eq:hyppert0} and \eqref{eq:g+9,0} imply that there exits a constant $C$ such that for each $(n,m) \in \bz^2$,
\begin{equation}
\label{hyp:decV}
\Vert V_h(n,m)\Vert_{\mathcal B(\mathcal{G})} \leq C \bra n \ket^{-\nu_1} \bra m \ket^{-\nu_2}.
\end{equation}
Using \eqref{rep:V} and the Cauchy-Schwartz inequality, one gets for any $\varphi \in \ell^2(\bz_h,\mathcal{G})$
\begin{equation}
\label{eq:ac}
\begin{split}
\Vert \widetilde V \varphi \Vert^2 & = \sum_{n\in\bz} \big\Vert \sum_{m \in \bz} \bra n \ket^{\nu_0/2} V_h(n,m) \bra m \ket^{\nu_0/2} \varphi(hm) \big\Vert_{\mathcal{G}}^2 \\
& \le \sum_{(n,m)\in\bz^2}  \| \bra n \ket^{\nu_0/2} V_h(n,m) \bra m \ket^{\nu_0/2} \|_{\mathcal B(\mathcal{G})}^2 \Vert \varphi \Vert^2.
\end{split}
\end{equation}
It follows from \eqref{hyp:decV} that
\begin{center}$
\displaystyle
\| \widetilde V \|^2 \le \sum_{(n,m)\in\bz^2} \| \bra n \ket^{\nu_0/2} V_h(n,m) \bra m \ket^{\nu_0/2} \|_{\mathcal B(\mathcal{G})}^2 \le
C \sum_{(n,m)\in\bz^2} \bra n \ket^{-2\nu_1+\nu_0} \bra m \ket^{-2\nu_2+\nu_0} < \infty.
$\end{center}

b) For $(n,m) \in \bz^2$, define in $\mathcal{G}$ the operator 
\begin{equation} 
\label{rep:v}
\widetilde V_h(n,m) := \sum_{(j,k)\in\Lambda^2} \bra j \ket^{\beta_0/2} v_{jk}^h(n,m) \bra k \ket^{\beta_0/2} | \e_j \rangle \langle \e_k |,
\end{equation}
which is bounded. Indeed, using the Cauchy-Schwartz inequality one gets for each $q \in \mathcal{G}$
\begin{align*}
& \Vert \widetilde V_h(n,m) q \Vert_{\mathcal{G}}^2 = \sum_{j \in \Lambda} \big\vert \bra \e_j,\widetilde V_h(n,m) q \ket_{\mathcal{G}} \big\vert^2 = \sum_{j \in \Lambda} \big\vert \sum_{k \in \Lambda} \bra j \ket^{\beta_0/2} v_{jk}^h(n,m) \bra k \ket^{\beta_0/2} \bra \e_k,q \ket_{\mathcal{G}} \big\vert^2 \\
& \le \sum_{j\in\Lambda} \Big( \sum_{k\in\Lambda} \vert \bra j \ket^{\beta_0/2} v_{jk}^h(n,m) \bra k \ket^{\beta_0/2} \vert^2 \Big) \sum_{k \in \Lambda} \vert \bra \e_k,q \ket_{\mathcal{G}} \vert^2 = \Vert q \Vert_{\mathcal{G}}^2 \sum_{(j,k)\in\Lambda^2} \bra j \ket^{\beta_0} \vert v_{jk}^h(n,m) \vert^2 \bra k \ket^{\beta_0} \\
& \le C \bra n \ket^{-2\nu_1} \bra m \ket^{-2\nu_2} \Vert q \Vert_{\mathcal{G}}^2 \sum_{(j,k)\in\Lambda^2} \bra j \ket^{\beta_0} G_1^2(i,j) \bra k \ket^{\beta_0} \le C \bra n \ket^{-2\nu_1} \bra m \ket^{-2\nu_2} \Vert q \Vert_{\mathcal{G}}^2.
\end{align*}
It follows that 
\begin{equation}
\label{hyp:dectildV}
\Vert \widetilde V_h(n,m)\Vert_{\mathcal B(\mathcal{G})} \leq C \bra n \ket^{-\nu_1} \bra m \ket^{-\nu_2}.
\end{equation}
Now, for every $k_\ast \in \Lambda$, one has
\begin{align*}
& \big( p \widetilde V_h(n,m) p \big) \e_{k_\ast} = \sum_{j \in \Lambda} \bra j \ket^{-\beta_0/2} | \e_j \rangle \bra \e_j,\widetilde V_h(n,m) p \e_{k_\ast} \ket_{\mathcal{G}} \\
& = \sum_{j \in \Lambda} \bra j \ket^{-\beta_0/2} \bra k_\ast \ket^{-\beta_0/2} | \e_j \rangle \bra \e_j,\widetilde V_h(n,m) \e_{k_\ast} \ket_{\mathcal{G}}.
\end{align*}
Since 
\begin{align*}
\bra \e_j,\widetilde V_h(n,m) \e_{k_\ast} \ket_{\mathcal{G}} & = \sum_{(j',k)\in\Lambda^2}  \bra j' \ket^{\beta_0/2} v_{j'k}^h(n,m) \bra k \ket^{\beta_0/2} \langle \e_j,\e_{j'} \rangle_{\mathcal{G}} \langle \e_k,\e_{k_\ast} \rangle_{\mathcal{G}} \\
& = \bra j \ket^{\beta_0/2} v_{jk_\ast}^h(n,m) \bra k_\ast \ket^{\beta_0/2},
\end{align*}
then we have
$$
\big( p \widetilde V_h(n,m) p \big) \e_{k_\ast} = \sum_{j \in \Lambda} v_{jk_\ast}^h(n,m) | \e_j \rangle = V_h(n,m) \e_{k_\ast},
$$
Therefore, for each $(n,m) \in \bz^2$, one has 
\begin{equation}
\label{eq:v}
p \widetilde V_h(n,m) p = V_h(n,m). 
\end{equation}
Together with \eqref{rep:V}, this implies that
\begin{equation}
\label{rep:tildV}
\widetilde V = \sum_{(n,m)\in\bz^2} | \delta_n \rangle \langle \delta_m | \otimes \bra n \ket^{\nu_0/2} p \widetilde V_h(n,m) p \bra m \ket^{\nu_0/2} = (I_{\ell^2(\bz)} \otimes p) \mathscr V (I_{\ell^2(\bz)} \otimes p),
\end{equation}
where
\begin{equation}
\label{rep:calV}
\mathscr V = \sum_{(n,m)\in\bz^2} | \delta_n \rangle \langle \delta_m | \otimes \bra n \ket^{\nu_0/2} \widetilde V_h(n,m) \bra m \ket^{\nu_0/2}.
\end{equation}
By putting together \eqref{rep:V+} and \eqref{rep:tildV}, one obtains \eqref{dec-V}.
Using \eqref{hyp:dectildV} and arguing as in \eqref{eq:ac}, one can show that $\mathscr V$ is a bounded operator.

\medskip

ii) Let us show that $\mathscr V \ge 0$ if $V \ge 0$. First, notice that the vectors of the basis $(\delta_n \otimes e_j)_{(n,j)\in {\mathbb Z}\times \Lambda}$ of $\ell^2(\bz_h,\mathcal{G})$ are eigenvectors of the operator $M_\psi \otimes p$. Indeed, one has
$$
(M_\psi \otimes p)(\delta_n \otimes \e_j) = \bra n \ket^{-\nu_0/2} \bra j \ket^{-\nu_0/2} (\delta_n \otimes \e_j), \qquad (n,j)\in {\mathbb Z}\times \Lambda.
$$
This implies that the Range of $M_\psi \otimes p$ is dense in $\ell^2(\bz_h,\mathcal{G})$. Then, for each $\varphi \in \ell^2(\bz_h,\mathcal{G})$, there exists a sequence of vectors $\varphi_q \in {\rm Range}(M_\psi \otimes p)$, $q \ge 0$, such that 
\begin{equation}
\label{eq:lim}
\lim_{q \to +\infty} \Vert \varphi_q - \varphi \Vert = 0.
\end{equation}
Hence, for each $q \ge 0$, $\varphi_q = (M_\psi \otimes p) \varphi_q'$ for some $\varphi_q' \in \ell^2(\bz_h,\mathcal{G})$. Noting that $M_\psi \otimes p$ is a positive operator and using \eqref{dec-V}, it follows that 
\begin{equation}
\label{eq:lim+}
\begin{split}
\bra \mathscr V \varphi_q,\varphi_q \ket & = 
\bra \mathscr V (M_\psi \otimes p) \varphi_q',(M_\psi \otimes p) \varphi_q'\ket \\
& = \bra (M_\psi \otimes p) \mathscr V (M_\psi \otimes p) \varphi_q', \varphi_q' \ket = \bra V \varphi_q', \varphi_q' \ket \ge 0.
\end{split}
\end{equation}
Now, since for $q \ge 0$
\begin{align*}
\vert \bra \mathscr V \varphi_q,\varphi_q \ket - \bra \mathscr V \varphi,\varphi \ket \vert & = \vert \bra \mathscr V (\varphi_q - \varphi),\varphi_q \ket + \bra \mathscr V \varphi,\varphi_q - \varphi \ket \vert \\
& \le \Vert \mathscr V \Vert \Vert \varphi_q - \varphi \Vert \Vert \varphi_q \Vert + \Vert \mathscr V \Vert \Vert \varphi \Vert \Vert \varphi_q - \varphi \Vert,
\end{align*}
one deduces from \eqref{eq:lim} and \eqref{eq:lim+} that
$$
\bra \mathscr V \varphi,\varphi \ket = \lim_{q \to +\infty} \bra \mathscr V \varphi_q,\varphi_q \ket \ge 0.
$$
Moreover, since the operator $M_\psi \otimes p$ is Hilbert-Schmidt according to Remark \ref{rem:perV}, then $V$ is trace class and \eqref{ntrV} follows by \eqref{dec-V} and the boundedness of $\mathscr V$.

\medskip

iii) As matrix $\mathscr V = \{ a(n,m) \}_{(n,m) \in \bz^2}$, one has from \eqref{rep:calV} that for each $(n,m) \in \bz^2$
$
a(n,m) = \bra n \ket^{\nu_0/2} \widetilde V_h(n,m) \bra m \ket^{\nu_0/2}.
$
Then, \eqref{matrcalV} follows immediately from \eqref{hyp:dectildV}.
\end{proof}

\section{Preliminary results}\label{sec7}

\subsection{Extensions of the kernel of $(-\Delta_h-z)^{-1}$ to the real axis}

For further references, we provide more details on our choice of analytic determinations of the complex logarithm and square-root functions. 
First, it can be checked that the map ${\exp} : s \in \bc \mapsto e^s \in \bc^\ast$ is a surjective group homomorphism with kernel 
$\ker (\exp) = 2i\pi \bz$. 
It follows that its restriction ${\exp} : s \in \{s \in \bc : -\pi \le {\rm Im} \, s < \pi \} \mapsto e^s \in \bc^\ast$ is a bijective map. Since the image of the axis $\{s \in \bc : {\rm Im} \, s = -\pi \}$ is the semi-axis $(-\infty,0)$, then that of the (open) domain $\{s \in \bc : -\pi < {\rm Im} \, s < \pi \}$ is the domain $\bc \smallsetminus (-\infty,0]$ so that
$$
\exp :  s \in \{s \in \bc : -\pi < {\rm Im} \, s < \pi \} \mapsto e^s \in \bc \smallsetminus (-\infty,0]
$$
is a holomorphic bijective map with non-vanishing derivative. 
The corresponding inverse map 
\begin{equation}
\label{det:ln}
{\rm Ln} :  s \in \bc \smallsetminus (-\infty,0] \mapsto {\rm Ln}(s) \in \{s \in \bc : -\pi < {\rm Im} \, s < \pi \}
\end{equation}
is then holomorphic and will define our complex logarithm determination. It corresponds to the principal value of the logarithm function. Hence, one can define the complex analytic square-root determination using the analytic function ${\rm Ln}$ by 
\begin{equation}
\label{det:rc}
\sqrt{\cdot} = e^{\frac12{\rm Ln}} :  s \in \bc \smallsetminus (-\infty,0] \mapsto e^{\frac12{\rm Ln}(s)} \in \{s \in \bc : {\rm Re} \, s > 0 \}.
\end{equation}
It corresponds to the principal value of the square-root function. Note that \eqref{det:ln} and \eqref{det:rc} correspond to employ the principal value of the argument ${\rm Arg}$ which takes values in $(-\pi,\pi]$ so that 
\begin{equation}
{\rm Ln}(s) = {\rm Ln} \vert s \vert + i {\rm Arg}(s).
\end{equation}

Let $z \in \bc \smallsetminus [0,\frac4{h^2}]$ and $R(z) := (-\Delta_h - z)^{-1}$ be the resolvent of $-\Delta_h$. One has
\begin{equation}
\label{eq:fourconv}
\mathscr F \big( R(z) \phi \big)(\theta) = \frac{(\mathscr F \phi)(\theta)}{f(\theta)-z},
\end{equation}
where $f(\theta)$ is given by \eqref{eq:fonctf}. It can be checked that 
\begin{equation}
\label{eq:fourconv+}
\mathscr F (h^2e^{i \alpha \vert \cdot \vert})(\theta) = -\frac{i}{2h^{-2}} \frac{\sin(h\alpha)}{\sin^2(h\theta/2) - \sin^2(h\alpha/2)}, \quad {\rm Im} \,\alpha > 0.
\end{equation}
It follows from identities \eqref{eq:fourconv} and \eqref{eq:fourconv+} that for ${\rm Im} \, z > 0$, $R(z)$ is an operator with convolution kernel given by $r(z,h(n-m))$,
where
\begin{equation}\label{KER}
r(z,hk) = \frac{ih^2}2 \frac{e^{i\alpha(z)h\vert k \vert}}{\sin(h\alpha(z))} = \frac{ie^{2i\vert k \vert {\rm Arcsin}(\frac{h}{2} \sqrt {z}) } } {\sqrt z \sqrt{4/{h^2} -z}} := R(z,k), \quad k \in \mathbb Z.
\end{equation}
 Here, $\alpha(z) = \frac2h{\rm Arcsin}(\frac{h}{2}\sqrt{z})$ is the unique solution to the equation $\frac{2  - 2\cos(h\alpha)}{h^2} = \frac{4}{h^2} \sin^2 \big( \frac{h\alpha}{2} \big) = z$ lying in the region $\{\alpha \in \mathbb C : -\frac{\pi}h \le {\rm Re} \, \alpha \leq \frac{\pi}h : {\rm Im} \, \alpha > 0\}$, where ${\rm Arcsin}$ is principal value of the real arcsine ($\arcsin$) function obtained by employing the above analytic determinations ${\rm Ln}$ and $\sqrt{\cdot} = e^{\frac12{\rm Ln}}$. Namely, for $s \in \bc \smallsetminus \big( (-\infty,-1] \cup [1,+\infty) \big)$, one has
 \begin{equation}
\label{det:arcs}
 {\rm Arcsin} \, s = \frac1i {\rm Ln} \big( is + \sqrt{1-s^2} \big) = w.
 \end{equation}
 It can be easily checked that $w$ given by \eqref{det:arcs} is solution to the equation $\sin(w) = \frac{e^{iw} - e^{-iw}}{2i} = s$. In particular, if $s = x \in \br$ with $\vert x \vert < 1$, then 
$
\label{det:arc}
 {\rm Arcsin} \, x = {\rm Arg} \big( ix + \sqrt{1-x^2} \big) \in (-\tfrac{\pi}2,\tfrac{\pi}2)
 $
 coincides with the real classical $\arcsin$ inverse fonction, i.e. ${\rm Arcsin} \, x = \arcsin x$. 
 
 The next result follows immediately taking into account the above considerations. 
 
 \begin{Proposition}
 \label{prop:extnoyau} 
 One has 
 \begin{equation}
 \label{limit:noy}
 \lim_{\varepsilon \searrow 0} R(\lambda + i\varepsilon,n-m) = 
 \begin{cases}
 \frac{e^{2i\vert n-m \vert {\rm Arcsin}(\frac{ih}{2}\sqrt{-\lambda} ) } } {\sqrt{-\lambda} \sqrt{4/{h^2}-\lambda}} \quad {\rm if} \quad \lambda < 0 \\
 \frac{ie^{2i\vert n-m \vert {\rm Arcsin}(\frac{h}{2}\sqrt{\lambda} ) } } {\sqrt{\lambda} \sqrt{4/{h^2}-\lambda}} \quad {\rm if} \quad \lambda \in (0,4/{h^2}).
 \end{cases}
 \end{equation}
 \end{Proposition} 


For $\lambda \in \br \smallsetminus \big( \{0 \} \cup [\frac4{h^2},+\infty) \big)$, one defines $R(\lambda)$ as the operator acting in $\ell^2(\bz_h)$ with the convolution kernel $R(\lambda,n-m)$ where
 \begin{equation}
 \label{limit:noy}
 R(\lambda,n-m) := \lim_{\varepsilon \searrow 0} R(\lambda + i\varepsilon,n-m).
  \end{equation}

\subsection{Estimates of the weighted resolvents}\label{ss:est}

\subsubsection{Hilbert-Schmidt bounds}

Let $T(z)$ be the weighted resolvent defined by \eqref{def:Tpond}. For ${\rm \-Im} \, z > 0$, thanks to \eqref{KER}, one defines $R(z)$ as the operator acting in $\ell^2(\bz_h)$ with convolution kernel $R(z,n-m)$. So, according to \eqref{limit:noy}, one extends $T(z)$ to $\bc \smallsetminus \mathscr E_Q$ by setting 
\begin{equation}
 \label{eq:Tdaint}
T_s(\lambda) :=  \mathscr V^{1/2} ( M_\psi R(\lambda - \mu_s) M_\psi \otimes K \pi_s K^\ast ) \mathscr V^{1/2}, \quad \lambda \in (\mu_s,\tfrac{4}{h^2} + \mu_s), \quad 0 \le s \le d,
\end{equation}
and
 \begin{equation}
 \label{eq:Tdaint-}
 \begin{split}
T(\lambda) := \sum_{s \, : \, \lambda - \mu_s \in (0,\tfrac{4}{h^2})} T_s(\lambda) +
\sum_{s' \, : \, \lambda - \mu_{s'} \notin [0,\tfrac{4}{h^2}]} \mathscr V^{1/2} ( M_\psi R(\lambda - \mu_{s'}) M_\psi \otimes K \pi_{s'} K^\ast ) \mathscr V^{1/2},
 \end{split}
 \end{equation}
 $\lambda \in \sigma(H_Q) \smallsetminus \mathscr E_Q
 $. Introduce $\bc^+ := \{ z \in \bc : {\rm Im} \, z > 0 \}$ and $\overline{\bc^+} := \{ z \in \bc : {\rm Im} \, z \ge 0 \}$. 
 
 \begin{Proposition}
 \label{prop:Ts2}
 Let $V$ satisfy Assumption \ref{eq:hyppert}. Then:
 \begin{enumerate}
 \item[{\rm i)}] For any $z \in \bc \smallsetminus \mathscr E_Q$, the operator $T(z) \in {\mathfrak S_2}(\ell^2(\bz_h,\mathcal{G}))$. Furthermore, if $Q = 0$, then
\begin{equation}
\label{estS2:Q0}
 \Vert T(\lambda) \Vert_{\mathfrak S_2} \le \dfrac{ \Vert K \Vert_{\mathfrak S_2}^2 \Vert \mathscr V \Vert}{\lambda^{1/2} (\frac4{h^2}-\lambda)^{1/2}} \sum_{n\in\bz} \bra n \ket^{-\nu_0}, \qquad \lambda \in (0,\tfrac{4}{h^2}),
\end{equation}
and for $Q \neq 0$ with $\lambda \in \sigma(H_Q) \smallsetminus \mathscr E_Q$, we have
\begin{equation}
\label{estS2:Qg}
 \begin{split}
& \Vert T(\lambda) \Vert_{\mathfrak S_2} \le \Vert K \Vert_{\mathfrak S_2}^2 \Vert \mathscr V \Vert  \sum_{n\in\bz} \bra n \ket^{-\nu_0}  \\
& \cdot \Bigg( \sum_{s \, : \, \lambda - \mu_s \in (0,\tfrac{4}{h^2})} \dfrac{1} {(\lambda - \mu_s)^{1/2} (\frac4{h^2}-\lambda + \mu_s)^{1/2}} + \sum_{s' \, : \, \lambda - \mu_{s'} \notin [0,\tfrac{4}{h^2}]} \dfrac{1} {{\rm dist}(\lambda - \mu_{s'},[0,\tfrac{4}{h^2}])} \Bigg).
\end{split}
\end{equation}
\item[{\rm ii)}] 
The operator-valued function $z \in \overline{\bc^+} \smallsetminus \mathscr E_Q \mapsto T(z) \in \mathfrak S_2$ is continuous.
\end{enumerate}
 \end{Proposition}
 
 \begin{proof}
 i) Let $z \in \bc \smallsetminus \sigma(H_Q)$. Then, $(H_Q-z)^{-1}$ is bounded with 
 $\Vert (H_Q-z)^{-1} \Vert \le \frac{1}{{\rm dist}(z,\sigma(H_Q))}$. 
 By Remark \ref{rem:perV} $M_\psi \in \mathfrak S_2$ and by Assumption \ref{eq:hyppert} $K \in \mathfrak S_2$. Then $T(z) \in \mathfrak S_2$ with 
  \begin{align*}
 \Vert T(z) \Vert_{\mathfrak S_2} \le \Vert {\mathscr M} \Vert_{\mathfrak S_2}^2 \Vert (H_Q-z)^{-1} \Vert 
 &\le \dfrac{\Vert M_\psi \Vert_{\mathfrak S_2}^2 \Vert K \Vert_{\mathfrak S_2}^2 \Vert \mathscr V \Vert}{{\rm dist}(z,\sigma(H_Q))} 
 = \dfrac{\Vert K \Vert_{\mathfrak S_2}^2 \Vert \mathscr V \Vert}{{\rm dist}(z,\sigma(H_Q))}\displaystyle \sum_{n\in\bz} \bra n \ket^{-\nu_0}.
  \end{align*}
 
 Let us show \eqref{estS2:Q0} assuming $Q=0$ and $\lambda \in (0,\frac4{h^2})$. The operator $R(\lambda)$ admits the convolution kernel $R(\lambda,n-m)$ given by \eqref{limit:noy}. Using Proposition \ref{prop:extnoyau}, 
 we get 
 \begin{center}$
 \Vert M_\psi R(\lambda) M_\psi \Vert_{\mathfrak S_2} \le \dfrac{\sum_{n\in\bz} \bra n \ket^{-\nu_0}}{\lambda^{1/2} (\frac4{h^2}-\lambda)^{1/2}}.
 $ \end{center}
 This together with \eqref{eq:Tdaint} implies that
  \begin{center}$
 \Vert T(\lambda) \Vert_{\mathfrak S_2} \le \dfrac{\Vert K \Vert_{\mathfrak S_2}^2 \Vert \mathscr V \Vert} {\lambda^{1/2} (\frac4{h^2}-\lambda)^{1/2}} \displaystyle \sum_{n\in\bz} \bra n \ket^{-\nu_0} \Vert.
 $ \end{center}
 
Now, let us show \eqref{estS2:Qg} assuming $Q \neq 0$ and $\lambda \in \sigma(H_Q) \smallsetminus \mathscr E_Q$. Using \eqref{eq:Tdaint-}, one treats the two sums of the r.h.s and by arguing as above and noting that $\Vert K \pi_s K^\ast \Vert_{\mathfrak S_2}^2 \le \Vert K \Vert_{\mathfrak S_2}^2$, one obtains
\begin{equation*}
\begin{split}
 \Vert & T(\lambda) \Vert_{\mathfrak S_2} \le \Vert K \Vert_{\mathfrak S_2}^2 \Vert \mathscr V \Vert \sum_{n\in\bz} \bra n \ket^{-\nu_0}  \\
& \cdot \Bigg( \sum_{s \, : \, \lambda - \mu_s \in (0,\tfrac{4}{h^2})} \dfrac{1} {(\lambda - \mu_s)^{1/2} (\frac4{h^2}-\lambda + \mu_s)^{1/2}} + \sum_{s \, : \, \lambda - \mu_s \notin [0,\tfrac{4}{h^2}]} \dfrac{1} {{\rm dist}(\lambda - \mu_s,[0,\tfrac{4}{h^2}])} \Bigg).
\end{split}
\end{equation*}
 
 ii) According to the point i), the map $z \in \overline{\bc^+} \smallsetminus \mathscr E_Q \mapsto T(z) \in \mathfrak S_2$ is well defined. Otherwise, since the map $z \mapsto (H_Q-z)^{-1}$ is holomorphic in $\bc \smallsetminus \sigma(H_Q)$, then the continuity of $z \in \overline{\bc^+} \smallsetminus \sigma(H_Q) \mapsto T(z) \in \mathfrak S_2$ follows immediately. Indeed, as $\vert z - z_0 \vert \to 0$ with $z$, $z_0 \in \overline{\bc^+} \smallsetminus \sigma(H_Q)$, 
 \begin{align*}
\Vert T(z) - T(z_0) \Vert_{\mathfrak S_2} & = \Vert {\mathscr M} (H_Q - z)^{-1} {\mathscr M}^\ast - {\mathscr M} (H_Q - z_0)^{-1} {\mathscr M}^\ast \Vert_{\mathfrak S_2} \\
& \le \Vert {\mathscr M} \Vert_{\mathfrak S_2}^2 \Vert (H_Q - z)^{-1} - (H_Q - z_0)^{-1} \Vert \to 0.
 \end{align*}

Now, we focus on the continuity in the set $\sigma(H_Q) \smallsetminus \mathscr E_Q$ by treating first the case $Q=0$ so that $\mathscr E_0 = \{0,\frac4{h^2} \}$. Let $z_0 = \lambda_0 \in (0,\frac4{h^2})$ and $0 < \delta \ll 1$. Then, for 
 $
z \in D_\delta(\lambda_0) := \{ z \in \overline{\bc^+} \smallsetminus \mathscr E_0 : \vert z - \lambda_0 \vert \le \delta \},
$ 
one has
\begin{equation}
\label{est:contT}
 \Vert T(z) - T(\lambda_0) \Vert_{\mathfrak S_2} \le \Vert \mathscr V \Vert \Vert K \Vert_{\mathfrak S_2}^2 \Vert M_\psi \big[ R(z)-R(\lambda_0) \big] M_\psi \Vert_{\mathfrak S_2}.
 \end{equation}
 The operator $M_\psi \big[ R(z)-R(\lambda_0) \big] M_\psi$ has the kernel
  $
 \bra n \ket^{-\nu_0/2} \big( R(z,n-m)-R(\lambda_0,n-m) \big) \bra m \ket^{-\nu_0/2},$ 
 so that
 \begin{equation}
 \begin{split}
 \label{est:contT+}
 \Vert M_\psi \big[ R(z)-R(\lambda_0) \big] M_\psi \Vert_{\mathfrak S_2}^2 
 \le \sum_{(n,m) \in \bz^2} \bra n \ket^{-\nu_0} \bra m \ket^{-\nu_0} \vert R(z,n-m)-R(\lambda_0,n-m) \vert^2.
 \end{split}
 \end{equation}
 The map $z \in D_\delta(\lambda_0) \mapsto \frac{1}{\vert z \vert^{1/2} \vert \frac4{h^2}-z \vert^{1/2}} \in \br$ is continuous. Since $D_\delta(\lambda_0)$ is compact, then there exists $a_0 \in D_\delta(\lambda_0)$ such that
 $$
 \sup_{z \in D_\delta(\lambda_0) \atop (n,m) \in \bz^2}  \vert R(z,n-m)-R(\lambda_0,n-m) \vert 
 \le \frac{1}{\vert a_0 \vert^{1/2} \vert \frac4{h^2}-a_0 \vert^{1/2}} + \frac{1}{\lambda_0^{1/2} (\frac4{h^2}-\lambda_0)^{1/2}} =: C(a_0,\lambda_0).
 $$
 That is the map $(n,m) \in \bz^2 \mapsto \bra n \ket^{-\nu_0} \bra m \ket^{-\nu_0} \vert R(z,n-m)-R(\lambda_0,n-m) \vert$ is uniformly dominated w.r.t $z \in D_\delta(\lambda_0)$ by the map $(n,m) \in \bz^2 \mapsto C (a_0,\lambda_0) \bra n \ket^{-\nu_0} \bra m \ket^{-\nu_0}$. Now, using
 $$
 \sum_{(n,m) \in \bz^2} \bra n \ket^{-\nu_0} \bra m \ket^{-\nu_0} = \Big( \sum_{n\in\bz} \bra n \ket^{-\nu_0} \Big)^2 < \infty,
 $$
\eqref{est:contT+}, Lebesgue's dominated convergence theorem and \eqref{est:contT}, one gets 
 $\Vert T(z) - T(\lambda_0) \Vert_{\mathfrak S_2} \to 0$ as $\vert z - \lambda_0 \vert \to 0.$
This shows that $\lambda \in (0,\tfrac{4}{h^2}) \mapsto T(\lambda) \in \mathfrak S_2$ is continuous. 

For the case $Q \neq 0$, let us show that the map $\lambda \in \sigma(H_Q) \smallsetminus \mathscr E_Q \mapsto T(\lambda) \in \mathfrak S_2$ is continuous. So, for $\lambda_0 \in \sigma(H_Q) \smallsetminus \mathscr E_Q = \bigcup_{s=1}^d (\mu_s,\frac{4}{h^2} + \mu_s)$ fixed, we want to prove that
$\Vert T(\lambda) - T(\lambda_0) \Vert_{\mathfrak S_2} \to 0$ as $\vert \lambda - \lambda_0 \vert \to 0$. Using \eqref{eq:Tdaint-}, write for $\lambda \in \sigma(H_Q) \smallsetminus \mathscr E_Q$ and sufficiently close to $\lambda_0$
\begin{equation}
 \label{eq:Tdiff}
 \begin{split}
T(\lambda) - T(\lambda_0) & = \sum_{s \, : \, \lambda_0 - \mu_s \in (0,\tfrac{4}{h^2})} \mathscr V^{1/2} ( M_\psi [ R(\lambda - \mu_{s}) - R(\lambda_0 - \mu_{s}) ]M_\psi \otimes K \pi_s K^\ast ) \mathscr V^{1/2} \\
& + \sum_{s' \, : \, \lambda_0 - \mu_{s'} \notin [0,\tfrac{4}{h^2}]} \mathscr V^{1/2} ( M_\psi [ R(\lambda - \mu_{s'})  - R(\lambda_0 - \mu_{s'}) ] M_\psi \otimes K \pi_s K^\ast ) \mathscr V^{1/2}.
 \end{split}
 \end{equation}
 In the r.h.s of \eqref{eq:Tdiff}, $\lambda$ is close to $\lambda_0$ so that $\lambda \in (\mu_s,\tfrac{4}{h^2} + \mu_s)$ for all $s$ and $\lambda \notin (\mu_{s'},\tfrac{4}{h^2} + \mu_{s'})$ for all $s'$. By arguing as in the case $Q=0$, one can prove that the map 
\begin{center}
$z \in (0,\tfrac{4}{h^2}) \mapsto \mathscr V^{1/2} ( M_\psi R(z) M_\psi \otimes K \pi_s K^\ast ) \mathscr V^{1/2} \in \mathfrak S_2$
\end{center}
is continuous. Moreover, clearly, the map
\begin{center}
$z \in \br \smallsetminus [0,\tfrac{4}{h^2}] \mapsto \mathscr V^{1/2} ( M_\psi R(z) M_\psi \otimes K \pi_s K^\ast ) \mathscr V^{1/2} \in \mathfrak S_2$
\end{center}
is continuous. This together with \eqref{eq:Tdiff} implies that
\begin{equation}
 \label{eq:Tdiff+}
 \begin{split}
& \Vert T(\lambda) - T(\lambda_0) \Vert_{\mathfrak S_2} \le \\ & \sum_{s \, : \, \lambda_0 - \mu_s \in (0,\tfrac{4}{h^2})} \Big\Vert \mathscr V^{1/2} ( M_\psi [ R(\lambda - \mu_{s}) - R(\lambda_0 - \mu_{s}) ]M_\psi \otimes K \pi_s K^\ast ) \mathscr V^{1/2} \Big\Vert_{\mathfrak S_2} \\
& + \sum_{s' \, : \, \lambda_0 - \mu_{s'} \notin [0,\tfrac{4}{h^2}]} \Big\Vert \mathscr V^{1/2} ( M_\psi [ R(\lambda - \mu_{s'})  - R(\lambda_0 - \mu_{s'}) ] M_\psi \otimes K \pi_s K^\ast ) \mathscr V^{1/2} \Big\Vert_{\mathfrak S_2},
 \end{split}
 \end{equation}
tends to $0$ as $\lambda \to \lambda_0$.
 \end{proof}


\begin{Corollary}
\label{cor:exT}
Let $V$ satisfy Assumption \ref{eq:hyppert} and $\lambda \in \br \smallsetminus \mathscr E_Q$. Then, the limit $T(\lambda+i0)$ exists in $\mathfrak S_2$ with $T(\lambda + i0) = T(\lambda)$.
\end{Corollary}

\subsubsection{Trace class bounds}

We want to establish the existence of $T(\lambda + i0)$, $\lambda \in \br \smallsetminus \mathscr E_Q$, in the trace class $\mathfrak S_1$. However, the proof is less evident than the one of Corollary \ref{cor:exT} obtained directly from Proposition \ref{prop:Ts2}. 
The first step consists of establishing the following simple result, whose proof is 
similar to the one of Proposition \ref{prop:Ts2} in many points.

\begin{Proposition}
 \label{prop:TS1-}
 Let $V$ satisfy Assumption \ref{eq:hyppert}. Then:
 \begin{enumerate}
 \item[{\rm i)}] For any $z \in \bc \smallsetminus \sigma(H_Q)$, the operator $T(z) \in {\mathfrak S_1}(\ell^2(\bz_h,\mathcal{G}))$ with 
  \begin{center}$
 \Vert T(z) \Vert_{\mathfrak S_1} \le \dfrac{\Vert M_\psi \Vert_{\mathfrak S_2}^2 \Vert K \Vert_{\mathfrak S_2}^2 \Vert \mathscr V \Vert}{{\rm dist}(z,\sigma(H_Q) \smallsetminus \mathscr E_Q)}.
$ \end{center}
\item[{\rm ii)}] 
The operator-valued function $z \in \bc \smallsetminus \sigma(H_Q) \mapsto T(z) \in \mathfrak S_1$ is holomorphic.
\end{enumerate}
 \end{Proposition}
 
 \begin{proof}
 i) Let $z \in \bc \smallsetminus \sigma(H_Q)$. Since the operators $M_\psi$, $K$ are Hilbert-Schmidt and $(H_Q-z)^{-1}$ is bounded, then $T(z) \in \mathfrak S_1$ with 
  \begin{center}$
 \Vert T(z) \Vert_{\mathfrak S_1} \le \Vert {\mathscr M} \Vert_{\mathfrak S_2}^2 \Vert (H_Q-z)^{-1} \Vert \le \dfrac{\Vert M_\psi \Vert_{\mathfrak S_2}^2 \Vert K \Vert_{\mathfrak S_2}^2 \Vert \mathscr V \Vert}{{\rm dist}(z,\sigma(H_Q) \smallsetminus \mathscr E_Q)}.
 $ \end{center}

 ii) Thanks to the point i), the map $z \in \bc \smallsetminus \sigma(H_Q) \mapsto T(z) \in \mathfrak S_1$ is well defined. Moreover,
 $$
 \Vert T'(z) \Vert_{\mathfrak S_1} = \Vert \mathscr M (H_Q-z)^{-2} \mathscr M^\ast \Vert_{\mathfrak S_1} \le \Vert \mathscr M (H_0-z)^{-1} \Vert_{\mathfrak S_2} \Vert (H_Q-z)^{-1} \mathscr M^\ast \Vert_{\mathfrak S_2},
 $$
 and as $\vert z - z_0 \vert \to 0$ with $z$, $z_0 \in \bc \smallsetminus \sigma(H_Q)$, one has
 \begin{align*}
\Big\Vert \frac{T(z) - T(z_0)}{z-z_0} - T'(z_0) \Big\Vert_{\mathfrak S_1}  
\le \Vert {\mathscr M} \Vert_{\mathfrak S_2}^2 \Big\Vert \frac{(H_Q-z)^{-1} - (H_Q-z_0)^{-1}}{z-z_0} - (H_Q-z_0)^{-2} \Big\Vert \to 0.
 \end{align*}
Thus the claim follows.
 \end{proof}

The second step consists of treating the case $\lambda \in \sigma(H_Q) \smallsetminus \mathscr E_Q$ which is more delicate. We first need to find a suitable integral decomposition of
\begin{equation}
\label{eq:Tzr}
 T_s(z) = \mathscr V^{1/2} ( M_\psi R(z-\mu_s) M_\psi \otimes K \pi_s K^\ast ) \mathscr V^{1/2}, \quad z \in \bc^+,
 \end{equation}
 $0 \le s \le d$.
 To simplify the notations, we introduce the operators $a(\theta) : \ell^2(\bz_h) \to \bc$ defined by
\begin{equation}
\label{op:A}
a(\theta) := \frac{1}{\sqrt{\tau}} \big \bra e^{-i(\cdot)\theta} \psi \vert, \qquad \tau = \frac{2\pi}h,
\end{equation}
and $a(\theta)^\ast : \bc \to \ell^2(\bz_h)$ the rank one operator given by
\begin{equation}
a(\theta)^\ast \zeta = \frac{\zeta}{\sqrt{\tau}} e^{-i(\cdot)\theta} \psi.
\end{equation}

\begin{Proposition}
\label{prop:intrepT}
Let $V$ satisfy Assumption \ref{eq:hyppert}. Then, for $z \in \bc^+$ and $0 \le s \le d$, one has
\begin{equation}
\label{eq:intrepT}
T_s(z) = \int_0^{\pi/h} u_s(\theta)^\ast u_s(\theta) \frac{d \theta}{f(\theta) - z + \mu_s},
\end{equation}
where $f(\theta)$ is given by \eqref{eq:fonctf} and $u_s(\theta) : \ell^2(\bz_h,\mathcal{G}) \to \bc^2 \otimes \mathcal{G}$ is the operator defined by
\begin{equation}
\label{eq:ut1}
u_s(\theta) := (\mathscr A(\theta) \otimes \pi_sK^\ast) \mathscr V^{1/2}, \quad \theta \in \mathbb T.
\end{equation}
Here, $\mathscr A(\theta) : \ell^2(\bz_h) \to \bc^2$ is the operator defined by
\begin{equation}
\label{eq:d}
\mathscr A(\theta) \phi = \begin{pmatrix}
\frac{1}{\sqrt{\tau}} \big \bra e^{-i(\cdot)\theta} \psi,\phi \big \ket_{\ell^2(\bz_h)} \\
\frac{1}{\sqrt{\tau}} \big \bra e^{i(\cdot)\theta} \psi,\phi \big \ket_{\ell^2(\bz_h)}
\end{pmatrix} = \begin{pmatrix}
a(\theta) \phi \\
a(-\theta) \phi
\end{pmatrix},
\end{equation}
where the operator $a(\theta)$ is given by \eqref{op:A}.
\end{Proposition}

\begin{proof}
Consider ${\mathscr F} : \ell^2({\mathbb Z}_h) \rightarrow {\rm L}^2({\mathbb T})$ the discrete Fourier transform defined by \eqref{eq:Tfour}. 
For any $\varphi$, $\Phi \in \ell^2(\bz_h)$ and $z \in \bc^+$, one has
\begin{equation}
\label{eq:a}
\begin{split}
& \bra \Phi,M_\psi R(z - \mu_s)M_\psi \varphi \ket_{\ell^2(\bz_h)}  = \bra M_\psi \Phi,R(z - \mu_s)M_\psi \varphi \ket_{\ell^2(\bz_h)} \\
& =  \bra {\mathscr F}(M_\psi \Phi),{\mathscr F} R(z - \mu_s) {\mathscr F}^{-1} {\mathscr F}(M_\psi \varphi) \ket_{{\rm L}^2({\mathbb T})} \\
&  = \frac1\tau \int_{-\pi/h}^{\pi/h} \frac{1}{f(\theta)-z + \mu_s} {\mathscr F}(M_\psi \varphi)(\theta) \overline{{\mathscr F}(M_\psi \Phi)(\theta)}d\theta \\
& = \frac1\tau \int_0^{\pi/h} \frac{1}{f(\theta) - z + \mu_s} \Big( {\mathscr F}(M_\psi \varphi)(-\theta) \overline{{\mathscr F}(M_\psi \Phi)(-\theta)} + {\mathscr F}(M_\psi \varphi)(\theta) \overline{{\mathscr F}(M_\psi \Phi)(\theta)} \Big) d\theta.
\end{split}
\end{equation}
For $\theta \in {\mathbb T}$, one has 
\begin{equation}
\label{eq:b}
{\mathscr F}(M_\psi \varphi)(-\theta) = \sum_{n \, \in \, {\mathbb Z}} e^{ihn\theta} \bra n \ket^{-\nu_0/2} \varphi(hn) = \big \bra e^{-i(\cdot)\theta} \psi,\varphi \big \ket_{\ell^2(\bz_h)},
\end{equation}
and then
\begin{equation}
\label{eq:c}
\overline{{\mathscr F}(M_\psi \Phi)(-\theta)} = \overline{\big \bra e^{-i(\cdot)\theta} \psi,\Phi \big \ket_{\ell^2(\bz_h)}}.
\end{equation}
By putting together \eqref{eq:a}-\eqref{eq:c} and using \eqref{eq:d}, one gets
\begin{equation}
\label{eq:e}
\begin{split}
& \bra \Phi,M_\psi R(z - \mu_s)M_\psi \varphi \ket_{\ell^2(\bz_h)} \\
& = \int_0^{\pi/h} \frac{1}{f(\theta) - z + \mu_s} \big( a(\theta)\varphi \cdot \overline{a(\theta)\Phi} + a(-\theta)\varphi \cdot \overline{a(-\theta)\Phi} \big) d\theta \\
& = \int_0^{\pi/h} \bra \mathscr A(\theta) \Phi, \mathscr A(\theta) \varphi \ket_{\bc^2} \frac{d \theta}{f(\theta) - z + \mu_s} \\
& = \int_0^{\pi/h} \bra \Phi, \mathscr A(\theta)^\ast \mathscr A(\theta) \varphi \ket_{\ell^2(\bz_h)} \frac{d \theta}{f(\theta) - z + \mu_s}.
\end{split}
\end{equation}
It follows that $M_\psi R(z - \mu_s)M_\psi$ admits the integral representation
\begin{equation}
\label{eq:repR}
M_\psi R(z - \mu_s)M_\psi = \int_0^{\pi/h} \mathscr A(\theta)^\ast \mathscr A(\theta) \frac{d \theta}{f(\theta) - z + \mu_s},
\end{equation}
where the operator $\mathscr A(\theta)^\ast : \bc^2 \to \ell^2(\bz_h)$ is given by
\begin{equation}
\mathscr A(\theta)^\ast \begin{pmatrix}
\zeta_1 \\
\zeta_2
\end{pmatrix} = 
a(\theta)^\ast \zeta_1 + a(-\theta)^\ast \zeta_2,
\end{equation}
so that $\mathscr A(\theta)^\ast \mathscr A(\theta) : \ell^2(\bz_h) \to \ell^2(\bz_h)$ is the rank two operator
\begin{equation}
\label{eq:ut2}
\begin{split}
\mathscr A(\theta)^\ast \mathscr A(\theta) = a(\theta)^\ast a(\theta) + a(-\theta)^\ast a(-\theta) = \frac{1}{\tau} \Big( \big \vert e^{-i(\cdot)\theta} \psi \big \ket \big \bra e^{-i(\cdot)\theta} \psi \big \vert + \big \vert e^{i(\cdot)\theta} \psi \big \ket \big \bra e^{i(\cdot)\theta} \psi \big \vert \Big).
\end{split}
\end{equation}
Then, by combining \eqref{eq:Tzr} and \eqref{eq:repR}, the integral representation \eqref{eq:intrepT} of $T_s(z)$ holds. 
\end{proof}

By performing the change of variable $\zeta = f(\theta)=\frac4{h^2} \sin^2(\frac{\theta h}2)$ in \eqref{eq:intrepT}, one gets 
\begin{equation}
\label{eq:raj}
\begin{split}
T_s(z) & = \int_0^{\pi/h} u_s(\theta)^\ast u_s(\theta) \frac{d \theta}{\frac4{h^2} \sin^2(\frac{\theta h}2) - z+\mu_s} \\
 & = \int_0^{4/{h^2}}  \frac{u_s ( \frac2h \arcsin(\frac{h}2 \sqrt{\zeta}) )^\ast u_s ( \frac2h \arcsin(\frac{h}2 \sqrt{\zeta}) )}{h\sqrt\zeta \sqrt{4/{h^2}-\zeta}} \frac{d \zeta}{\zeta-z+\mu_s} \\
 & = \int_0^{4/{h^2}} \frac{{\mathscr U}_s(\zeta) d \zeta}{\zeta-z+\mu_s} = \int_{-2/h^2}^{2/h^2} \frac{{\mathscr U}_s(\zeta+\frac2{h^2}) d \zeta}{\zeta+\frac2{h^2}-z+\mu_s}, 
\end{split}
\end{equation}
where
\begin{equation}
\label{eq:holdf1}
{\mathscr U}_s(\zeta) := u_s(\zeta)^\ast u_s(\zeta),
\end{equation}
and
\begin{equation}
\label{eq:holdf2}
u_s(\zeta) := \frac{u_s( \frac2h \arcsin(\frac{h}2 \sqrt{\zeta}) )}{h^{1/2}\zeta^{1/4}(\frac4{h^2}-\zeta)^{1/4}}.
\end{equation}


Next, we
establishe Proposition \ref{prop:exTS1} below. Firstly, the following lemma holds.

%

\begin{Lemma}
\label{lem:hold2}
Let $V$ satisfy Assumption \ref{eq:hyppert} with $\nu_0 > 3$. Then, 
the map $\zeta \in (0,\frac4{h^2}) \mapsto {\mathscr U}_s(\zeta)$, $0 \le s \le d$, is locally $\alpha$-H\"older in the $\mathfrak S_1$-norm with $\alpha=1$.
\end{Lemma}

\begin{proof}
For $\zeta \in (0,\frac4{h^2})$, thanks to \eqref{eq:holdf1} and \eqref{eq:holdf2}, one has
\begin{equation}
\label{eq:fctghraj}
{\mathscr U}_s(\zeta) = \frac{u_s( \frac2h g(\zeta) )^\ast u_s( \frac2h g(\zeta) )}{h\zeta^{1/2}(\frac4{h^2}-\zeta)^{1/2}}, \quad g(\zeta) := \arcsin \Big( \frac{h}{2}\sqrt{\zeta} \Big), \quad \zeta \in (0,\tfrac{4}{h^2}).
\end{equation}
Fix $\zeta_0 \in (0,\frac4{h^2})$ and consider $(\zeta_0 - \delta,\zeta_0+\delta) \subset (0,\frac4{h^2})$ a neighborhood of $\zeta_0$, $\delta > 0$ small enough. For every $\zeta_1$, $\zeta_2 \in (\zeta_0 - \delta,\zeta_0+\delta)$, one has
\begin{equation}
\label{eq:hold3}
\begin{split}
& \Big\Vert {\mathscr U}_s(\zeta_1) - {\mathscr U}_s(\zeta_2) \Big\Vert_{\mathfrak S_1} \le \Big\Vert \frac{u_s( \frac2h g(\zeta_1) )^\ast u_s( \frac2h g(\zeta_1) ) - u_s( \frac2h g(\zeta_2) )^\ast u_s( \frac2h g(\zeta_2) )}{h\zeta_1^{1/2}(\frac4{h^2}-\zeta_1)^{1/2}} \Big\Vert_{\mathfrak S_1} \\
& +\Big\vert \frac{1}{h\zeta_1^{1/2}(\frac4{h^2}-\zeta_1)^{1/2}} - \frac{1}{h\zeta_2^{1/2}(\frac4{h^2}-\zeta_2)^{1/2}} \Big\vert \Big\Vert u_s\Big( \frac2h g(\zeta_2) \Big)^\ast u_s\Big( \frac2h g(\zeta_2) \Big) \Big\Vert_{\mathfrak S_1}.
\end{split}
\end{equation}

a) Let us treat the first term of the r.h.s. of \eqref{eq:hold3}. We have
\begin{equation}
\label{eq:hold4}
\begin{split}
\Big\Vert & \frac{u_s( \frac2h g(\zeta_1) )^\ast u_s( \frac2h g(\zeta_1) ) - u_s( \frac2h g(\zeta_2) )^\ast u_s( \frac2h g(\zeta_2) )}{h\zeta_1^{1/2}(\frac4{h^2}-\zeta_1)^{1/2}} \Big\Vert_{\mathfrak S_1} \\
& \hspace{1cm} \le \frac{1}{h\zeta_1^{1/2}(\frac4{h^2}-\zeta_1)^{1/2}} \int_{\min(\zeta_1,\zeta_2)}^{\max(\zeta_1,\zeta_2)} \Big\Vert \partial_\zeta  \Big[ u_s\Big( \frac2h g(\zeta) \Big)^\ast u_s\Big( \frac2h g(\zeta) \Big) \Big] \Big\Vert_{\mathfrak S_1} d\zeta.
\end{split}
\end{equation}
From \eqref{eq:ut1}, one gets for $\zeta \in (0,\frac4{h^2})$
\begin{equation}
\label{eq:hold41}
\begin{split}
 u_s\Big( \frac2h g(\zeta) \Big)^\ast u_s\Big( \frac2h g(\zeta) \Big) 
= \mathscr V^{1/2} \Big[ \mathscr A \Big( \frac2h g(\zeta) \Big)^\ast \mathscr A \Big( \frac2h g(\zeta) \Big) \otimes K \pi_s K^\ast \Big] \mathscr V^{1/2},
\end{split}
\end{equation}
so that
\begin{equation}
\label{eq:hold42}
\Big\Vert \partial_\zeta  \Big[ u_s \Big( \frac2h g(\zeta) \Big)^\ast u_s \Big( \frac2h g(\zeta) \Big) \Big] \Big\Vert_{\mathfrak S_1} \le \Vert \mathscr V \Vert \Vert K \Vert_{\mathfrak S_2}^2 \Big\Vert \partial_\zeta  \Big[ \mathscr A \Big( \frac2h g(\zeta) \Big)^\ast \mathscr A \Big( \frac2h g(\zeta) \Big) \Big] \Big\Vert_{\mathfrak S_1}.
\end{equation}
By using \eqref{eq:ut2}, one obtains for $\phi \in \ell^2(\bz_h)$
\begin{equation}
\label{eq:impr2av}
\begin{split}
& \frac{\pi}{h} \mathscr A \Big( \frac2h g(\zeta) \Big)^\ast \mathscr A \Big( \frac2h g(\zeta) \Big) \phi(hn) \\
& = \frac12 \Big( e^{-2in g(\zeta)} \bra (hn)h^{-1} \ket^{-\nu_0/2} \sum_{m\in\bz} e^{2im g(\zeta)} \bra (hm)h^{-1} \ket^{-\nu_0/2} \phi(hm) \\
& + e^{2in g(\zeta)} \bra (hn)h^{-1} \ket^{-\nu_0/2} \sum_{m\in\bz} e^{-2im g(\zeta)} \bra (hm)h^{-1} \ket^{-\nu_0/2} \phi(hm) \Big) \\
& = \sum_{m\in\bz}\cos[2h(n+m)h^{-1} g(\zeta)] \bra (hn)h^{-1} \ket^{-\nu_0/2} \bra (hm)h^{-1} \ket^{-\nu_0/2} \phi(hm).
\end{split}
\end{equation}
Consequently, for any $\phi \in \ell^2(\bz_h)$, by setting $C(\zeta) := -\frac{1}{\sqrt \zeta \sqrt{4/h^2-\zeta}}$,  one gets
\begin{equation*}
\begin{split}
& \frac{\pi}{h} \partial_\zeta  \Big[ \mathscr A \Big( \frac2h g(\zeta) \Big)^\ast \mathscr A \Big( \frac2h g(\zeta) \Big) \Big] \phi(hn) = C(\zeta) \bra (hn)h^{-1} \ket^{-\nu_0/2} \\
& \hspace{2cm} \times \sum_{m\in\bz}h(n+m)h^{-1} \sin[2h(n+m)h^{-1} g(\zeta)] \bra (hm)h^{-1} \ket^{-\nu_0/2} \phi(hm) \\
& = C(\zeta) \Big( \bra (hn)h^{-1} \ket^{-\nu_0/2} (hn)h^{-1} \sin[2(hn)h^{-1} g(\zeta)] \\
& \hspace{2cm} \times \sum_{m\in\bz}\cos[2(hm)h^{-1} g(\zeta)] \bra (hm)h^{-1} \ket^{-\nu_0/2} \phi(hm) \\
& + \bra (hn)h^{-1} \ket^{-\nu_0/2} (hn)h^{-1} \cos[2(hn)h^{-1} g(\zeta)] \\
& \hspace{2cm} \times \sum_{m\in\bz}\sin[2(hm)h^{-1} g(\zeta)] \bra (hm)h^{-1} \ket^{-\nu_0/2} \phi(hm) \\
& + \bra (hn)h^{-1} \ket^{-\nu_0/2} \sin[2(hn)h^{-1} g(\zeta)] \\
& \hspace{2cm} \times \sum_{m\in\bz}\cos[2(hm)h^{-1} g(\zeta)] \bra (hm)h^{-1} \ket^{-\nu_0/2} (hm)h^{-1} \phi(hm) \\
& + \bra (hn)h^{-1} \ket^{-\nu_0/2} \cos[2(hn)h^{-1} g(\zeta)] \\
& \hspace{2cm} \times \sum_{m\in\bz}\sin[2(hm)h^{-1} g(\zeta)] \bra (hm)h^{-1} \ket^{-\nu_0/2} (hm)h^{-1} \phi(hm) \Big).
\end{split}
\end{equation*}
Il follows, using the inequality $\Vert A_1A_2 \Vert_{\mathfrak S_1} \le \Vert A_1 \Vert_{\mathfrak S_2} \Vert A_2 \Vert_{\mathfrak S_2}$ for $A_1$, $A_2 \in \mathfrak S_2$, that
\begin{equation*}
\begin{split}
& \frac{\pi}{h} \Big\Vert \partial_\zeta  \Big[ \mathscr A \Big( \frac2h g(\zeta) \Big)^\ast \mathscr A \Big( \frac2h g(\zeta) \Big) \Big] \Big\Vert_{\mathfrak S_1} \\
& \hspace{0.5cm} \le \vert C(\zeta) \vert \Big[ \Big( \sum_{n\in\bz} \sin^2[2n g(\zeta)] n^2 \bra n \ket^{-\nu_0} \Big)^{1/2} \Big( \sum_{n\in\bz} \cos^2[2n g(\zeta)] \bra n \ket^{-\nu_0} \Big)^{1/2} \\
& \hspace{1cm} +\Big( \sum_{n\in\bz} \cos^2[2n g(\zeta)] n^2 \bra n \ket^{-\nu_0} \Big)^{1/2} \Big( \sum_{n\in\bz} \sin^2[2n g(\zeta)] \bra n \ket^{-\nu_0} \Big)^{1/2} \\
& \hspace{1cm} + \Big( \sum_{n\in\bz} \sin^2[2n g(\zeta)] \bra n \ket^{-\nu_0} \Big)^{1/2} \Big( \sum_{n\in\bz} \cos^2[2n g(\zeta)] n^2 \bra n \ket^{-\nu_0} \Big)^{1/2} \\
& \hspace{1cm} +\Big( \sum_{n\in\bz} \cos^2[2n g(\zeta)] \bra n \ket^{-\nu_0} \Big)^{1/2} \Big( \sum_{n\in\bz} \sin^2[2n g(\zeta)] n^2 \bra n \ket^{-\nu_0} \Big)^{1/2} \Big],
\end{split}
\end{equation*}
so that
\begin{equation}
\label{eq:hold5}
\frac{\pi}{h} \Big\Vert \partial_\zeta  \Big[ \mathscr A \Big( \frac2h g(\zeta) \Big)^\ast \mathscr A \Big( \frac2h g(\zeta) \Big) \Big] \Big\Vert_{\mathfrak S_1} \le 4 \vert C(\zeta) \vert \sum_{n\in\bz} n^2 \bra n \ket^{-\nu_0}.
\end{equation}
Putting together \eqref{eq:hold4}, \eqref{eq:hold42} and \eqref{eq:hold5}, one gets finally 
\begin{equation}
\label{eq:hold4a}
\begin{split}
\Big\Vert & \frac{u_s( \frac2h g(\zeta_1) )^\ast u_s( \frac2h g(\zeta_1) ) - u_s( \frac2h g(\zeta_2) )^\ast u_s( \frac2h g(\zeta_2) )}{h\zeta_1^{1/2}(\frac4{h^2}-\zeta_1)^{1/2}} \Big\Vert_{\mathfrak S_1} \\
& \hspace{1cm} \le \frac4\pi (\max_{\zeta \in (\zeta_0 - \delta,\zeta_0+\delta)})^2 \frac{\Vert \mathscr V \Vert \Vert K \Vert_{\mathfrak S_2}^2 \sum_{n\in\bz} n^2 \bra n \ket^{-\nu_0}}{\zeta^{1/2} (\frac4{h^2}-\zeta)^{1/2}} \vert \zeta_1 - \zeta_2 \vert.
\end{split}
\end{equation}

b) Now, one treats the second term of the r.h.s. of \eqref{eq:hold3}. We have
\begin{equation}
\label{eq:hold6-}
\begin{split}
& \Big\vert \frac{1}{h\zeta_1^{1/2}(\frac4{h^2}-\zeta_1)^{1/2}} - \frac{1}{h\zeta_2^{1/2}(\frac4{h^2}-\zeta_2)^{1/2}} \Big\vert \Big\Vert u_s \Big( \frac2h g(\zeta_2) \Big)^\ast u_s \Big( \frac2h g(\zeta_2) \Big) \Big\Vert_{\mathfrak S_1} \\
& \le \frac1h \max_{\zeta \in (\zeta_0 - \delta,\zeta_0+\delta)} \Big\vert \partial_\zeta  \frac{1}{\zeta^{1/2}(\frac4{h^2}-\zeta)^{1/2}} \Big\vert \Big\Vert u_s \Big( \frac2h g(\zeta_2) \Big)^\ast u_s \Big( \frac2h g(\zeta_2) \Big) \Big\Vert_{\mathfrak S_1} \vert \zeta_1 - \zeta_2 \vert.
\end{split}
\end{equation}
From \eqref{eq:hold41}, it follows that for $\zeta \in (0,\frac4{h^2})$,
\begin{equation}
\label{eq:hold6}
\Big\Vert u_s \Big( \frac2h g(\zeta) \Big)^\ast u_s \Big( \frac2h g(\zeta) \Big) \Big\Vert_{\mathfrak S_1} \le \Vert \mathscr V \Vert \Vert K \Vert_{\mathfrak S_2}^2 \Big\Vert \mathscr A \Big( \frac2h g(\zeta) \Big)^\ast \mathscr A \Big( \frac2h g(\zeta) \Big)\Big\Vert_{\mathfrak S_1}.
\end{equation}
Using \eqref{eq:impr2av}, one obtains that for $\phi \in \ell^2(\bz_h)$
\begin{equation*}
\begin{split}
& \frac{\pi}{h} \mathscr A \Big( \frac2h g(\zeta) \Big)^\ast \mathscr A \Big( \frac2h g(\zeta) \Big) \phi(hn) \\
& = \cos[2(hn)h^{-1} g(\zeta)] \bra (hn)h^{-1} \ket^{-\nu_0/2} \sum_{m\in\bz} \cos[2(hm)h^{-1} g(\zeta)] \bra (hm)h^{-1} \ket^{-\nu_0/2} \phi(hm) \\
& + \sin[2(hn)h^{-1} g(\zeta)] \bra (hn)h^{-1} \ket^{-\nu_0/2} \sum_{m\in\bz} \sin[2(hm)h^{-1} g(\zeta)] \bra (hm)h^{-1} \ket^{-\nu_0/2} \phi(hm).
\end{split}
\end{equation*}
This gives that
\begin{equation}
\label{eq:hold7}
\begin{split}
\frac{\pi}{h} & \Big\Vert \mathscr A \Big( \frac2h g(\lambda) \Big)^\ast \mathscr A \Big( \frac2h g(\lambda) \Big) \Big\Vert_{\mathfrak S_1} \\
& \le \sum_{n\in\bz} \cos^2[2n g(\lambda)] \bra n \ket^{-\nu_0} + \sum_{n\in\bz} \sin^2[2n g(\lambda)] \bra n \ket^{-\nu_0} = \sum_{n\in\bz} \bra n \ket^{-\nu_0}.
\end{split}
\end{equation}
Using \eqref{eq:hold6-}, \eqref{eq:hold6} and \eqref{eq:hold7}, we finally get
\begin{equation}
\label{eq:hold7b}
\begin{split}
& \Big\vert \frac{1}{h\zeta_1^{1/2}(\frac4{h^2}-\zeta_1)^{1/2}} - \frac{1}{h\zeta_2^{1/2}(\frac4{h^2}-\zeta_2)^{1/2}} \Big\vert \Big\Vert u_s \Big( \frac2h g(\zeta_2) \Big)^\ast u_s \Big( \frac2h g(\zeta_2) \Big) \Big\Vert_{\mathfrak S_1} \\
& \le \frac1\pi \Vert \mathscr V \Vert \Vert K \Vert_{\mathfrak S_2}^2 \sum_{n\in\bz} \bra n \ket^{-\nu_0} \max_{\zeta \in (\zeta_0 - \delta,\zeta_0+\delta)} \Big\vert \partial_\zeta  \frac{1}{\zeta^{1/2}(\frac4{h^2}-\zeta)^{1/2}} \Big\vert \vert \zeta_1 - \zeta_2 \vert.
\end{split}
\end{equation}
Now, the lemma follows by putting together \eqref{eq:hold3}, \eqref{eq:hold4a} and \eqref{eq:hold7b}.
\end{proof}

As a direct consequence, applying Sokhotski-Plemelj formula \cite{plemelj formula}, the following corollary holds.

\begin{Corollary}
\label{cor:exTS1}
Let $z = \lambda + i\varepsilon$ with $\lambda \in (\mu_s,\frac4{h^2}+\mu_s)$ for a given $0 \le s \le d$.
Then,
\begin{equation}
\Big\Vert \int_0^{4/{h^2}} \frac{{\mathscr U}_s(\zeta) d \zeta}{\zeta-\lambda-i \varepsilon + \mu_s} - {\rm p. v.} \int_0^{4/{h^2}} \frac{{\mathscr U}_s(\zeta) d \zeta}{\zeta-\lambda+\mu_s} - i\pi {\mathscr U}_s(\lambda - \mu_s) \Big\Vert_{\mathfrak S_1} = \mathcal O(\varepsilon), \quad \varepsilon \searrow 0.
\end{equation}
\end{Corollary}

One can now state the next result showing the existence of $T_s(\lambda+i0)$ in $\mathfrak S_1$ for $\lambda \in (\mu_s,\frac4{h^2}+\mu_s)$. Notations are those introduced above. It follows from Corollary \ref{cor:exTS1} and \eqref{eq:raj} the following

\begin{Proposition}
\label{prop:exTS1}
Let $z = \lambda + i\varepsilon$ with $\lambda \in (\mu_s,\frac4{h^2}+\mu_s)$, $0 \le s \le d$ fixed. 
Then, as $\varepsilon \searrow 0$, 
\begin{equation}
T_s(\lambda + i \varepsilon) = \int_0^{4/{h^2}} \frac{{\mathscr U}_s(\zeta) d \zeta}{\zeta-\lambda -i \varepsilon + \mu_s} \quad \longrightarrow \quad {\rm p. v.} \int_0^{4/{h^2}} \frac{{\mathscr U}_s(\zeta) d \zeta}{\zeta-\lambda+\mu_s} + i\pi {\mathscr U}_s(\lambda - \mu_s).
\end{equation}
\end{Proposition}

The next corollary is a direct consequence of Corollary \ref{cor:exT}, Propositions \ref{prop:TS1-}, \ref{prop:exTS1} and \eqref{eq:Tdaint-}.

\begin{Corollary}
\label{cor:exTs1}
Let $V$ satisfy Assumption \ref{eq:hyppert} and $\lambda \in \br \smallsetminus \mathscr E_Q$. Then, the limit $T(\lambda+i0)$ exists in $\mathfrak S_1$ with $T(\lambda+i0) = T(\lambda)$. In particular, $T_s(\lambda+i0) = T_s(\lambda)$ for $\lambda \in (\mu_s,\frac4{h^2}+\mu_s)$, $0 \le s \le d$.
\end{Corollary}

For further references in the proof of the main results, let us point out the following remarks.

\begin{Remark}
\label{rem:noy0}
For $\lambda \in (\mu_s,\frac4{h^2} + \mu_s)$, one knows by definition that $T_s(\lambda)$ is given by \eqref{eq:Tdaint}, where $R(\lambda - \mu_s)$ admits the convolution kernel 
\begin{equation}
\label{eq:noy04}
R(\lambda - \mu_s,n-m) = \dfrac{ie^{2i\vert n-m \vert {\rm Arcsin}(\frac{h}{2} \sqrt{\lambda - \mu_s} ) } } {\sqrt{\lambda - \mu_s} \sqrt{4/{h^2} + \mu_s -\lambda}}.
\end{equation}
\end{Remark}

\begin{Remark}
\label{rem:noy2}
In the case $\lambda < \mu_s$, by the uniqueness of the limit and arguing as in the proof of ii) of Proposition \ref{prop:Ts2}, on can show that the formula \eqref{eq:Tdaint} for $T_s(\lambda)$ remains valid with $R(\lambda - \mu_s)$ admitting the convolution kernel 
\begin{equation}
\label{rem:noy0+}
R(\lambda - \mu_s,n-m) = \frac{e^{2i\vert n-m \vert {\rm Arcsin}(\frac{ih}{2} \sqrt{\mu_s-\lambda} ) } } {\sqrt{\mu_s-\lambda} \sqrt{4/{h^2} + \mu_s -\lambda}}.
\end{equation}
\end{Remark}


\section{Proof of the main results}\label{sec8}


One assumes that the potential $V$ satisfies Assumption \ref{eq:hyppert} throughout this section. 
To prove the next result for $0 \le s \le d$ fixed, one will exploit Proposition \ref{prop:exTS1} and Corollary \ref{cor:exTs1}. Nevertheless, note that it is also possible to use the formula 
\begin{equation*}
{\rm Im} \, T_s(\lambda) = \mathscr V^{1/2} ( M_\psi {\rm Im} \, R(\lambda-\mu_s) M_\psi \otimes K \pi_s K^\ast ) \mathscr V^{1/2}, \quad \lambda \in (\mu_s,\tfrac4{h^2}+\mu_s),
\end{equation*}
where ${\rm Im} \, R(\lambda-\mu_s)$ admits the convolution kernel 
\begin{equation*}
{\rm Im} \, R(\lambda-\mu_s,n-m) = \frac{\cos(2(n-m) {\rm Arcsin}(\frac{h}{2}\sqrt{\lambda-\mu_s} ) ) } {\sqrt{\lambda-\mu_s} \sqrt{4/{h^2}+\mu_s-\lambda}}.
\end{equation*}

\begin{Proposition}
\label{prop:impos}
For $\lambda \in (\mu_s,\frac4{h^2}+\mu_s)$, $0 \le s \le d$, we have $0 \le {\rm Im} \,T_s(\lambda) \in \mathfrak S_1$. Moreover, 
\begin{equation}
\label{eq:intrepTim}
{\rm Im} \,T_s(\lambda) = \frac{1}{\sqrt{\lambda-\mu_s} \sqrt{4/{h^2}+\mu_s-\lambda}} b_s(\lambda)^\ast b_s(\lambda),
\end{equation}
where $b_s(\lambda) : \ell^2(\bz_h,\mathcal{G}) \to \bc^2 \otimes \mathcal{G}$ is the operator defined by
\begin{equation}
b_s(\lambda) := (Y_s(\lambda) \otimes \pi_sK^\ast) \mathscr V^{1/2},
\end{equation}
with $Y_s(\lambda) : \ell^2(\bz_h) \to \bc^2$ defined by \eqref{eq:dim}.
\end{Proposition}

\begin{proof}
Thanks to Proposition \ref{prop:exTS1} and Corollary \ref{cor:exTs1}, 
\begin{equation}
\label{eq:pim2}
{\rm Im} \,T_s(\lambda) = \pi {\mathscr U}_s(\lambda-\mu_s) \in \mathfrak S_1, \qquad \lambda \in (\mu_s,4/{h^2}+\mu_s), \quad 0 \le s \le d,
\end{equation}
where w.r.t. the notations of Proposition \ref{prop:intrepT}, we have
\begin{equation}
{\mathscr U}_s(\lambda-\mu_s) = u_s(\lambda-\mu_s)^\ast u_s(\lambda-\mu_s), \qquad u_s(\lambda-\mu_s) = \frac{u_s(\frac2h g_s(\lambda))}{h^{1/2}(\lambda-\mu_s)^{1/4}(\frac4{h^2}+\mu_s-\lambda)^{1/4}}.
\end{equation}
It follows that ${\rm Im} \,T_s(\lambda) \ge 0$. From \eqref{eq:ut1}, one gets
\begin{equation}
\label{eq:impr1}
\begin{split}
{\rm Im} \,T_s(\lambda) & = \frac{\pi}{h\sqrt{\lambda-\mu_s} \sqrt{4/{h^2}+\mu_s-\lambda}} u_s \Big( \frac2h g_s(\lambda) \Big)^\ast u \Big( \frac2h g_s(\lambda) \Big) \\
& = \frac{\pi}{h\sqrt{\lambda-\mu_s} \sqrt{4/{h^2}+\mu_s-\lambda}} \mathscr V^{1/2} \Big[ \mathscr A \Big( \frac2h g_s(\lambda) \Big)^\ast \mathscr A \Big( \frac2h g_s(\lambda) \Big) \otimes K \pi_s K^\ast \Big] \mathscr V^{1/2}.
\end{split}
\end{equation}
Thanks to \eqref{eq:impr2av}, one has for $\phi \in \ell^2(\bz_h)$
\begin{equation}
\label{eq:impr2}
\begin{split}
& \frac{\pi}{h} \mathscr A \Big( \frac2h g_s(\lambda) \Big)^\ast \mathscr A \Big( \frac2h g_s(\lambda) \Big) \phi(hn) \\
& = \cos[2(hn)h^{-1} g_s(\lambda)] \bra (hn)h^{-1} \ket^{-\nu_0/2} \sum_{m\in\bz} \cos[2(hm)h^{-1} g_s(\lambda)] \bra (hm)h^{-1} \ket^{-\nu_0/2} \phi(hm) \\
& + \sin[2(hn)h^{-1} g_s(\lambda)] \bra (hn)h^{-1} \ket^{-\nu_0/2} \sum_{m\in\bz} \sin[2(hm)h^{-1} g_s(\lambda)] \bra (hm)h^{-1} \ket^{-\nu_0/2} \phi(hm) \\
& = Y_s(\lambda)^\ast Y_s(\lambda) \phi(hn),
\end{split}
\end{equation}
where the operator $Y_s(\lambda)^\ast$ is given by \eqref{eq:dimm}.
Putting together \eqref{eq:impr1} and \eqref{eq:impr2}, one gets \eqref{eq:intrepTim}.
\end{proof}

It is useful to recall the following standard properties of the counting functions $\mathscr N_\pm$ defined by \eqref{eq:fcompt}. If $T_1 = T_1^\ast$ and $T_2 = T_2^\ast$ belong to $\mathfrak S_\infty(\mathcal G)$, then one has the Weyl inequalities
\begin{equation}
\label{eq:inegweyl}
\mathscr N_\pm(x_1 + x_2,T_1 + T_2) \le \mathscr N_\pm(x_1,T_1) + \mathscr N_\pm(x_2,T_2), \quad x_1, x_2 > 0.
\end{equation}
If $T \in \mathfrak S_p(\mathcal G)$ for some $p \ge 1$, then 
\begin{equation}
\label{eq:inegsp}
\mathscr N_\pm(x,T) \le x^{-p} \Vert T \Vert_{\mathfrak S_p}^{-p}, \quad x > 0.
\end{equation} 
 

\subsection{Proof of Theorem \ref{theo1}}

We assume the conditions of Theorem \ref{theo1}. The following preliminary result holds.

\begin{Proposition}
\label{prop:prel}
As $\lambda \nearrow \mu_s$, $0 \le s \le d$, one has the estimates
\begin{equation*}
\mathscr N_\pm (1+\varepsilon,T_s(\lambda)) + \mathcal O(1) \le \mp \xi(\lambda;H^\mp,H_Q) \le \mathscr N_\pm (1-\varepsilon,T_s(\lambda)) + \mathcal O(1),
\end{equation*}
for any $\varepsilon \in (0,1)$.
\end{Proposition}
 
 \begin{proof}
 Fix $\varepsilon \in (0,1)$ and $0 \le s \le d$. First, note that 
\begin{equation}
\label{eq:pim1}
{\rm Re} \, T_s(\lambda) = T_s(\lambda) \quad {\rm and} \quad {\rm Im} \, T_s(\lambda) = 0, \quad \lambda \in \br \smallsetminus \sigma(\mu_s,\tfrac{4}{h^2}+\mu_s).
\end{equation}
In particular, for $\lambda \nearrow \mu_s$ we have 
${\rm Re} \, T(\lambda) = T_s(\lambda) + {\rm Re} \,( T(\lambda) - T_s(\lambda))$.
Then as $\lambda \nearrow \mu_s$, it follows from the Weyl inequalities \eqref{eq:inegweyl}, Corollary \ref{cor:exTs1} and Lemma \ref{lem:pushCV} that
 \begin{equation}
 \label{eq:nsr}
 \begin{split}
& \int_{\br} \mathscr N_\pm (1+ \varepsilon, T_s(\lambda))  \frac{dt}{\pi(1+t^2)} - \mathscr N_\mp ( \varepsilon/2, {\rm Re} \,( T(\lambda) - T_s(\lambda)) ) - \frac{2}{\pi \varepsilon} \Vert {\rm Im} \, T(\lambda) \Vert_{\mathfrak S_1}\\
& \le \int_{\br} \mathscr N_\pm(1,A(\lambda+i0)+tB(\lambda+i0))  \frac{dt}{\pi(1+t^2)} = \mp \xi(\lambda;H^\mp,H_Q) \\
& \le \int_{\br} \mathscr N_\pm (1-\varepsilon,T_s(\lambda))  \frac{dt}{\pi(1+t^2)} +\mathscr N_\pm ( \varepsilon/2, {\rm Re} \,( T(\lambda) - T_s(\lambda)) ) + \frac{2}{\pi \varepsilon} \Vert {\rm Im} \, T(\lambda) \Vert_{\mathfrak S_1}.
\end{split}
 \end{equation}
 For $\lambda$ sufficiently close to $\mu_s$, write
 \begin{equation}
\label{eqres:noy4+dd}
\begin{split}
T(\lambda) & = \mathscr V^{1/2} ( M_\psi R(\lambda - \mu_s) M_\psi \otimes K \pi_sK^\ast ) \mathscr V^{1/2} \\
& \qquad + \sum_{s' : \, \mu_s - \mu_{s'} \in (0,\tfrac{4}{h^2})} \mathscr V^{1/2} ( M_\psi R(\lambda - \mu_{s'}) M_\psi \otimes K \pi_{s'} K^\ast ) \mathscr V^{1/2}  \\
& \qquad \qquad + \sum_{s' : \, \mu_s - \mu_{s'} \notin [0,\tfrac{4}{h^2}]} \mathscr V^{1/2} ( M_\psi R(\lambda - \mu_{s'}) M_\psi \otimes K \pi_{s'} K^\ast ) \mathscr V^{1/2},
\end{split}
\end{equation}
It follows that as $\lambda \searrow \mu_s$, we have ${\rm Re} \,( T(\lambda) - T_s(\lambda)) = {\rm Re} \, Z_s(\lambda)$ with
\begin{equation}
\label{eq:diff0-4d}
\begin{split}
Z_s(\lambda) := \sum_{s' : \, \mu_s - \mu_{s'} \in (0,\tfrac{4}{h^2})} & \mathscr V^{1/2} ( M_\psi R(\lambda - \mu_{s'}) M_\psi \otimes K \pi_{s'} K^\ast ) \mathscr V^{1/2}  \\
& + \sum_{s' : \, \mu_s - \mu_{s'} \notin [0,\tfrac{4}{h^2}]} \mathscr V^{1/2} ( M_\psi R(\lambda - \mu_{s'}) M_\psi \otimes K \pi_{s'} K^\ast ) \mathscr V^{1/2},
\end{split}
\end{equation}
and
\begin{equation}
\label{eq:diff0-4dd}
{\rm Im} \, T(\lambda) = {\rm Im} \, Z_s(\lambda) = \sum_{s' : \, \mu_s - \mu_{s'} \in (0,\tfrac{4}{h^2})} {\rm Im} \, T_{s'}(\lambda).
\end{equation}
Note that the operator $Z_s(\mu_s)$ is well defined and is Hilbert-Schmidt. Similarly to \eqref{eq:Tdiff+}, we can prove that $\lim_{\lambda \nearrow \mu_s} \Vert Z_{s}(\lambda) - Z_{s}(\mu_s) \Vert_{\mathfrak S_2} = 0$, so that 
\begin{equation}
\label{eq:limitdiffg}
\lim_{\lambda \nearrow \mu_s} \Vert {\rm Re} \, ( T(\lambda) - T_s(\lambda)) - {\rm Re} \, Z_{s}(\mu_s) \Vert_{\mathfrak S_2} = 0.
\end{equation}
Therefore, by inequalities \eqref{eq:inegweyl} we have 
\begin{equation}
\label{eq:limitdiffg1}
\begin{split}
& \mathscr N_\pm ( \varepsilon/2, {\rm Re} \,( T(\lambda) - T_s(\lambda)) ) \\
&  \le \mathscr N_\pm ( \varepsilon/4, {\rm Re} \, ( T(\lambda) - T_s(\lambda)) - {\rm Re} \, Z_s(\mu_s) ) + \mathscr N_\pm ( \varepsilon/4, {\rm Re} \, Z_s(\mu_s) ) \\
& = \mathscr N_\pm ( \varepsilon/4, {\rm Re} \, Z_s(\mu_s) ) = \mathcal O(1) \quad as \quad \lambda \nearrow \mu_s.
\end{split}
\end{equation}
Otherwise, for $\mu_s \in (\mu_{s'},\frac{4}{h^2}+\mu_{s'})$, we have $\lim_{\lambda \nearrow \mu_s} \Vert {\rm Im} \, T_{s'}(\lambda) - {\rm Im} \, T_{s'}(\mu_s) \Vert_{\mathfrak S_1} = 0$ thanks to identity \eqref{eq:pim2} and Lemma \ref{lem:hold2}. This together with \eqref{eq:diff0-4dd} implies that
\begin{equation}
\label{eq:limitdiffg2}
\lim_{\lambda \nearrow \mu_s} \Vert {\rm Im} \, T(\lambda) - {\rm Im} \, Z_{s}(\mu_s) \Vert_{\mathfrak S_1} = 0.
\end{equation}
Finally, bearing in mind the identity 
$ \int_\br \frac{dt}{\pi(1+t^2)} = 1,$
the proposition follows from \eqref{eq:nsr}, \eqref{eq:limitdiffg1} and \eqref{eq:limitdiffg2}.
\end{proof}

Now, one shows in the next result that $\mathscr N_\pm (x,T_s(\lambda))$ can be bounded as $\lambda \to \mu_s$ fixed, from below and from above by expressions involving $\mathcal L_s(\lambda)$, up to $\mathcal O(1)$. Here, $\mathcal L_s(\lambda) : \ell^2(\bz_h) \otimes \mathcal{G} \to \ell^2(\bz_h) \otimes \mathcal{G}$ is the trace class operator defined by 
\begin{equation}
\label{eq:opdom0}
\mathcal L_s(\lambda) = \frac h2 \frac{L_{s}^\ast L_{s}}{\sqrt{\mu_s -\lambda}}, \qquad \lambda \nearrow \mu_s,
\end{equation}
where $L_{s}$ is defined by \eqref{eq:L0}.

\begin{Proposition}
\label{prop:dom0}
Let $\nu_0 > 3$ in Assumption \ref{eq:hyppert}. Then, as $\lambda \nearrow \mu_s$, $0 \le s \le d$, we have
\begin{equation*}
 \mathscr N_+ ((1+\varepsilon)x,\mathcal L_s(\lambda)) + \mathcal O(1) \le \mathscr N_+ (x,T_s(\lambda)) \le \mathscr N_+ ((1-\varepsilon)x,\mathcal L_s(\lambda)) + \mathcal O(1),
\end{equation*}
and
\begin{equation*}
 \mathcal O(1) \le \mathscr N_- (x,T_s(\lambda)) \le  \mathcal O(1),
\end{equation*}
for any $\varepsilon \in (0,1)$ and $x > 0$.
\end{Proposition}

\begin{proof}
Fix $\mu_s \in \mathscr E_Q$. As above, the main idea of the proof is to approximate the operator $T_s(\lambda) - \mathcal L_s(\lambda)$ in norm, as $\lambda \nearrow \mu_s$, by a compact operator independent of $\lambda$.

a) Let $\lambda \nearrow \mu_s$ in a small neighborhood containing $\mu_s$ as unique threshold. The convolution kernel $R(\lambda-\mu_s,n-m)$ given by \eqref{rem:noy0+} can be decomposed as
\begin{equation*}
\begin{split}
R(\lambda - \mu_s,n-m) & = \frac{h}{2\sqrt{\mu_s - \lambda}} + \Big( \frac{1}{\sqrt{\mu_s - \lambda} \sqrt{4/h^2 + \mu_s -\lambda}} - \frac{h}{2\sqrt{\mu_s -\lambda}} \Big) \\
& + \frac{e^{2i\vert n-m \vert {\rm Arcsin}(\frac{ih}{2} \sqrt{\mu_s - \lambda} ) } -1} {\sqrt{\mu_s -\lambda} \sqrt{4/h^2 + \mu_s -\lambda}}.
\end{split}
\end{equation*}
Together with Remark \ref{rem:noy2} and \eqref{eq:Tdaint-}, this implies that
\begin{equation}
\label{eq:diff}
T_s(\lambda) - \mathcal L_s(\lambda) = \mathscr V^{1/2} (\Xi_{\nu_0}^{(\lambda)}  \otimes K \pi_s K^\ast) \mathscr V^{1/2} + \mathcal I_s(\lambda),
\end{equation}
where $\Xi_{\nu_0}^{(\lambda)} : \ell^2(\bz_h) \to \ell^2(\bz_h)$ is the summation kernel operator defined by
\begin{equation*}
(\Xi_{\nu_0}^{(\lambda)} \varphi)(hn) := \sum_{m\in\bz} \bra n \ket^{-\nu_0/2} \frac{e^{2i\vert n-m \vert {\rm Arcsin}(\frac{ih}{2} \sqrt{\mu_s-\lambda} ) } -1} {\sqrt{\mu_s-\lambda} \sqrt{4/h^2+\mu_s-\lambda}} \bra m \ket^{-\nu_0/2} \varphi(hm),
\end{equation*}
$\varphi \in \ell^2(\bz_h)$,
\begin{equation*}
\mathcal I_s (\lambda) := \Big( \frac{1}{\sqrt{\mu_s-\lambda} \sqrt{4/h^2+\mu_s-\lambda}} - \frac{h}{2\sqrt{\mu_s-\lambda}} \Big) L_{s}^\ast L_{s},
\end{equation*}
with the operator $L_{s}$ given by \eqref{eq:L0}.
Since $\frac{1}{\sqrt{\mu_s-\lambda} \sqrt{4/h^2+\mu_s-\lambda}} - \frac{h}{2\sqrt{\mu_s-\lambda}} = \mathcal O(\sqrt{\mu_s-\lambda})$ as $\lambda \nearrow \mu_s$, and the operator $L_{s}$ is independent of $\lambda$, it follows that
\begin{equation}
\label{eq:I0}
\lim_{\lambda \nearrow \mu_s} \Vert \mathcal I_s (\lambda) \Vert_{\mathfrak S_2} = \lim_{\lambda \nearrow \mu_s} \mathcal O(\sqrt{\mu_s-\lambda}) \Vert L_{s}^\ast L_{s} \Vert_{\mathfrak S_2} = 0.
\end{equation}
Define the operator
\begin{equation}
\label{eq:opS0op}
\mathcal T_s = \mathscr V^{1/2} (\Xi_{\nu_0}^{(0)} \otimes K\pi_sK^\ast) \mathscr V^{1/2},
\end{equation}
where $\Xi_{\nu_0}^{(0)} : \ell^2(\bz_h) \to \ell^2(\bz_h)$ is the summation kernel operator given by
\begin{equation}
\label{eq:opS0}
(\Xi_{\nu_0}^{(0)} \phi)(hn) := - \frac{h^2}2 \sum_{m\in\bz} \bra n \ket^{-\nu_0/2} \vert n-m \vert \bra m \ket^{-\nu_0/2} \phi(hm), \qquad \phi \in \ell^2(\bz_h).
\end{equation}
Since $\nu_0 > 3$, then $\sum_{n,m} \vert \bra n \ket^{-\nu_0/2} \vert n-m \vert \bra m \ket^{-\nu_0/2} \vert^2 < \infty$ and $\Xi_{\nu_0}^{(0)}$ belongs to $\mathfrak S_2(\ell^2(\bz_h))$. In particular, the operator $\mathcal T_s$ is compact in $\ell^2(\bz_h,\mathcal{G})$. By using the Lebesgue dominated convergence theorem and the convolution kernels of the operators $\Xi_{\nu_0}^{(\lambda)}$ and $\Xi_{\nu_0}^{(0)}$, one gets
\begin{equation}
\label{eq:limdiff}
\lim_{\lambda \nearrow \mu_s} \Vert \Xi_{\nu_0}^{(\lambda)} - \Xi_{\nu_0}^{(0)} \Vert_{\mathfrak S_2}^2 = 0.
\end{equation}
Putting together \eqref{eq:diff}, \eqref{eq:I0} and \eqref{eq:limdiff}, one obtains
\begin{equation}
\label{eq:approxdiff}
\lim_{\lambda \nearrow \mu_s} \Vert T(\lambda) - \mathcal L_s(\lambda) - \mathcal T_s \Vert_{\mathfrak S_2} = 0.
\end{equation}

b) Now, consider $\lambda \nearrow \mu_s$ as above, $\varepsilon \in (0,1)$ and $x > 0$. Using Weyl's inequalities \eqref{eq:inegweyl}, one gets
\begin{equation*}
\begin{split}
 \mathscr N_\pm ((1+\varepsilon)x,\mathcal L_s(\lambda)) - \mathscr N_\mp (\varepsilon x, & T(\lambda) - \mathcal L_s(\lambda)) \le \mathscr N_\pm (x,T(\lambda)) \\
 & \le \mathscr N_\pm ((1-\varepsilon)x,\mathcal L_s(\lambda)) + \mathscr N_\pm (\varepsilon x,T(\lambda) - \mathcal L_s(\lambda)),
 \end{split}
 \end{equation*}
 Since $\mathcal L_s(\lambda)$ is a positive operator, then 
 \begin{equation*}
 \mathscr N_- (t,\mathcal L_s(\lambda)) = 0, \quad \forall \, t > 0.
 \end{equation*}
 Therefore, to get the proposition, it suffices to prove that for every $\varepsilon \in (0,1)$ and $x > 0$, 
 \begin{equation}
 \label{eq:born}
 \mathscr N_\pm (\varepsilon x, T(\lambda) - \mathcal L_s(\lambda)) = \mathcal O(1), \qquad \lambda \nearrow \mu_s.
 \end{equation}
 This follows by arguing as in \eqref{eq:limitdiffg1}.
\end{proof}

In the next result, for $0 \le s \le d$, $\mathcal L_{4,s}(\lambda) : \ell^2(\bz_h) \otimes \mathcal{G} \to \ell^2(\bz_h) \otimes \mathcal{G}$ is the trace class operator defined by 
\begin{equation}
\label{eq:opdom4}
\mathcal L_{4,s}(\lambda) = - \frac{h}2 \frac{L_{4,s}^\ast L_{4,s}}{\sqrt{\lambda-4/{h^2} - \mu_s}}, \qquad \lambda \searrow \tfrac{4}{h^2} + \mu_s, 
\end{equation}
where $L_{4,s}$ is defined by \eqref{eq:L4}.
The proof is similar to that of Proposition \ref{prop:dom0} and then will be shortened. Only the main quantities will be specified.

\begin{Proposition}
\label{prop:dom4}
Suppose $\nu_0 > 3$ in Assumption \ref{eq:hyppert}. Then, as $\lambda \searrow \frac4{h^2} + \mu_s$, $0 \le s \le d$, we have
\begin{equation*}
 \mathcal O(1) \le \mathscr N_+ (x,T(\lambda)) \le  \mathcal O(1),
\end{equation*}
and
\begin{equation*}
 \mathscr N_- ((1+\varepsilon)x,\mathcal L_{4,s}(\lambda)) + \mathcal O(1) \le \mathscr N_- (x,T(\lambda)) \le \mathscr N_- ((1-\varepsilon)x,\mathcal L_{4,s}(\lambda)) + \mathcal O(1),
\end{equation*}
for any $\varepsilon \in (0,1)$ and $x > 0$.
\end{Proposition}

\begin{proof}
For $\lambda \searrow \frac4{h^2} + \mu_s$, one can reduce the analysis near the threshold $\mu_s$ with $\lambda \nearrow \mu_s$, exploiting the symmetry between the two thresholds through the unitary operator $J$ given by  \eqref{op:rel0}.

Indeed, as $\lambda \searrow \frac4{h^2} + \mu_s$,  write
\begin{equation*}
\label{eqres:noy4+d}
\begin{split}
T(\lambda) & = \mathscr V^{1/2} ( J M_\psi J^{\ast}R(\lambda - \mu_s)J M_\psi J^\ast \otimes K \pi_sK^\ast ) \mathscr V^{1/2} \\
& \qquad + \sum_{s' : \, \mu_s - \mu_{s'} \in (0,\tfrac{4}{h^2})} \mathscr V^{1/2} ( J M_\psi J^\ast R(\lambda - \mu_{s'}) J M_\psi J^\ast \otimes K \pi_{s'} K^\ast ) \mathscr V^{1/2}  \\
& \qquad \qquad + \sum_{s' : \, \mu_s - \mu_{s'} \notin [0,\tfrac{4}{h^2}]} \mathscr V^{1/2} ( J M_\psi J^\ast R(\lambda - \mu_{s'}) J M_\psi J^\ast \otimes K \pi_{s'} K^\ast ) \mathscr V^{1/2},
\end{split}
\end{equation*}
where for $\bullet \in \{s,s'\}$, we have $J^{\ast}R(\lambda - \mu_\bullet)J = -R(\tfrac{4}{h^2} + \mu_\bullet - \lambda)$. It follows that 
\begin{equation}
\label{eqres:noy4+1d}
\begin{split}
T(\lambda) & = -\mathscr V^{1/2} ( J M_\psi R(\tfrac{4}{h^2} + \mu_s - \lambda) M_\psi J^\ast \otimes K \pi_sK^\ast ) \mathscr V^{1/2} \\
& \qquad - \sum_{s' : \, \mu_s - \mu_{s'} \in (0,\tfrac{4}{h^2})} \mathscr V^{1/2} ( J M_\psi R(\tfrac{4}{h^2} + \mu_{s'} - \lambda) M_\psi J^\ast \otimes K \pi_{s'} K^\ast ) \mathscr V^{1/2}  \\
& \qquad \qquad - \sum_{s' : \, \mu_s - \mu_{s'} \notin [0,\tfrac{4}{h^2}]} \mathscr V^{1/2} ( J M_\psi R(\tfrac{4}{h^2} + \mu_{s'} - \lambda) M_\psi J^\ast \otimes K \pi_{s'} K^\ast ) \mathscr V^{1/2}.
\end{split}
\end{equation}
Now, one can observe that we have $\frac{4}{h^2} + \mu_{s} - \lambda = (\frac{4}{h^2} + 2\mu_{s} - \lambda) - \mu_s$ so that $\frac{4}{h^2} + 2\mu_{s} - \lambda \nearrow \mu_s$ as $\lambda \searrow \frac4{h^2} + \mu_s$. Using \eqref{eqres:noy4+1d}, the claim follows by arguing as in Propositions \ref{prop:prel} and \ref{prop:dom4}.
\end{proof}

For $x > 0$ and $\lambda \nearrow \mu_s$, one has
\begin{equation}
\label{eq:top0}
\begin{split}
\mathscr N_+ & (x,\mathcal L_s(\lambda)) = \mathscr N_+ \Big( x, \frac{h L_{s}^\ast L_{s}}{2\sqrt{\mu_s-\lambda}} \Big) = \mathscr N_+ \Big( x,\frac{h L_{s} L_{s}^\ast}{2\sqrt{\mu_s-\lambda}} \Big) \\ 
& = \mathscr N_+ \Big( x,\frac{h (\bra \psi \vert \otimes \pi_s K^\ast) \mathscr V (\bra \psi \vert^\ast \otimes K \pi_s) }{2\sqrt{\mu_s-\lambda}} \Big) = \mathscr N_+ ( x,\omega_0(\lambda)),
\end{split}
\end{equation}
where the operator $\omega_0(\lambda)$ is given by \eqref{eq:opdom04}. Similarly, we show that for any $x > 0$ and $\lambda \searrow \frac4{h^2} + \mu_s$,
\begin{equation}
\label{eq:top4}
\mathscr N_- (x,\mathcal L_{4,s}(\lambda)) = \mathscr N_- \Big( x, -\frac{h L_{4,s}^\ast L_{4,s}}{2\sqrt{\lambda-4/h^2-\mu_s}} \Big) = \mathscr N_- ( x,\omega_{4,s}(\lambda)),
\end{equation}
where the operator $\omega_{4,s}(\lambda)$ is given by \eqref{eq:opdom04}. Now, Theorem \ref{theo1} follows directly from Propositions \ref{prop:dom0} and \ref{prop:dom4} together with identities \eqref{eq:top0} and \eqref{eq:top4}.

\subsection{Proof of Theorem \ref{theo2}}

This section concerns the case $\lambda \in \searrow \mu_s$ and $\lambda \nearrow \frac4{h^2}+\mu_s$, $0 \le s \le d$. We assume that $V$ satisfies Assumption \ref{eq:hyppert}. In the next result, one shows the boundedness of $\mathscr N_\pm (x,{\rm Re} \,T(\lambda))$ as $\lambda \to \lambda_s \in \{\mu_s,\frac4{h^2}+\mu_s \}$, $s$ fixed.

\begin{Proposition}
\label{prop:dom04}
Suppose $\nu_0 > 3$ in Assumption \ref{eq:hyppert}. Then, for any $x > 0$,
\begin{equation*}
 \mathscr N_\pm (x,{\rm Re} \, T(\lambda)) =  \mathcal O(1), 
\end{equation*}
as $\lambda \searrow \mu_s$ and $\lambda \nearrow \frac4{h^2}+\mu_s$ for $0 \le s \le d$ fixed.
\end{Proposition}

\begin{proof}
The idea is to approximate the operator ${\rm Re} \, T(\lambda) $ in norm, as $\lambda \searrow \mu_s$ and $\lambda \nearrow \frac4{h^2}+\mu_s$, by a compact operator independent of $\lambda$.

a) First, let us focus on the case $\lambda \searrow \mu_s$ fixed. Write
\begin{equation}
\label{eqres:noy4+df}
\begin{split}
T(\lambda) & = \mathscr V^{1/2} ( M_\psi R(\lambda - \mu_s) M_\psi \otimes K \pi_sK^\ast ) \mathscr V^{1/2} \\
& \qquad + \sum_{s' : \, \mu_s - \mu_{s'} \in (0,\tfrac{4}{h^2})} \mathscr V^{1/2} ( M_\psi R(\lambda - \mu_{s'}) M_\psi \otimes K \pi_{s'} K^\ast ) \mathscr V^{1/2}  \\
& \qquad \qquad + \sum_{s' : \, \mu_s - \mu_{s'} \notin [0,\tfrac{4}{h^2}]} \mathscr V^{1/2} ( M_\psi R(\lambda - \mu_{s'}) M_\psi \otimes K \pi_{s'} K^\ast ) \mathscr V^{1/2}.
\end{split}
\end{equation}
Thanks to \eqref{eq:noy04}, the operator ${\rm Re} \, R(\lambda-\mu_s)$ admits the convolution kernel
\begin{equation*}
\label{eq:Renoy04}
{\rm Re} \, R(\lambda-\mu_s,n-m) = -\frac{\sin ( 2\vert n-m \vert {\rm Arcsin}(\frac{h}{2} \sqrt{\lambda-\mu_s} )) } {\sqrt{\lambda-\mu_s} \sqrt{4/{h^2}+\mu_s-\lambda}}.
\end{equation*}
So, one obtains
\begin{equation}
\label{eq:diff04}
{\rm Re} \, T(\lambda) = \mathscr V^{1/2} (E_{\nu_0}^{(\lambda)} \otimes K \pi_s K^\ast) \mathscr V^{1/2} + \mathcal Z_s(\lambda),
\end{equation}
where $E_{\nu_0}^{(\lambda)} : \ell^2(\bz_h) \to \ell^2(\bz_h)$ is the convolution kernel operator given by
\begin{equation*}
(E_{\nu_0}^{(\lambda)} \phi)(hn) := \sum_{m\in\bz} \bra n \ket^{-\nu_0/2} {\rm Re} \, R(\lambda,n-m) \bra m \ket^{-\nu_0/2} \phi(hm), \quad \phi \in \ell^2(\bz_h),
\end{equation*}
and 
\begin{equation}
\label{eq:diff0-4}
\begin{split}
\mathcal Z_s(\lambda) & := \sum_{s' : \, \mu_s - \mu_{s'} \in (0,\tfrac{4}{h^2})} \mathscr V^{1/2} ( M_\psi {\rm Re} \, R(\lambda - \mu_{s'}) M_\psi \otimes K \pi_{s'} K^\ast ) \mathscr V^{1/2}  \\
& \qquad \qquad + \sum_{s' : \, \mu_s - \mu_{s'} \notin [0,\tfrac{4}{h^2}]} \mathscr V^{1/2} ( M_\psi {\rm Re} \, R(\lambda - \mu_{s'}) M_\psi \otimes K \pi_{s'} K^\ast ) \mathscr V^{1/2}.
\end{split}
\end{equation}
By using Lebesgue's dominated convergence theorem, one gets
\begin{equation}
\label{eq:limdiff04}
\lim_{\lambda \searrow \mu_s} \Vert E_{\nu_0}^{(\lambda)} - \Xi_{\nu_0}^{(0)} \Vert_{\mathfrak S_2}^2 = 0,
\end{equation}
where $\Xi_{\nu_0}^{(0)}$ is the operator given by \eqref{eq:opS0}. Similarly to \eqref{eq:Tdiff+}, we can prove that
\begin{equation}
\label{eq:imdiff0-4}
\lim_{\lambda \searrow \mu_s} \Vert \mathcal Z_{s}(\lambda) - \mathcal Z_{s}(\mu_s) \Vert_{\mathfrak S_2} = 0.
\end{equation}
It follows from \eqref{eq:diff04}, \eqref{eq:limdiff04} and \eqref{eq:imdiff0-4} that 
\begin{equation*}
\lim_{\lambda \searrow \mu_s} \Vert {\rm Re} \, T(\lambda) - \mathscr V^{1/2} (\Xi_{\nu_0}^{(0)} \otimes K \pi_s K^\ast) \mathscr V^{1/2} - \mathcal Z_s(\mu_s) \Vert_{\mathfrak S_2} = 0,
\end{equation*}
Now, the claim follows by arguing as in part b) of the proof of Proposition \ref{prop:dom0}. 

b) The case $\lambda \nearrow \frac4{h^2} + \mu_s$ can be proved as follows. Write
\begin{equation}
\label{eqres:noy04+d}
\begin{split}
T(\lambda) &= \mathscr V^{1/2} ( J M_\psi R_J(\lambda-\mu_s) M_\psi J^\ast \otimes K \pi_s K^\ast ) \mathscr V^{1/2} \\
& \qquad + \sum_{s' : \, \mu_s - \mu_{s'} \in (0,\tfrac{4}{h^2})} \mathscr V^{1/2} ( J M_\psi R_J(\lambda - \mu_{s'}) M_\psi J^\ast \otimes K \pi_{s'} K^\ast ) \mathscr V^{1/2}  \\
& \qquad \qquad + \sum_{s' : \, \mu_s - \mu_{s'} \notin [0,\tfrac{4}{h^2}]} \mathscr V^{1/2} ( J M_\psi R_J(\lambda - \mu_{s'}) M_\psi J^\ast \otimes K \pi_{s'} K^\ast ) \mathscr V^{1/2},
\end{split}
\end{equation}
where 
$R_J(\lambda-\mu_\bullet) := J^{\ast}R(\lambda-\mu_\bullet)J = -R(\tfrac{4}{h^2} + \mu_\bullet - \lambda)$, $\bullet = s, s'$, so that the operator $R_J(\lambda-\mu_s)$ admits the kernel 
\begin{equation}
\label{eq:noy04+d}
R_J(\lambda-\mu_s,n-m) = -\frac{ie^{2i\vert n-m \vert {\rm Arcsin}(\frac{h}{2}\sqrt{4/{h^2}+\mu_s-\lambda} ) } } {\sqrt{\lambda-\mu_s} \sqrt{4/{h^2}+\mu_s-\lambda}}.
\end{equation}
It follows that
\begin{equation}
\label{eqres:noy04++d}
\begin{split}
{\rm Re} \, T(\lambda) & = \mathscr V^{1/2} ( J M_\psi {\rm Re} \, R_J(\lambda-\mu_s) M_\psi J^\ast \otimes K \pi_s K^\ast ) \mathscr V^{1/2} \\
& \qquad + \sum_{s' : \, \mu_s - \mu_{s'} \in (0,\tfrac{4}{h^2})} \mathscr V^{1/2} ( J M_\psi {\rm Re} \, R_J(\lambda - \mu_{s'}) M_\psi J^\ast \otimes K \pi_{s'} K^\ast ) \mathscr V^{1/2}  \\
& \qquad \qquad + \sum_{s' : \, \mu_s - \mu_{s'} \notin [0,\tfrac{4}{h^2}]} \mathscr V^{1/2} ( J M_\psi {\rm Re} \, R_J(\lambda - \mu_{s'}) M_\psi J^\ast \otimes K \pi_{s'} K^\ast ) \mathscr V^{1/2},
\end{split}
\end{equation}
with ${\rm Re} \, R_J(\lambda - \mu_s)$ admitting the convolution kernel 
\begin{equation}
\label{eq:noy04++d}
{\rm Re} \, R_J(\lambda,n-m) = \frac{\sin(2\vert n-m \vert {\rm Arcsin}(\frac{h}{2}\sqrt{4/{h^2}+\mu_s-\lambda} ) ) } {\sqrt{\lambda-\mu_s} \sqrt{4/{h^2}+\mu_s-\lambda}}.
\end{equation}
Now, the rest of the proof follows as in a) above. 
\end{proof}

The next result uses in particular the identities (see e.g. \cite[Section 5.4]{ferrai})
\begin{equation}
\label{eq:eqarctan}
\int_{\br} \mathscr N_\pm (x,t\mathcal T)  \frac{dt}{\pi(1+t^2)} = \frac1{\pi} {\rm Tr} \arctan(x^{-1}\mathcal T), \qquad x > 0,
\end{equation}
where $0 \le \mathcal T = \mathcal T^\ast \in \mathfrak S_1$.

\begin{Proposition}
\label{prop:prel04}
Let $\nu_0 > 3$ in Assumption \ref{eq:hyppert}. As $\lambda \searrow \mu_s$ and $\lambda \nearrow \frac4{h^2}+\mu_s$, $0 \le s \le d$, the following bounds hold:
 \begin{align*}
 \frac1{\pi} {\rm Tr} & \arctan((x(1+\varepsilon))^{-1} {\rm Im} \, T_s(\lambda)) + \mathcal O(1) \\
& \le \int_{\br} \mathscr N_\pm(x,{\rm Re} \, T(\lambda)+t{\rm Im} \, T(\lambda))  \frac{dt}{\pi(1+t^2)} \\
&\hspace{0.6cm} \le \frac1{\pi} {\rm Tr} \arctan((x(1-\varepsilon))^{-1} {\rm Im} \, T_s(\lambda)) + \mathcal O(1),
 \end{align*}
for any $\varepsilon \in (0,1)$ and $x > 0$.
\end{Proposition}

\begin{proof}
It follows from the Weyl inequalities \eqref{eq:inegweyl} that for any $\varepsilon \in (0,1)$ and $x > 0$, 
\begin{equation}
\label{eq:ch1}
\begin{split}
 \mathscr N_\pm & ((1+\varepsilon)x,t{\rm Im} \, T_s(\lambda)) - \mathscr N_\mp ( \varepsilon x, t({\rm Im} \,T(\lambda) - {\rm Im} \,T_s(\lambda)) + {\rm Re} \, T(\lambda) ) \\
 & \le \mathscr N_\pm (x,{\rm Re} \, T(\lambda)+t{\rm Im} \, T(\lambda)) \\
& \hspace{1cm} \le \mathscr N_\pm ((1-\varepsilon)x,t{\rm Im} \, T_s(\lambda)) + \mathscr N_\pm ( \varepsilon x, t({\rm Im} \,T(\lambda) - {\rm Im} \,T_s(\lambda)) + {\rm Re} \, T(\lambda)).
 \end{split}
 \end{equation}
Now, focus on the case $\lambda \searrow \mu_s$. The case $\lambda \nearrow \frac4{h^2}+\mu_s$ follows in a similar way. We have
\begin{equation}
\label{eq:ch2}
\begin{split}
\mathscr N_\pm ( \varepsilon x, & t({\rm Im} \,T(\lambda) - {\rm Im} \,T_s(\lambda)) + {\rm Re} \, T(\lambda) ) \\
& \le \mathscr N_\pm ( \varepsilon x/2, t({\rm Im} \,T(\lambda) - {\rm Im} \,T_s(\lambda)) ) + \mathscr N_\pm ( \varepsilon x/2, {\rm Re} \, T(\lambda) ).
\end{split}
\end{equation}
 Let us treat the first term of the r.h.s. of \eqref{eq:ch2}. Using \eqref{eqres:noy4+df}, we get as $\lambda \searrow \mu_s$
 \begin{equation*}
\label{eq:ch3}
{\rm Im} \,T(\lambda) - {\rm Im} \,T_s(\lambda) = \sum_{s' : \, \mu_s - \mu_{s'} \in (0,\tfrac{4}{h^2})} {\rm Im} \,T_{s'}(\lambda).
\end{equation*}
Therefore, using \eqref{eq:eqarctan} one gets
\begin{equation}
\label{eq:ch4}
\begin{split}
\int_{\br} \mathscr N_\pm & \big( \varepsilon x/2, t({\rm Im} \,T(\lambda) - {\rm Im} \,T_s(\lambda)) \big)  \frac{dt}{\pi(1+t^2)} \\
& = \frac1{\pi} \sum_{s' : \, \mu_s - \mu_{s'} \in (0,\tfrac{4}{h^2})} {\rm Tr} \arctan \big( (\varepsilon x/2)^{-1} {\rm Im} \,T_{s'}(\lambda) \big) \\
& = \frac1{\pi} \sum_{s' : \, \mu_s - \mu_{s'} \in (0,\tfrac{4}{h^2})} {\rm Tr} \arctan \big( (\varepsilon x/2)^{-1} {\rm Im} \,T_{s'}(\mu_s) \big) + \mathcal O(1), \quad as \quad\lambda \searrow \mu_s.
\end{split}
\end{equation}
The asymptotic in \eqref{eq:ch4} follows by arguing as in the proof of 
Proposition \ref{prop:append} in the appendix. 
By putting together \eqref{eq:ch1}-\eqref{eq:ch4}, Proposition \ref{prop:dom04}, \eqref{eq:eqarctan} and Proposition \ref{prop:impos}, we get the claim.
\end{proof}

Applying Proposition \ref{prop:prel04} with $x=1$, one obtains immediately the following result.

\begin{Corollary}
\label{cor:prel}
Let $\nu_0 > 3$ in Assumption \ref{eq:hyppert}. Then, for any $\varepsilon \in (0,1)$,
 \begin{align*}
\frac1{\pi} & {\rm Tr} \arctan((1+\varepsilon)^{-1} {\rm Im} \, T_s(\lambda)) + \mathcal O(1) \\
 & \le \mp \xi(\lambda;H^\mp,H_Q) \\
& \le \frac1{\pi} {\rm Tr} \arctan((1-\varepsilon)^{-1} {\rm Im} \, T_s(\lambda)) + \mathcal O(1),
 \end{align*}
as $\lambda \searrow \mu_s$ and $\lambda \nearrow \frac4{h^2}+\mu_s$ for $0 \le s \le d$ fixed.
\end{Corollary}

Now, for $x>0$, $t \in \mathbb R$ and $\lambda \in (\mu_s,\frac4{h^2}+\mu_s)$, $0 \le s \le d$, if follows from Proposition \ref{prop:impos} that
\begin{equation*}
\label{eq:eqarctan+}
\mathscr N_\pm(x,t{\rm Im} \, T_s(\lambda)) 
= \mathscr N_\pm \Big( x,\frac{tb_s(\lambda) b_s(\lambda)^\ast}{\sqrt{\lambda-\mu_s} \sqrt{4/h^2+\mu_s-\lambda}} \Big) = \mathscr N_\pm(x,t\Omega_s(\lambda)).
\end{equation*}
This together with \eqref{eq:eqarctan} and Corollary \ref{cor:prel} gives Theorem \ref{theo2}.

\section{Appendix: proof of the asymptotics \eqref{as:lifs-kre0} and \eqref{as:lifs-kre4}}\label{sec9}

The aim of this section is to prove identities \eqref{as:lifs-kre0} and \eqref{as:lifs-kre4}. The operators $\Omega_s(\lambda)$ and $\Omega_{\bullet,s}(\lambda)$, $\bullet \in \{ 0,4 \}$ are respectively given by \eqref{eq:opeomega} and Remark \ref{rem:corr}.

\begin{Proposition}
\label{prop:append}
Let $V$ satisfy Assumption \ref{eq:hyppert} with $\nu_0 > 3$. Then, for any $x > 0$ and any fixed threshold $\mu_s \in \mathscr E_Q$, $0 \le s \le d$, we have
\begin{equation}
\label{as:lifs-kre0app}
{\rm Tr} \big( \arctan(x^{-1} \Omega_s(\lambda)) - \arctan(x^{-1} \Omega_{0,s}(\lambda)) \big) = \mathcal O(1), \qquad \lambda \searrow \mu_s,
\end{equation}
and
\begin{equation}
\label{as:lifs-kre4app}
{\rm Tr} \big( \arctan(x^{-1} \Omega_s(\lambda)) - \arctan(x^{-1} \Omega_{4,s}(\lambda)) \big) = \mathcal O(1), \qquad \lambda \nearrow \tfrac4{h^2}+\mu_s.
\end{equation}
\end{Proposition}

\begin{proof}
We only give the proof of \eqref{as:lifs-kre0app} since the one of \eqref{as:lifs-kre4app} follows in a similar way. 

Using \eqref{eq:eqarctan} and \eqref{eq:eqarctan+} one gets for $x > 0$ and $\lambda \in (\mu_s,\frac4{h^2}+\mu_s)$,
\begin{equation}
\label{eq:app1}
{\rm Tr} \big( \arctan(x^{-1} \Omega_s(\lambda)) \big) = {\rm Tr} \big( \arctan(x^{-1} \widetilde \Omega_s(\lambda)) \big),
\end{equation}
where according to 
Proposition \ref{prop:impos}, we have
\begin{equation}
\label{eq:imTapp-app}
\widetilde \Omega_s(\lambda) = {\rm Im} \,T_s(\lambda) = \frac{1}{\sqrt{\lambda-\mu_s} \sqrt{4/{h^2}+\mu_s-\lambda}} \mathscr V^{1/2} \big[ Y_s(\lambda)^\ast Y_s(\lambda) \otimes K \pi_s K^\ast \big] \mathscr V^{1/2},
\end{equation}
with the operator $Y_s(\lambda) : \ell^2(\bz_h) \to \bc^2$ defined by \eqref{eq:dim}. Similarly, one shows that
\begin{equation}
\label{eq:app1+}
{\rm Tr} \big( \arctan(x^{-1} \Omega_{0,s}(\lambda)) \big) = {\rm Tr} \big( \arctan(x^{-1} \widetilde\Omega_{0,s}(\lambda)) \big),
\end{equation}
where
\begin{equation}
\label{eq:imTapp-app+}
\widetilde\Omega_{0,s}(\lambda) = \frac{h}{2\sqrt{\lambda-\mu_s}} \mathscr V^{1/2} \big[ Y_{0}^\ast Y_{0} \otimes K \pi_s K^\ast \big] \mathscr V^{1/2},
\end{equation}
with the operator $Y_{0} : \ell^2(\bz_h) \to \bc^2$ defined by \eqref{eq:Y0nu}. It follows from the Lifshits-Krein trace formula \eqref{eq:fctdec} and identities \eqref{eq:app1}-\eqref{eq:imTapp-app+} that
\begin{equation}
\label{as:lifs-kre0app+1}
\begin{split}
\big \vert {\rm Tr} & \big( \arctan(x^{-1} \Omega_s(\lambda)) - \arctan(x^{-1} \Omega_{0,s}(\lambda)) \big) \big \vert \\
& \le \int_\br \vert \xi (s;x^{-1} \widetilde \Omega_s(\lambda),x^{-1} \widetilde\Omega_{0,s}(\lambda)) \vert ds \le \frac1x \Vert \widetilde \Omega_s(\lambda) - \widetilde\Omega_{0,s}(\lambda) \Vert_{\mathfrak S_1},
\end{split}
\end{equation}
(see \cite[Theorem 8.2.1]{yaf}). Thanks to \eqref{eq:imTapp-app} and \eqref{eq:imTapp-app+}, one has
\begin{equation}
\label{eq:app2+}
\begin{split}
\Vert & \widetilde \Omega_s(\lambda) - \widetilde\Omega_{0,s}(\lambda) \Vert_{\mathfrak S_1} \\
& \le \Vert  \mathscr V \Vert \Vert K \pi_s K^\ast \Vert_{\mathfrak S_1} \Big\Vert \frac{1}{\sqrt{\lambda-\mu_s} \sqrt{4/{h^2}+\mu_s-\lambda}} Y_s(\lambda)^\ast Y_s(\lambda) - \frac{h}{2\sqrt{\lambda-\mu_s}} Y_{0}^\ast Y_{0} \Big\Vert_{\mathfrak S_1}.
\end{split}
\end{equation}
Then, to conclude, it suffices to show that
\begin{equation}
\label{eq:app3+}
\Big\Vert \frac{1}{\sqrt{\lambda-\mu_s} \sqrt{4/{h^2}+\mu_s-\lambda}} Y_s(\lambda)^\ast Y_s(\lambda) - \frac{h}{2\sqrt{\lambda-\mu_s}} Y_{0}^\ast Y_{0} \Big\Vert_{\mathfrak S_1} = \mathcal O(\sqrt{\lambda-\mu_s}),\qquad \lambda \searrow \mu_s,
\end{equation}
as follows. Firstly, one can observe that for $\lambda \in (\mu_s,\frac4{h^2}+\mu_s)$,
\begin{equation}
\label{eq:app4+}
\begin{split}
& \Big\Vert \frac{1}{\sqrt{\lambda-\mu_s} \sqrt{4/{h^2}+\mu_s-\lambda}} Y_s(\lambda)^\ast Y_s(\lambda) - \frac{h}{2\sqrt{\lambda-\mu_s}} Y_{0}^\ast Y_{0} \Big\Vert_{\mathfrak S_1} \\
& \le \Big\vert \frac{1}{\sqrt{\lambda-\mu_s} \sqrt{4/{h^2}+\mu_s-\lambda}} - \frac{h}{2\sqrt{\lambda-\mu_s}} \Big\vert \Vert Y_s(\lambda)^\ast Y_s(\lambda) \Vert_{\mathfrak S_1} + \frac{h}{2\sqrt{\lambda-\mu_s}} \Vert Y_s(\lambda)^\ast Y_s(\lambda) - Y_{0}^\ast Y_{0} \Vert_{\mathfrak S_1}.
\end{split}
\end{equation}

a) Let us treat the first term of the r.h.s. of \eqref{eq:app4+}. It can be checked that $\frac{1}{\sqrt{\lambda-\mu_s} \sqrt{4/h^2+\mu_s-\lambda}} - \frac{h}{2\sqrt{\lambda-\mu_s}} = \mathcal O(\sqrt{\lambda-\mu_s})$ as $\lambda \searrow \mu_s$. Furthermore, it follows from \eqref{eq:impr2} that
\begin{equation*}
\Vert Y_s(\lambda)^\ast Y_s(\lambda) \Vert_{\mathfrak S_1} \le \sum_{n\in\bz} \cos^2[2n g_s(\lambda)] \bra n \ket^{-\nu_0} + \sum_{n\in\bz} \sin^2[2n g_s(\lambda)] \bra n \ket^{-\nu_0} = \sum_{n\in\bz} \bra n \ket^{-\nu_0}.
\end{equation*}
Consequently, one gets
\begin{equation}
\label{eq:app5+}
\Big\vert \frac{1}{\sqrt{\lambda-\mu_s} \sqrt{4/{h^2}+\mu_s-\lambda}} - \frac{h}{2\sqrt{\lambda}} \Big\vert \Vert Y_s(\lambda)^\ast Y_s(\lambda) \Vert_{\mathfrak S_1} = \mathcal O(\sqrt{\lambda-\mu_s}), \qquad \lambda \searrow \mu_s.
\end{equation}

b) Now, let us treat the second term of the r.h.s. of \eqref{eq:app4+}. A direct computation shows that for any $\phi \in \ell^2(\bz_h)$ and $n \in \bz$, 
$$
Y_{0}^\ast Y_{0} \phi(hn) = \bra (hn)h^{-1} \ket^{-\nu_0/2} \sum_{m\in\bz} \bra (hm)h^{-1} \ket^{-\nu_0/2} \phi(hm).
$$
This together with \eqref{eq:impr2} gives
\begin{equation*}
\begin{split}
& (Y_s(\lambda)^\ast Y_s(\lambda) - Y_{0,s}^\ast Y_{0,s})\phi(hn) \\
& = -2\sin^2[(hn)h^{-1} g_s(\lambda)] \bra (hn)h^{-1} \ket^{-\nu_0/2} \sum_{m\in\bz} \bra (hm)h^{-1} \ket^{-\nu_0/2} \phi(hm) \\
& - \bra (hn)h^{-1} \ket^{-\nu_0/2} \sum_{m\in\bz} 2\sin^2[(hm)h^{-1} g_s(\lambda)] \bra (hm)h^{-1} \ket^{-\nu_0/2} \phi(hm) \\
& + 2\sin^2[(hn)h^{-1} g_s(\lambda)] \bra (hn)h^{-1} \ket^{-\nu_0/2} \sum_{m\in\bz} 2\sin^2[(hm)h^{-1} g_s(\lambda)] \bra (hm)h^{-1} \ket^{-\nu_0/2} \phi(hm) \\
& + \sin[2(hn)h^{-1} g_s(\lambda)] \bra (hn)h^{-1} \ket^{-\nu_0/2} \sum_{m\in\bz} \sin[2(hm)h^{-1} g_s(\lambda)] \bra (hm)h^{-1} \ket^{-\nu_0/2} \phi(hm).
\end{split}
\end{equation*}
It follows that
\begin{equation*}
\begin{split}
& \Vert Y_s(\lambda)^\ast Y_s(\lambda) - Y_{0,s}^\ast Y_{0,s} \Vert_{\mathfrak S_1} \\
&  \hspace{0.5cm} \le 4 \Big( \sum_{n\in\bz} \sin^4[n g_s(\lambda)] \bra n \ket^{-\nu_0} \Big)^{1/2} \Big( \sum_{n\in\bz}  \bra n \ket^{-\nu_0} \Big)^{1/2} \\
& \hspace{2cm} + 4 \sum_{n\in\bz} \sin^4[n g_s(\lambda)] \bra n \ket^{-\nu_0} + \sum_{n\in\bz} \sin^2[2n g_s(\lambda)] \bra n \ket^{-\nu_0} \\
& \le \big( 4 \vert g_s(\lambda) \vert + 8 g_s^2(\lambda) \big) \sum_{n\in\bz} n^2 \bra n \ket^{-\nu_0} \underset{\lambda \searrow 0}{\sim} 2h \sqrt{\lambda-\mu_s} \sum_{n\in\bz} n^2 \bra n \ket^{-\nu_0}. 
\end{split}
\end{equation*}
Therefore,
\begin{equation}
\label{eq:app6+}
\frac{h}{2\sqrt{\lambda-\mu_s}} \Vert Y_s(\lambda)^\ast Y_s(\lambda) - Y_{0}^\ast Y_{0} \Vert_{\mathfrak S_1} = \mathcal O(1), \qquad \lambda \searrow \mu_s.
\end{equation}
One obtains immediately the claim by putting together \eqref{eq:app4+}, \eqref{eq:app5+} and \eqref{eq:app6+}.
\end{proof}

\begin{Remark}
If $\nu_0 > 5$, then we see from the proof that more precise estimates in Proposition \ref{prop:append} may be obtained so that
\begin{equation}
\label{as:lifs-kre0app+}
{\rm Tr} \big( \arctan(x^{-1} \Omega_s(\lambda)) - \arctan(x^{-1} \Omega_{0,s}(\lambda)) \big) = \mathcal O(\sqrt{\lambda-\mu_s}) = O(1), \qquad \lambda \searrow \mu_s,
\end{equation}
and
\begin{equation}
\label{as:lifs-kre4app+}
{\rm Tr} \big( \arctan(x^{-1} \Omega_s(\lambda)) - \arctan(x^{-1} \Omega_{4,s}(\lambda)) \big) = \mathcal O(\sqrt{4/{h^2}+\mu_s-\lambda}) = O(1), \quad \lambda \nearrow \tfrac4{h^2}+\mu_s.
\end{equation}
\end{Remark}

\medskip

{\bf Acknowledgements:} The authors deeply thank Vincent Bruneau for the careful reading and useful suggestions and comments. A. Taarabt has been supported by the Chilean grant Fondecyt 1230949, O. Bourget has been supported by the Chilean grant Fondecyt 1211576.



\def\cprime{$'$} \def\polhk#1{\setbox0=\hbox{#1}{\ooalign{\hidewidth
  \lower1.5ex\hbox{`}\hidewidth\crcr\unhbox0}}}
  \def\polhk#1{\setbox0=\hbox{#1}{\ooalign{\hidewidth
  \lower1.5ex\hbox{`}\hidewidth\crcr\unhbox0}}}
  \def\polhk#1{\setbox0=\hbox{#1}{\ooalign{\hidewidth
  \lower1.5ex\hbox{`}\hidewidth\crcr\unhbox0}}} \def\cprime{$'$}
  \def\cprime{$'$} \def\polhk#1{\setbox0=\hbox{#1}{\ooalign{\hidewidth
  \lower1.5ex\hbox{`}\hidewidth\crcr\unhbox0}}}
  \def\polhk#1{\setbox0=\hbox{#1}{\ooalign{\hidewidth
  \lower1.5ex\hbox{`}\hidewidth\crcr\unhbox0}}} \def\cprime{$'$}
  \def\cprime{$'$} \def\cprime{$'$}


\end{document}